\documentclass{gtart}

\newcommand{\pathtotrunk}{./}

\usepackage{amsmath,amssymb,amsfonts}
\usepackage{ifpdf}

\usepackage{ifthen}

\usepackage[plainpages=false,hypertexnames=false,pdfpagelabels]{hyperref}
\usepackage{breakurl}
\newcommand{\arxiv}[1]{\href{http://arxiv.org/abs/#1}{\tt arXiv:\nolinkurl{#1}}}

\ifpdf
\usepackage[pdftex,all,color]{xy}
\else
\usepackage[all,color]{xy}
\fi
\SelectTips{cm}{}%

\theoremstyle{plain}

\newtheorem{prop}{Proposition}[section]
\newtheorem{thm}[prop]{Theorem}
\newtheorem{lem}[prop]{Lemma}
\newtheorem{cor}[prop]{Corollary}
\newtheorem*{cor*}{Corollary}               
\theoremstyle{definition}
\newtheorem{defn}[prop]{Definition}         
\newtheorem*{defn*}{Definition}             

\newenvironment{rem}{\noindent\textsl{Remark.}}{}  
\numberwithin{equation}{section}

\newcounter{comment}
\newcommand{\noop}[1]{}

\ifpdf
\newcommand{\pathtodiagrams}{\pathtotrunk diagrams/pdf/}
\else
\newcommand{\pathtodiagrams}{\pathtotrunk diagrams/eps/}
\fi

\newcommand{\mathfig}[2]{{\hspace{-3pt}\begin{array}{c}%
  \raisebox{-2.5pt}{\includegraphics[width=#1\textwidth]{\pathtodiagrams #2}}%
\end{array}\hspace{-3pt}}}
\newcommand{\reflectmathfig}[2]{{\hspace{-3pt}\begin{array}{c}%
  \raisebox{-2.5pt}{\reflectbox{\includegraphics[width=#1\textwidth]{\pathtodiagrams #2}}}%
\end{array}\hspace{-3pt}}}
\newcommand{\rotatemathfig}[3]{{\hspace{-3pt}\begin{array}{c}%
  \raisebox{-2.5pt}{\rotatebox{#2}{\includegraphics[height=#1\textwidth]{\pathtodiagrams #3}}}%
\end{array}\hspace{-3pt}}}
\newcommand{\placefig}[2]{\includegraphics[width=#1\linewidth]{\pathtodiagrams #2}}

\def\bd{\partial}

\newcommand{\Integer}{\mathbb Z}
\newcommand{\Z}{\mathbb Z}
\newcommand{\hZ}{\frac{\Z}{2}}

\newcommand{\Id}{\boldsymbol{1}}
\renewcommand{\imath}{\mathfrak{i}}
\renewcommand{\jmath}{\mathfrak{j}}

\newcommand{\su}[1]{\mathfrak{su}_{#1}}

\newcommand{\juxta}{\bullet}
\def\refl#1{\overline{#1}}

\newcommand{\To}{\rightarrow}
\newcommand{\Into}{\hookrightarrow}

\newcommand{\Iso}{\cong}

\newcommand{\IsoTo}{\overset{\Iso}{\To}}

\newcommand{\clockwise}{\circlearrowright}
\newcommand{\counterclockwise}{\circlearrowleft}
\newcommand{\anticlockwise}{\circlearrowleft}

\newcommand{\Alt}{\operatorname{Alt}}

\newcommand{\htpy}{\simeq}

\newcommand{\cone}[3]{C\left(#1 \xrightarrow{#2} #3\right)}

\newcommand{\compose}{\circ}

\newcommand{\psmallmatrix}[1]{\left(\begin{smallmatrix} #1 \end{smallmatrix}\right)}

\newcommand{\directSum}{\oplus}
\newcommand{\DirectSum}{\bigoplus}
\newcommand{\tensor}{\otimes}

\newcommand{\directSumStack}[2]{{\begin{matrix}#1 \\ \DirectSum \\#2\end{matrix}}}
\newcommand{\directSumStackThree}[3]{{\begin{matrix}#1 \\ \DirectSum \\#2 \\ \DirectSum \\#3\end{matrix}}}

\newcommand{\OO}{{\mathcal{O} \To \mathcal{O}}}
\newcommand{\OP}{{\mathcal{O} \To \mathcal{P}}}
\newcommand{\PO}{{\mathcal{P} \To \mathcal{O}}}
\newcommand{\PP}{{\mathcal{P} \To \mathcal{P}}}

\newcommand{\tripletorus}{\mathfig{0.1}{cobordisms/triple_torus}}
\newcommand{\R}{\mathcal R}

\newcommand{\Kom}[1]{\text{Kom}\left(#1\right)}

\newcommand{\Mat}[1]{\text{Mat}\left(#1\right)}

\newcommand{\Hom}[3]{\operatorname{Hom}_{#1}\left(#2,#3\right)}

\newcommand{\Un}{{\mathbf{Un}}}
\newcommand{\Or}{{\mathbf{Or}}}
\newcommand{\Dis}{{\mathbf{Dis}}}
\newcommand{\Tang}{{\mathbf{Tang}}}
\newcommand{\Ab}{{\mathbf{Ab}}}
\newcommand{\PD}{{\mathbf{PD}}}
\newcommand{\OrPD}{{\mathbf{OrPD}}}
\newcommand{\DisPD}{{\mathbf{DisPD}}}
\newcommand{\OrTang}{{\mathbf{OrTang}}}
\newcommand{\DisAb}{{\mathbf{DisAb}}}
\newcommand{\DisTang}{{\mathbf{DisTang}}}
\newcommand{\UnAb}{{\mathbf{UnAb}}}
\newcommand{\AltAb}{{\mathbf{AltAb}}}
\newcommand{\C}{{\mathcal{C}}}
\newcommand{\MatC}{\Mat{\C}}
\newcommand{\KomC}{\Kom{\C}}
\newcommand{\MatUnAb}{\Mat{\UnAb}}
\newcommand{\MatDisAb}{\Mat{\DisAb}}
\newcommand{\KomUnAb}{\Kom{\UnAb}}
\newcommand{\KomDisAb}{\Kom{\DisAb}}

\newcommand{\grading}[1]{{\color{blue}\{#1\}}}

\newcommand{\parenthesise}[1]{\ifthenelse{\equal{#1}{}}{}{\left(#1\right)}}

\def\clap#1{\hbox to 0pt{\hss#1\hss}}

\usepackage{breakurl}

\usepackage{color}

\ifpdf
\usepackage[pdftex]{graphicx}
\else
\usepackage[dvips]{graphicx}
\fi

\def\llbracket{\left[\!\!\left[}
\def\rrbracket{\right]\!\!\right]}
\newcommand{\Kh}[1]{\llbracket#1\rrbracket}
\newcommand{\KhU}[1]{\Kh{#1}_U}
\newcommand{\KhD}[1]{\Kh{#1}_D}

\title[Fixing functoriality]{Fixing the functoriality of Khovanov homology.}

\author{David~Clark}
\address{
    Department of Mathematics, University of California, San Diego 92093-0112
} \email{dclark@math.ucsd.edu}
\urladdr{http://math.ucsd.edu/\~{}dclark}

\author{Scott~Morrison}
\address{
   Microsoft Station Q, University of California, Santa Barbara 93106-6105
} \email{scott@math.berkeley.edu} \urladdr{http://tqft.net/}

\author{Kevin~Walker}
\email{kevin@canyon23.net} \urladdr{http://canyon23.net/math/}

\date{
  First edition: December 9, 2006.
  This edition: \today.
}


\thanks{
  Electronic versions: \url{http://tqft.net/functoriality}
  and arXiv:math.GT/0701339.
}

\begin{document}

\begin{abstract}
We describe a modification of Khovanov homology \cite{MR1740682}, 
in the spirit of Bar-Natan \cite{MR2174270}, 
which makes the theory properly functorial with respect to link
cobordisms.

This requires introducing `disorientations' in the category of
smoothings and abstract cobordisms between them used in Bar-Natan's
definition. Disorientations have `seams' separating oppositely oriented
regions, coming with a preferred normal direction. The
seams satisfy certain relations (just as the underlying
cobordisms satisfy relations such as the neck cutting relation).

We construct explicit chain maps for the various Reidemeister moves,
then prove that the compositions of chain maps associated to each
side of each of Carter and Saito's movie moves
\cite{MR1238875,MR1445361} 
always agree. These calculations are greatly simplified by following
arguments due to Bar-Natan and Khovanov, which ensure that the two
compositions must agree, up to a sign. We set up this argument in
our context by proving a result about duality in Khovanov homology,
generalising previous results about mirror images of knots to a
`local' result about tangles.
Along the way, we reproduce Jacobsson's sign table \cite{MR2113903} 
for the original `unoriented theory', with a few disagreements.
\end{abstract}

\maketitle

\newpage%
\tableofcontents

\section{Introduction}
\label{sec:intro}
Khovanov homology \cite{MR1740682, MR1928174, MR2174270} is a ``categorified'' invariant:
it assigns to a link a graded module (or a complex of such)
rather than a ``scalar'' object such as a number or a polynomial.
Thus we expect not merely a module for each link, but also a functor which
assigns module isomorphisms to each isotopy between links.
(This isomorphism should depend only on the isotopy class of the isotopy.)
That is, given two links and a specific isotopy between them, we want an explicit
isomorphism between their Khovanov invariants, not merely the knowledge that
the Khovanov invariants are isomorphic.
Unfortunately, the original unoriented version of
Khovanov homology gives slightly less than this --- the isomorphisms assigned to
isotopies are well-defined only up to sign.

Unoriented Khovanov homology also gives more: the functor extends to surface cobordisms
in $B^3 \times I$ (but still with a sign ambiguity) \cite{MR2113903}.
More precisely, let ${\cal L}$ be the above category of oriented
links and (isotopy classes of) isotopies between them, and let
${\cal C}$ be the category whose morphisms are (isotopy classes of)
oriented surfaces properly embedded in $B^3\times I$.
If we associate to each isotopy between links the
track of the isotopy in $B^3\times I$, we get a functor ${\cal L} \Into {\cal C}$,
and the $Kh$ functor on ${\cal L}$ is the pull-back of an extended
$Kh$ functor on ${\cal C}$.
The extended $Kh$ also has a sign ambiguity.

The aim of this paper is to fix the above sign issues.

\begin{figure}[ht]
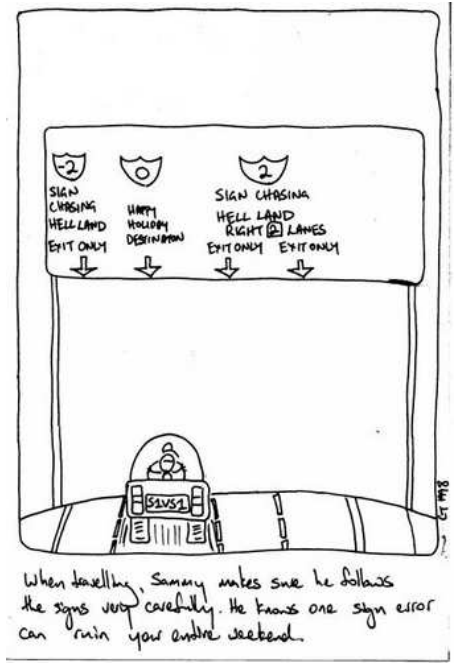

\begin{center}
\ifpdf
\placefig{0.55}{sammy/sammy}
\else
\placefig{0.45}{sammy/sammy_small}
\fi
\end{center}
\caption{Sammy the Graduate Student \cite{sammy}; used with permission.}
\end{figure}

For motivation, consider the `pre-categorified' situation.
Unoriented Khovanov homology is based on the unoriented
Kauffman bracket polynomial, with skein relation shown in Figure \ref{fig:kauffman_skein}
(with a further writhe correction, which introduces a dependence on the orientations of the link).
Closely related is the quantum $\su{2}$ polynomial,
which has a skein theory based on piecewise oriented (or ``disoriented") tangles,
as shown in Figure \ref{fig:disoriented_skein}
(see \cite{MR1117149}).
The two polynomials (and their associated TQFTs) differ only by a sprinkling of signs.
The Kauffman bracket has the advantage of simpler (unoriented) objects and
trivial Frobenius-Schur indicators,
while the quantum $\su{2}$ polynomial has the advantage of producing
positive-definite TQFTs (that is, TQFTs with nicer signs).

\begin{figure}[ht]
\begin{align*}
 \mathfig{0.065}{knot_pieces/unoriented_crossing} & =
 \mathfig{0.08}{knot_pieces/two_strand_unoriented} - q \rotatemathfig{0.08}{90}{knot_pieces/two_strand_unoriented} \\
 \mathfig{0.07}{knot_pieces/circle} & = q + q^{-1}
\end{align*}
\caption{A version of the Kauffman skein
relations.}\label{fig:kauffman_skein}
\end{figure}

\begin{figure}[ht]
\begin{align*}
 \mathfig{0.065}{knot_pieces/positive_crossing} & = q \mathfig{0.08}{knot_pieces/two_strand_identity} - q^2 \mathfig{0.08}{knot_pieces/disoriented_smoothing} \\
 \mathfig{0.065}{knot_pieces/negative_crossing} & = - q^{-2} \mathfig{0.08}{knot_pieces/disoriented_smoothing} + q^{-1} \mathfig{0.08}{knot_pieces/two_strand_identity} \\
 \mathfig{0.07}{knot_pieces/disoriented_strand1} & = - \mathfig{0.07}{knot_pieces/disoriented_strand2} \\
 \mathfig{0.07}{knot_pieces/two_disorientations} & = \mathfig{0.07}{knot_pieces/no_disorientations} \\
\end{align*}
\caption{The `disoriented' $\su{2}$ skein
relations.}\label{fig:disoriented_skein}
\end{figure}

Our strategy is to categorify the disoriented skein relation
of the quantum $\su{2}$ polynomial,
rather than the unoriented Kauffman skein relation.
We introduce the appropriate category of disoriented surface cobordisms,
and then imitate Bar-Natan's approach.
We find that disorientations also lead to nicer signs in the categorified setting:

\begin{thm}
\label{thm:isofunctoriality} There is a functor $Kh$ from the
category of oriented links in $S^3$ and (isotopy classes of)
isotopies between them to the category whose objects are (graded)
complexes of disoriented smoothings and abstract disoriented
cobordisms between smoothings (modulo local relations) and whose
morphisms are (graded) chain isomorphisms. Its graded Euler
characteristic, appropriately interpreted, gives the Jones
polynomial. It agrees with the original unoriented version of $Kh$,
modulo the sign ambiguity for isotopies in that theory.
\end{thm}

\begin{thm}
\label{thm:cobfunctoriality}
The above functor extends to the category of oriented links in $B^3$ and oriented
surface cobordisms (modulo isotopy) in $B^3 \times I$.
\end{thm}

We split the statement into two theorems because functoriality
with respect to isotopies of links would be expected of any
link invariant taking values in a category,
while functoriality with respect to surface cobordisms is a special feature
of Khovanov homology.

These results are also discussed in \cite{caprau}. The proofs of these statements given there are partially independent,
relying on our preparatory Lemmas at the beginning of \S \ref{ssec:movie_moves}. Further, the proofs in \cite{caprau},
just as with our proofs in the first arXiv version of this paper, omit checking some of the variations of certain movie moves. (See \S \ref{ssec:changelog} below and \S \ref{sssec:MM6} for details.)
\footnote{The results in \cite{caprau} are described as specialising to those here by setting a formal variable `$a$' equal to $0$, but this appears to be incorrect;
the paragraph after Lemma \ref{lem:ttorus} makes clear that our construction is agnostic to the value of the triple torus surface. The variable `$a$' in \cite{caprau} is simply some multiple of this surface.}

We actually get much more than a functor on cobordisms.
We can construct a 4-category (or, if you prefer, a 4-dimensional version of a planar algebra)
whose 3-morphisms are tangles in $B^3$ and whose 4-morphisms are elements of
appropriate Khovanov homology modules.
This 4-category enjoys the following duality or ``Frobenius reciprocity'' type
property:

\begin{thm}
\label{thm:frobrecip}
Given oriented tangles $P$, $Q$ and $R$, there is a duality isomorphism
between the spaces of chain maps up to homotopy
$$F : \Hom{Kh}{\Kh{P \juxta Q}}{\Kh{R}} \IsoTo \Hom{Kh}{\Kh{P}}{\Kh{R \juxta \refl{Q}}}.$$
The duality isomorphisms are coherent in the following sense (although this is not proved in the current version of this paper). 
To each such isomorphism we can associate an isotopy of links in $S^3$ ---
roughly speaking we slide $Q$ from the bottom of $S^3$ to the top.
Then two composable sequences of duality isomorphisms give the same result if
the associated isotopies in $S^3 \times I$ are isotopic.
\end{thm}


The paper is organized as follows.

Section \ref{sec:newcon} defines the invariant.
We introduce the appropriate category of disoriented cobordisms,
associate a chain complex based on this category to each oriented planar tangle diagram,
and associate a morphism of complexes to each Reidemeister and Morse move.

Section \ref{sec:moviemoves} verifies that our construction is well-defined.
We show that if two different sequences of Reidemeister and Morse moves are
related by movie moves, then the associated morphisms of chain complexes are equal.
Along the way, we prove the first part of the above duality result (Theorem \ref{thm:frobrecip}).

Section \ref{sec:oddsandends}, as its title suggests, contains miscellaneous results.
We show that setting $\omega = 1$ in our construction recovers the signs
from \cite{MR2113903}.
We show that modulo signs, our invariant agrees with the original
unoriented version.
We give an example calculation, showing that in the new construction, the cobordisms which `attach a handle to a strand' on either side
of a crossing give homotopic chain maps, whereas the old construction gave maps homotopic only with a sign.
Finally, we discuss the possibility of extending the invariant from
oriented tangles to disoriented tangles.

\subsection*{Acknowledgements}
David Clark would like to thank Justin Roberts for his encouragement and countless useful discussions, and Magnus Jacobsson for some helpful correspondence.

Scott Morrison would like to thank Dror Bar-Natan, for many useful discussions about Khovanov homology and his local cobordism model, and in particular
for sharing the idea that surfaces with piecewise orientations and some sort of seams might be useful in Khovanov homology. He'd also like to thank
Noah Snyder of UC Berkeley for an interesting discussion regarding the isomorphism between the usual Khovanov invariant of a knot, and the variation
defined here.

Kevin Walker thanks the NSF for support in the form of a Focused Research Group grant.
He also thanks Paul Melvin, Rob Kirby and Mike Freedman for helpful conversations.

We'd like to thank Chris Tuffley for allowing us to use his `Sammy the Graduate Student' comic \cite{sammy}, Scott Carter and Masahico Saito for allowing us to
reuse some of their diagrams from \cite{MR1238875}, and to offer our heartfelt apologies to Dror Bar-Natan, from whom we actually stole the movie move diagrams, callously failing to mention his addition of the `film-strip' edges to the same.

\subsection*{Changelog}
\label{ssec:changelog}
You're reading the `\textbf{v2}' version of this paper, which is available on the arXiv at \url{http://arxiv.org/abs/math.GT/0701339v2}, or
possibly a subsequent version. The following lists public versions of this work, and describes the differences between them.

\begin{itemize}
\item May 6 2007. Talk at Knots in Washington, slides available at \url{http://tqft.net/kiw}.
\item September 7 2007. Talk at Categorification in Uppsala, slides available at \url{http://tqft.net/uppsala}.
\item November 17 2007. Talk at Columbia Gauge Theory seminar, slides available at \url{http://tqft.net/columbia}.
\item[v0] December 9 2006. First public version, distributed by email and at \url{http://tqft.net/functoriality}.
\item[v1] January 12 2007. First arXiv version of the paper, available at \url{http://arxiv.org/abs/math.GT/0701339v1}.
\begin{itemize}
\item Corrected mistake in calculation of R1 chain maps; images switched.
\item Removed second section on duality. This section may appear separately later, as part of a paper on functoriality in $S^3$.
\end{itemize}
\item[v2] January ?? 2008. Second arXiv version of the paper, available at \url{http://arxiv.org/abs/math.GT/0701339v2}.
\begin{itemize}
\item Added David Clark as a coauthor.
\item More carefully described all 8 variations of the R3 move, along with the inverse maps and mirror image maps, in \S \ref{sec:R3-maps}.
\item Fixing incorrect cobordism diagrams for the $R3_{hml}$ move. 
\item Dealt correctly with all variations of MM6, in \S \ref{sssec:MM6}, including the `interleaved' variations, which had not been noticed in \textbf{v1}, or
in Caprau's subsequent paper on the functoriality of Khovanov homology, \cite{caprau}.
\item Removed the argument claiming to deal correctly with all 48 variations of MM10 directly; it was incorrect. The redundancy argument is still valid,
however, and there's now an illustration of the relevant 3-cell.
\item Included example of `sliding a handle past a crossing', in \S \ref{ssec:handle-sliding}
\end{itemize}
\end{itemize}

\section{The new construction}
\label{sec:newcon}
\subsection{Disorientations}
\label{ssec:disorientations}%

In this paper we follow the Bar-Natan approach of defining Khovanov
homology in terms of surface cobordism categories
--- categories whose objects are (possibly crossingless) tangles
in $B^3$ and whose morphisms are
surface cobordisms between tangles.
We'll deal with three sorts of tangles and surfaces:
unoriented (and possibly non-orientable), oriented, and disoriented.
We assume reader is familiar with the former two categories.

A disoriented 1- or 2-manifold is a piecewise oriented manifold where each
component of the interface between differently oriented domains is equipped
with a preferred normal direction.
In figures, we indicate this normal direction with a fringe pointing in the preferred direction. 
We'll call the interface between differently oriented domains a
disorientation seam.

\begin{figure}[ht]
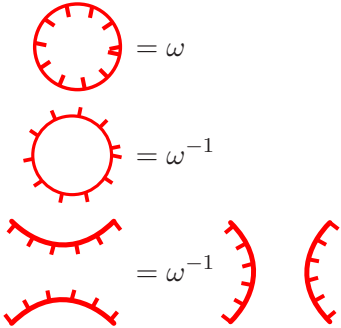

\begin{align*}
 \mathfig{0.09}{disorientation_relations/inward_circle} & = \omega \\
 \mathfig{0.1}{disorientation_relations/outward_circle} & = \omega^{-1}\\
 \mathfig{0.12}{disorientation_relations/parallel_fringes_together} & = \omega^{-1} \mathfig{0.12}{disorientation_relations/parallel_fringes_apart}
\end{align*}
\caption{Disorientation relations.}
\label{fig:disorientation_relations}
\end{figure}

We almost always (and usually without comment) consider disoriented
surfaces modulo the local fringe relations illustrated in Figure
\ref{fig:disorientation_relations}. If $\omega$ is a primitive
fourth root of unity ($\omega^2 = -1$), we will see below that we
get a version of Khovanov homology that satisfies functoriality. If
$\omega = 1$, then we reproduce the original unoriented version of
Khovanov homology, simply because the disorientations become
irrelevant. (We keep track of factors of $\omega$ explicitly, rather
than just writing $\omega = i$ everywhere, so that we can do
calculations in both the old and the new setup in parallel.)

\subsection{Cobordism categories}
\label{ssec:cobcats}%

The main goal of this paper is to construct a functor from $\OrTang$, the category
of oriented tangles and oriented cobordisms in $B^4$, to $\KomDisAb$,
a category of chain complexes based on abstract disoriented cobordisms between disoriented
crossingless planar diagrams.
Along the way we'll meet several other variant cobordism categories.
In this subsection we introduce the various categories we'll need.
The categories will be given compound names like $\OrTang$, $\KomDisAb$ and $\KomUnAb$;
we'll start by explaining the meanings of the components of the names.

The manifolds in the categories (1-manifolds for objects, 2-manifolds for morphisms)
can be unoriented, oriented or disoriented, which we denote by
$\Un$, $\Or$ and $\Dis$.
In all cases, we think of the objects as 1-manifolds embedded in $B^2\times I = B^3$,
with specified endpoints along the circle $\bd B^2 \times \{\frac{1}{2}\} \subset \bd B^3$.

We now introduce three categories of tangles. The first one, $\Tang$, is the one of real interest; it denotes the category whose
objects are arbitrary tangles in $B^3$ and whose morphisms
are isotopy classes of surface cobordisms embedded in $B^3\times I = B^4$.

The second, $\PD$, should be thought of as a `combinatorial model'
of $\Tang$. The objects of $\PD$ are tangles in $B^3$ which are in
general position with respect to the projection $p_z : B^3 \Iso B^2
\times I \To B^2$. The morphisms of the category can be
described by generators and relations. The generators are
\begin{itemize}%
\item Isotopies through tangles in general position.
\item Morse moves; birth or death of a circle, or a saddle move.
\item Reidemeister moves.
\end{itemize}
One should think of these generators as those isotopies which have
at most one `singular time slice'; that is, one moment at which the
projection of the link to $B^2$ is not generic, and the only the
simplest types of singularity are allowed to occur. These simplest
singularities are, of course, simply the Morse and Reidemeister
moves.

\begin{figure}[ht]
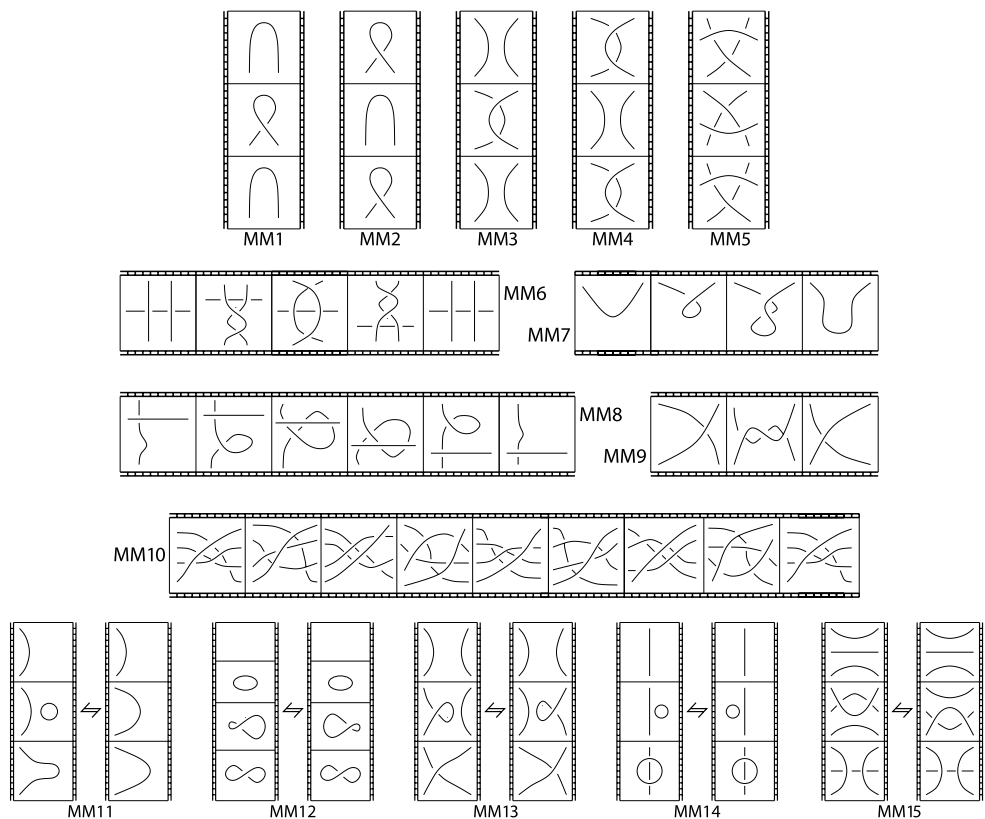

$$\mathfig{0.55}{movie_moves/MM1-5}$$
$$\mathfig{0.77}{movie_moves/MM6-10}$$
$$\mathfig{0.98}{movie_moves/MM11-15}$$
\caption{Carter and Saito's unoriented movies
moves.}\label{fig:unoriented_movie_moves}
\end{figure}

The first relation we impose is a boring one; composing an `isotopy
through general position tangles' with any other morphism simply
gives a morphism of the same type, given by gluing the isotopies
together. We then impose more relations, the movie moves of Carter
and Saito \cite{MR1238875, MR1445361} (see also Roseman
\cite{MR1634466}). The unoriented versions of these moves are shown
in Figure \ref{fig:unoriented_movie_moves} (thanks to Carter and
Saito for originally drawing these diagrams!), using the numbering
scheme introduced by Bar-Natan in \cite{MR2174270}. Note that we
also need to consider variations involving mirror images and/or
crossing changes.

They prove a theorem to the effect that two unoriented cobordisms
between unoriented tangles represented by compositions of Morse and
Reidemeister moves are isotopic if and only if those compositions
are related by a sequence of movie moves. To describe the relations we impose in $\OrPD$, we need the oriented
version of this, which, by much the same argument as they gave,
requires a separate version of each unoriented movie move for each
possible orientation of the strands (subject to some constraints;
movies involving saddles must have strands oriented appropriately so
the saddles are valid morphisms).

Finally, note that in $\DisPD$ there are both additional
Reidemeister moves (sliding a disorientation through a crossing) and
additional movie moves, involving this new Reidemeister move. As, in
this version of the paper, we're not discussing the extension of
Khovanov homology for $\DisTang$, we'll omit most of the details of
this, except what appears in \S \ref{ssec:confusions}.

Actually, we need to add a little more data to the objects in $\PD$;
a specified ordering on the crossings.
(The chain complexes we eventually assign to diagrams will
vary in boring but important ways
according to the ordering of the crossings.)
In addition to the morphisms
described above (Reidemeister and Morse moves), we need to add
`reordering morphisms', which are all isomorphisms. Further, we need
to modify our notion of the Reidemeister moves so that the source
and target tangles have (arbitrarily) ordered crossings -- but all
such different Reidemeister moves differ simply by pre- or
post-composition with reordering isomorphisms.

Finally, $\Ab$ denotes a category whose objects are
tangles without any crossings (think of them as embedded in $B^2\times \{\frac{1}{2}\} \subset B^3$).
The morphisms are abstract surfaces ({\it not} embedded in $B^4$),
modulo relations given below.
The $\Ab$ categories also have linearized morphisms spaces: morphisms
are linear combinations of cobordisms sharing the same range and domain, with coefficients in some ring containing $\frac{1}{2}$.
The relations we impose in $\Ab$ are
\begin{align}
\label{eq:cobordism_relations}
 \mathfig{0.075}{cobordisms/sphere} & = 0
 \qquad \qquad \mathfig{0.1}{cobordisms/torus} = 2 \\
\notag
 \mathfig{0.2}{cobordism_relations/cylinder} & = \frac{1}{2}
     \left(\mathfig{0.2}{cobordism_relations/neck_cutting_left} +
     \mathfig{0.2}{cobordism_relations/neck_cutting_right}\right)
\end{align}
Note that for $\DisAb$ these relations are imposed
away from the disorientation seams.
The last relation above is called the neck cutting relation.
In $\DisAb$ we of course also impose the fringe relations (Figure \ref{fig:disorientation_relations}, earlier).
We will see below that in $\DisAb$ it is unnecessary to set the 2-sphere equal to zero:
it follows from the fringe relations that
any connected, closed, orientable, disoriented surface whose Euler characteristic
is not a multiple of 4 is equivalent to zero.

The next lemma addresses the consistency of the neck cutting and fringe relations.

\begin{lem}
\label{lem:neckplusfringe}
Let $Y$ be a connected surface in $\DisAb$.
Then
\begin{enumerate}
\item If $Y$ has non-empty boundary, then $Y$ is not equivalent to zero.
\item If $Y$ is closed, orientable, and has odd genus, then $Y$ is not equivalent to zero.
If $Y'$ is a different disorientation on the same underlying surface, then
$Y$ and $Y'$ are equal up to a power of $\omega$.
\item If $Y$ is closed, orientable, and has even genus, then $Y$ is equivalent to zero.
\item There are no closed, nonorientable, disoriented surfaces.
\end{enumerate}
\end{lem}

\def\bdy{\partial}
\begin{proof}
First we consider the consistency of the fringe relations by themselves (no neck cutting).

Let $Y$ and $Y'$ be two disorientations on the same underlying surface $\Sigma$.
The disorientations seams of $Y$ and $Y'$ are properly embedded codimension 1
submanifolds of $\Sigma$ with oriented normal bundle (the orientation comes
from the direction of the fringe), and these determine cocycles
$a, a' \in C^1(\Sigma, \Z)$ which restrict to the same cocycle on $\bdy \Sigma$.
The fringe moves change these cocycles by coboundaries, so $Y$ and $Y'$
are related by a sequence of fringe moves if and only if
$a-a'$ is cohomologous to zero in $H^1(\Sigma, \bdy\Sigma; \Z)$.

Assume now that $Y$ and $Y'$ are related by two different sequences of fringe moves.
To each sequence we can associate a transversely oriented properly embedded surface in
$\Sigma\times I$ (consider the ``track" of the sequence of fringe moves), and so
we get cocycles $c_1, c_2 \in C^1(\Sigma\times I; \Z)$ which restrict to $a-a'$
on $\bdy(\Sigma\times I)$.

If $c_1-c_2$ is cohomologous to zero in $H^1(\Sigma\times I, \bdy(\Sigma\times I); \Z)$,
then the associated surfaces in $\Sigma\times I$ are related by a sequence of
elementary isotopies and
transversely oriented Morse moves.
Each such move changes the sequence of fringe moves, and one can verify that these
modifications do not change the factor of $\omega$ relating $Y$ and $Y'$.
(Figure \ref{fig:morse-moves-on-seams} shows the only nontrivial case.)

\begin{figure}[ht]
$$\xymatrix{
\mathfig{0.025}{disorientation_relations/morse_moves_1} \ar[r]^{\omega} & \mathfig{0.1}{disorientation_relations/morse_moves_2} \ar[r]^{\omega^{-1}} & \mathfig{0.1}{disorientation_relations/morse_moves_3}
}$$
\caption{The non-trivial sequence of Morse moves on seams.}\label{fig:morse-moves-on-seams}
\end{figure}
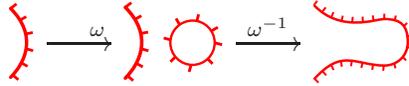

If $\Sigma$ has non-empty boundary, then $H^1(\Sigma\times I, \bdy(\Sigma\times I); \Z) = 0$,
and so $c_1$ and $c_2$ are always cohomologous.
So in this case the fringe relations are consistent: any two sequences of fringe moves relating $Y$ and $Y'$ yield the same factor of $\omega$.

If $\Sigma$ is closed then $H^1(\Sigma\times I, \bdy(\Sigma\times I); \Z) \cong \Z$, with the
generator corresponding to a transversely oriented surface in $\Sigma\times I$ parallel to
the boundary.
Put another way, in going from $Y$ to $Y'$ the disorientation seams can sweep out
$\Sigma$ an integer number of times.
If this integer changes by $\pm 1$, one can check that the factor of $\omega$ changes
by $\omega^{\pm\chi(\Sigma)}$.
Thus if $\Sigma$ has even genus ($\chi(\Sigma)$ congruent to 2 mod 4), then the fringe
relations are inconsistent and any disorientation is equivalent to $-1$ times itself and
hence to zero.  This proves part (3) of the Lemma.

Next we need to check the consistency of the combined fringe and neck cutting relations.
This boils down to checking that applying the neck cutting relation
on either side of a closed disorientation seam yields the same answer (Figure \ref{fig:cutting_disorientation}).
This implies part (1) of the Lemma, and using the next Lemma implies part (2).

\begin{figure}
\begin{align*}
 \mathfig{0.2}{cobordism_relations/disoriented_cylinder_left}
    & = \frac{1}{2} \left(\phantom{\omega^{-1}} \mathfig{0.2}{cobordism_relations/neck_cutting_left_d1} + \phantom{\omega^2}   \mathfig{0.2}{cobordism_relations/neck_cutting_right_d1}\right) \\
    & = \frac{1}{2} \left(\omega^{-1}           \mathfig{0.2}{cobordism_relations/neck_cutting_left_d2} + \omega^{\phantom{2}} \mathfig{0.2}{cobordism_relations/neck_cutting_right}\right) \\
    & = \frac{1}{2} \left(\omega^{-2}           \mathfig{0.2}{cobordism_relations/neck_cutting_left_d3} + \omega^2             \mathfig{0.2}{cobordism_relations/neck_cutting_right_d2}\right) \\
    & = \frac{1}{2} \left(\omega^{-1}           \mathfig{0.2}{cobordism_relations/neck_cutting_left}    + \omega^{\phantom{2}} \mathfig{0.2}{cobordism_relations/neck_cutting_right_d3}\right) \\
    & = \frac{1}{2} \left(\phantom{\omega^{-1}} \mathfig{0.2}{cobordism_relations/neck_cutting_left_d4} + \phantom{\omega^2}   \mathfig{0.2}{cobordism_relations/neck_cutting_right_d4}\right) \\
    & = \mathfig{0.2}{cobordism_relations/disoriented_cylinder_right}
\end{align*}
\caption{Checking that neck cutting on either side of a closed
disorientation seam yields the same answer.}
\label{fig:cutting_disorientation}
\end{figure}

Finally, a closed nonorientable disoriented surface must contain at least
one seam that is an
orientation reversing closed curve.
It is impossible to assign consistent fringe directions to such a curve.
This proves part (4).
\end{proof}

\begin{lem}
\label{lem:ttorus} Let $Y$ be an orientable surface in $\Un$, $\Or$ or $\Dis$
(and if in $\Dis$, further assume that the signed number of disorientation
fringes around each boundary component of $Y$ is zero).
Then $Y$ is equivalent to a
$\hZ[\omega]$-linear combination of disjoint unions of oriented
disks, punctured tori and closed genus 3 surfaces.
\end{lem}

\begin{proof}
Repeatedly apply the neck cutting relation, starting with curves parallel to the
disorientation seams and curves parallel to the boundary.
(In the disoriented case, the assumption about the boundary implies that after
applying fringe relations, we may assume that there is a curve parallel to each
boundary component which does not cross any seams.)
\end{proof}

As explained in \cite{MR2174270, math.GT/0603347}, setting the genus three
surface to zero in $\UnAb$
leads to the original version of Khovanov homology, while setting it to
a nonzero complex number gives something isomorphic to Lee homology \cite{MR2173845}.
Although it
makes very little difference for this paper, we'd like to encourage
leaving this surface unevaluated, as described in \cite{math.GT/0603347}. This makes the morphism spaces
into $\hZ[\tripletorus]$ modules. For convenience, we'll abbreviate
$\hZ[\tripletorus]$ simply as $\R$; although for the purposes of the
rest of the paper you can take $\R$ to be any ring with $2$
invertible, if you prefer.

Further, in all of these categories, we allow objects to carry an
integer, thought of as a `formal grading shift', just as in
\cite{MR2174270}. We'll denote this grading shift by a power of $q$.
We grade all of the morphism spaces, so that for a cobordism $C$
with source object $q^{m_1} D_1$ and target object $q^{m_2} D_2$,
each with $k$ boundary points, $\operatorname{deg}(C) = \chi(C) -
k/2 + m_2 - m_1$. It is not hard to see that these degrees are
additive under both composition and planar operations (in fact,
$\chi(C) - k/2$ and $m_2 - m_1$ are each additive separately). The
local relations in Equation \ref{eq:cobordism_relations} are clearly
degree homogeneous, so our grading makes sense on the quotient.

Given any category $\C$ with linear morphism spaces (called
`pre-additive' in \cite{MR2174270}), we can form a category $\MatC$
whose objects are tuples of objects of $\C$ (written as formal
direct sums), and whose morphisms are matrices of morphisms of $\C$.
Composition is given by multiplying matrices.

As an example to illustrate the grading and matrix conventions, let
us recall the `delooping' isomorphism described in
\cite{math.GT/0606318}. This is an isomorphism in $\Mat{\UnAb}$
(there is an identical isomorphism in $\Mat{\DisAb}$) between
$\mathfig{0.05}{knot_pieces/circle}$ and $q \emptyset \directSum
q^{-1} \emptyset$, given by the matrices
$\psmallmatrix{\rotatemathfig{0.025}{90}{cobordisms/cap_bdy_left}
\\ \frac{1}{2} \mathfig{0.025}{cobordisms/handle_bdy_down}}$ and
$\psmallmatrix{ \frac{1}{2}
\mathfig{0.025}{cobordisms/handle_bdy_up} &
\rotatemathfig{0.025}{90}{cobordisms/cap_bdy_right}}$.  That these
matrices are inverses follows immediately from the relations in Equation
\ref{eq:cobordism_relations} (and a quick calculation that the double torus is zero, by neck cutting). Observe that all the matrix entries
here are degree $0$ morphisms, once the grading shifts on the source
and target objects have been taken into account.

We can also form the category $\KomC$, whose objects are chain
complexes built out of $\MatC$, and whose morphisms are degree 0
chain maps modulo chain homotopy.

So, reviewing the nomenclature introduced thus far, we have:
\begin{itemize}
\item $\OrTang$ --- objects are oriented tangles in $B^3$, and morphisms are
oriented surface cobordisms in $B^4$.
\item $\OrPD$ --- objects are oriented tangles in $B^3$, with generic projection in the $z$ direction, and
morphisms are formal compositions (movies) of oriented surface
cobordisms, each of which has at most one `singular' moment, modulo
movie moves.
\item  $\UnAb$ --- objects are crossingless unoriented tangles in $B^3$, and morphisms are
linear combinations of abstract unoriented cobordisms, modulo local relations.
\item  $\DisAb$ --- objects are crossingless disoriented tangles in $B^3$, and morphisms are
linear combinations of abstract disoriented cobordisms, modulo local relations.
\item  $\KomDisAb$ --- objects are complexes in $\MatDisAb$, and
morphisms are chain maps modulo chain homotopy.
\end{itemize}

\begin{lem}
\label{lem:pdtangeq}%
In either the oriented or unoriented context, the functor $(i
\compose f):\PD \to \Tang$, which first forgets the ordering data on a
planar diagram in $\PD$, then includes the diagram into $\Tang$,
(recall tangles in $\PD$ have generic projections, whereas tangles
in $\Tang$ need not) is a natural isomorphism of categories. (The
same result is true in the disoriented context too, but we don't
need that for now.)
\end{lem}

\begin{proof}
First, we dispense with the ordering data on objects in $\PD$:
consider for a moment $\PD^{\text{unordered}}$, the same category as
$\PD$, but without the ordering data on crossings. The forgetful
functor $f$ is an equivalence of categories; its inverse (up to
natural isomorphisms) can arbitrarily specify the crossing ordering, after we've noticed
that all possible orderings on a diagram are
isomorphic.

Next, we construct a functor $j$ which is the inverse of the inclusion $i$ (up to natural isomorphism)
of $\PD^{\text{unordered}}$ into $\Tang$.
For every tangle $T$, choose an isotopy $I_T$ to a general position tangle $j(T)$
(object of $\PD$).
For every cobordism $Y: T_1 \to T_2$, $I_{T_2} Y T_{T_1}^{-1}$
is an isotopy from $j(T_1)$ to $j(T_2)$.
Up to `second order' isotopy, we can assume that $I_{T_2} Y T_{T_1}^{-1}$ is composed
of a sequence of Reidemeister moves and Morse moves.
Define $j(Y)$ to be this sequence of moves.

To show that $j(Y)$ is well-defined, we must show that choosing a
different second order isotopy above changes the sequence of
Reidemeister and Morse moves by movie moves. This is one of the
fundamental properties of movie moves. (Note that we have different
versions of movie moves for $\Un$, $\Or$.)

To complete the proof, it is easy to show that $\{I_T\}$ comprise an invertible
natural transformation between $ij$ and the identity functor on $\Tang$, and that
$\{j(I_T^{-1})\}$ comprise an invertible
natural transformation between $ji$ and the identity functor on $\PD$.
\end{proof}

The cobordism categories we've described above actually split up
into disjoint smaller categories, indexed by the number (and
possibly orientations, when relevant) of boundary points appearing on the equator
of $B^3$. These categories fit together as a canopolis (as
introduced in \cite{MR2147420}), that is, a planar algebra
\cite{math.QA/9909027} of categories. If you're unfamiliar with
planar algebras or canopolises, we've included a brief summary in
Appendix \ref{ssec:planar-algebras}. The planar operations are in
all cases simply given by gluing, both for objects and morphisms.

It's worth pointing out how the planar operations interact with the
ordering of crossings in objects of $\PD$. The internal discs of a
spaghetti and meatball diagram (indexing an operation of the planar
algebra) come with an ordering. When we glue together objects of
$\PD$ inside of one of these diagrams, we simply concatenate the
orderings specified inside each object.

The `matrix category' construction defining $\MatUnAb$ and
$\MatDisAb$ has an obvious analogue for canopolises; the planar
operations distribute over direct sums.

Similarly, taking complexes over a category extends to a
parallel construction for taking complexes over a canopolis. In any
canopolis $\C$, we can form a new canopolis $\Kom{\C}$ whose objects
complexes in $\C$ and whose morphisms are chain maps (or chain maps
up to homotopy). To apply a planar operation to a suitable
collection of complexes in $\Kom{\C}$, we take the formal tensor
product of the complexes (i.e. form a multicomplex, sprinkle signs,
and collapse), then apply the specified planar operation to each
object and differential. See Appendix
\ref{ssec:appendix:derived_canopolis} for more details. Notice that
this planar operation on complexes in $\Kom{\C}$ depends on the
ordering of the internal discs through the way that signs appear
when we take the tensor product of complexes, even when the original
canopolis was `symmetric'.

One consequence of these observations is that invariance for a local
model of a movie move implies invariance for that movie move
embedded in any larger tangle.





\subsection{Disoriented Khovanov homology}
\label{ssec:disoriented_khovanov_homology}
Our goal is to construct a map of canopolises (that is, a functor
for each category, compatible with planar operations) $\OrTang \to
\KomDisAb$. We follow closely Bar-Natan's approach, except that we
replace his target category $\KomUnAb$ with $\KomDisAb$. We'll write $\Kh{T}$ to denote the complex in $\DisAb$ associated to a tangle $T$.

It follows from Lemma \ref{lem:pdtangeq} that if we want to construct a functorial
invariant of $\OrTang$ it suffices to construct a functorial
invariant of $\OrPD$, and to do this it in turn suffices to
\begin{enumerate}
\item Construct a complex for each planar tangle diagram (equipped with an ordering of the crossings).
\item Construct a map of complexes for each Reidemeister move, each Morse move and each crossing reordering map.
\item Check that the relations coming from each oriented movie move are satisfied.
\end{enumerate}

We'll do the first two steps in this subsection and verify the movie move relations in \S \ref{ssec:movie_moves}.

\subsubsection{The complex}
\label{sec:complex}
The objects of $\OrPD$ are generated via planar algebra operations
by positive and negative crossings. We define the functor on single
crossings as follows:
\begin{equation}
\label{eq:invariant}
\xymatrix@R-1mm{
 \mathfig{0.08}{knot_pieces/positive_crossing} \ar@{|->}[r] & \Bigg( \bullet \ar[r] &
    q \mathfig{0.08}{knot_pieces/two_strand_identity} \ar[r]^{\mathfig{0.06}{cobordisms/saddle_seam_up}} & q^2 \mathfig{0.08}{knot_pieces/disoriented_smoothing} \Bigg) \\
 \mathfig{0.08}{knot_pieces/negative_crossing} \ar@{|->}[r] & \Bigg(
    q^{-2} \mathfig{0.08}{knot_pieces/disoriented_smoothing} \ar[r]^{\mathfig{0.06}{cobordisms/saddle_seam_down}} & q^{-1} \mathfig{0.08}{knot_pieces/two_strand_identity} \ar[r] & \bullet \Bigg) &
}
\end{equation}
In both cases, disorientation marks point to the right, relative to
the overall direction of the crossing. (This is just an arbitrary
convention; they could be equally well face to the left.)

Observe that a positive crossing is supported in homological heights
0 and 1, while a negative crossing is supported in heights -1 and 0.
We denote the grading shifts on objects simply by a multiplicative
factor of some power of $q$.

Next we must define the functor on morphisms of $\OrPD$. The
morphisms are generated (again, via planar operations) by
Reidemeister moves, Morse moves and the crossing reordering map
which switches the ordering of a pair of crossings. Note that Morse
moves (the cup, the saddle and the cap) are already morphisms of
$\DisAb$, and hence also morphisms (between one term complexes) of
$\KomDisAb$, so defining the functor on Morse moves is trivial.

When switching the ordering of a pair of crossings in a tangle, we
associate a chain map which is simply $\pm1$ on every object in the
complex. Following the homological conventions described in \S
\ref{ssec:appendix:permuting}, this map is $-1$ on objects in
which both crossings have been resolved in the disoriented way, and
$+1$ otherwise.

In the following sections, in which we describe the chain maps
associated to Reidemeister moves, we'll restrict our attention to
one particular ordering of the crossings in the source and target
tangle. The chain maps associated to other moves with other
orderings are simply obtained by pre- and post-composition with the
reordering maps from the previous paragraph.

Specifying the chain maps for the various Reidemeister moves will occupy
the remainder of this subsection.
Each of these chain maps will be invertible up to chain homotopy,
so by the end of this subsection we will have established the following weak result:
If two planar tangle diagrams are isotopic, then the complexes we assign to them
are isomorphic up to chain homotopy.
Full functoriality will not be established until we have verified the
movie move relations in \S \ref{ssec:movie_moves}.

\subsubsection{The R1 chain maps}
\label{sec:R1-maps}
\begin{figure}
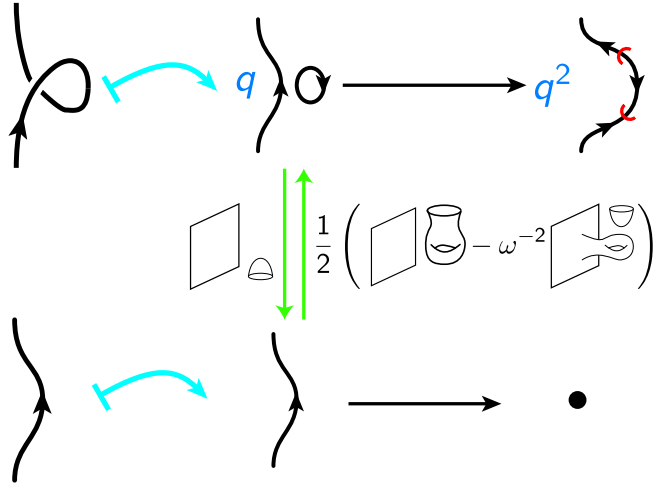

$$\mathfig{0.65}{reidemeister_maps/R1a_maps}$$
\caption{The R1a chain maps.}\label{fig:R1a_maps}
\end{figure}

\begin{figure}
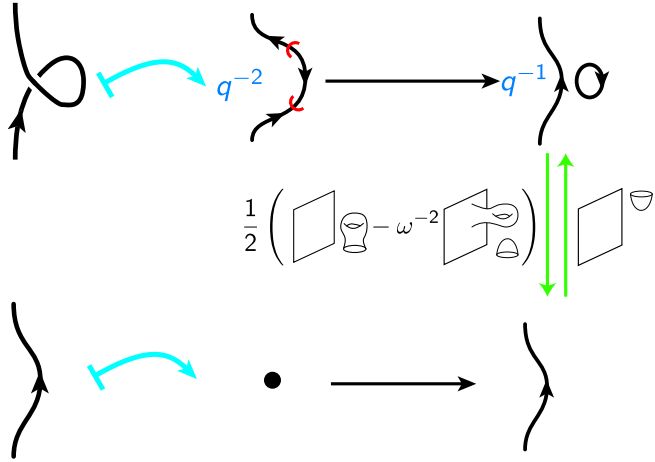

$$\mathfig{0.65}{reidemeister_maps/R1b_maps}$$
\caption{The R1b chain maps.}\label{fig:R1b_maps}
\end{figure}

The `twist' and `untwist' chain maps for the R1a and R1b moves are
shown in Figures \ref{fig:R1a_maps} and \ref{fig:R1b_maps}.  The
horizontal straight arrows are the differentials in the complex, and
the vertical (green) arrows show the chain map itself.


Being extra careful, we might want to distinguish two variations of
each of R1a and R1b, depending on whether the kink lies on the left
or the right side. However, the chain maps are just mirror images of
those shown here.

\subsubsection{The R2 chain maps}
\label{sssec:R2_maps}%
The Reidemeister 2 move comes in four
variations, which we'll call R2al, R2ar, R2b+ and R2b-.

\begin{align*}
 R2al & : \mathfig{0.2}{knot_pieces/R2al} \To \mathfig{0.15}{knot_pieces/two_strand_parallel} \\
 R2ar & : \mathfig{0.2}{knot_pieces/R2ar} \To \mathfig{0.15}{knot_pieces/two_strand_parallel} \\
 R2b+ & : \mathfig{0.2}{knot_pieces/R2b+} \To \mathfig{0.15}{knot_pieces/two_strand_antiparallel1} \\
 R2b- & : \mathfig{0.2}{knot_pieces/R2b-} \To \mathfig{0.15}{knot_pieces/two_strand_antiparallel2}
\end{align*}

Notice that we always chose to number the crossings so the negative
crossing comes first. This is, of course, an arbitrary choice, but
made so that the two R2a maps, and the two R2b maps, look as
similar to each other as possible.

\begin{figure}
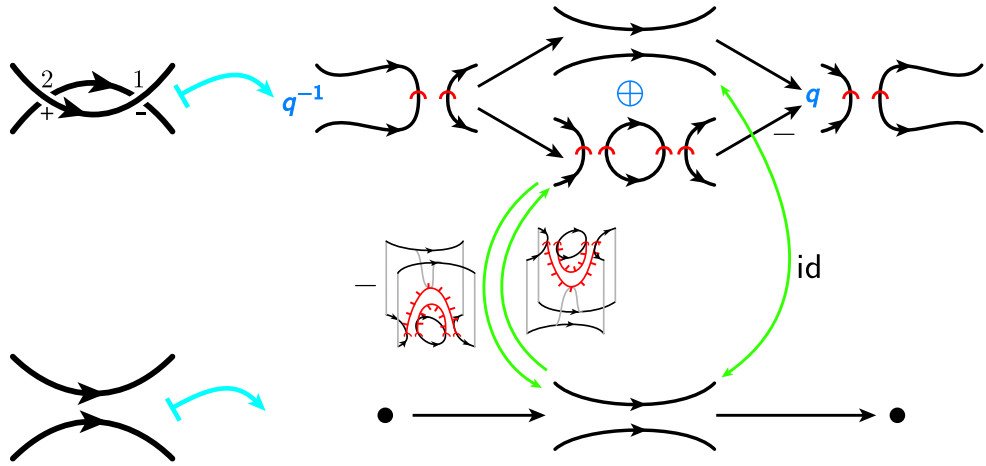

$$\mathfig{0.98}{reidemeister_maps/R2a_map}$$
\caption{The R2al chain map. (The R2ar chain map is
identical.)}%
\label{fig:R2a_maps}
\end{figure}

\begin{figure}
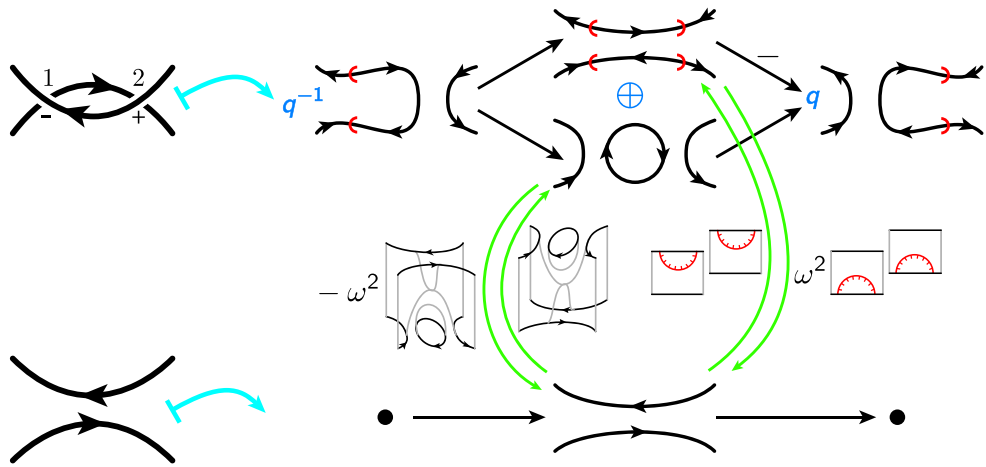

$$\mathfig{0.98}{reidemeister_maps/R2b_map}$$
\caption{The R2b- chain map. (The R2b+ chain map is the same,
but with all fringes reversed.)}%
\label{fig:R2b_maps}
\end{figure}

Explicit chain maps between the two sides of the Reidemeister R2al
and R2ar moves are shown in Figure \ref{fig:R2a_maps}, while maps
for the R2b- and R2b+ moves are shown in Figure
\ref{fig:R2b_maps}.


Calculations showing that these are indeed chain equivalences (and
showing how to discover them in the first place) have been relegated
to Appendix \ref{ssec:reidemeister_appendix}.

\subsubsection{The R3 chain maps}
\label{sec:R3-maps}
The work of this section is divided into three parts. First, we
explicitly describe a chain map for one variation of the $R3$ move, and
write down several properties of this chain map. Second, we state the
corresponding generalisations of these properties for the other seven
variations of the $R3$ move. Third, we describe an alternative chain map,
which is chain homotopic to the initial one, in each case.

We'll construct the chain maps for the first $R3$ move directly using the simplification
algorithm described by Bar-Natan in \cite{math.GT/0606318}; specifically, applying it to the complexes appearing on
either side of the Reidemeister move, we'll see that we obtain
(almost) exactly the same complexes. Composing the `simplifying' and
`unsimplifying' maps gives us the desired chain map. The result appears as Proposition \ref{prop:R3}.

We'll provide the chain maps for the other seven $R3$ moves less
explicitly, using the idea that all $R3$ moves are equivalent modulo $R2$
moves.

We'll state three lemmas (Lemmas \ref{lem:first-R3}, \ref{lem:second-R3},
\ref{lem:third-R3} for the first variation, and Lemmas
\ref{lem:first-R3-general}, \ref{lem:second-R3-general},
\ref{lem:third-R3-general} for the other seven variations) capturing the
features of these maps relevant to later movie move calculations, but
postpone the proofs until \S\ref{ssec:R3_variations}.

Sadly, the `categorified Kauffman trick' first described by Bar-Natan
\cite{MR2174270} doesn't work in the disoriented category; the
disorientation marks get in the way of using the second Reidemeister
move. With `vertigos' (as wished for in \S \ref{ssec:confusions}), this
method should recover its utility and give easier proofs of the
statements we need about the seven variations, by giving an easy direct
construction of the chain map in each case.

\begin{prop}
\label{prop:R3}%
There's a homotopy equivalence between the complexes associated to
either side of the Reidemeister move
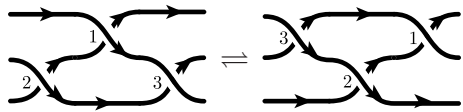
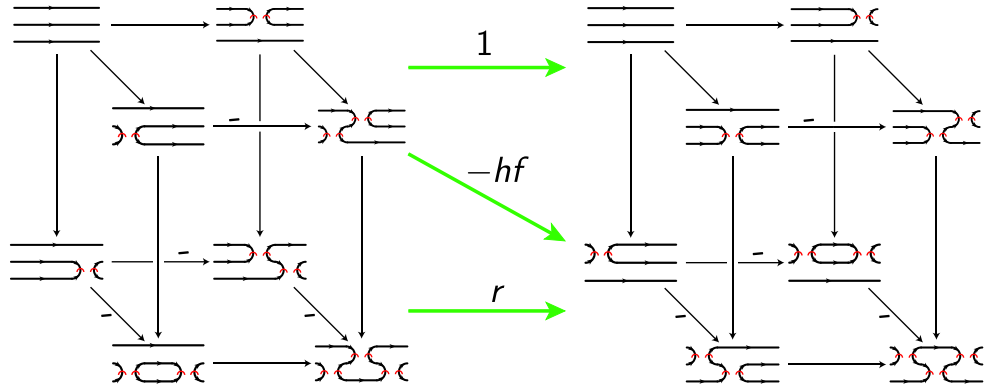
\begin{equation*}%
\label{eq:R3a}%
 \mathfig{0.2}{knot_pieces/R3al_horizontal} \rightleftharpoons \mathfig{0.2}{knot_pieces/R3ar_horizontal}
\end{equation*}%
given by
\begin{equation*}
\mathfig{0.98}{reidemeister_maps/R3_complex/cube_map}
\end{equation*}
The complex for each tangle is shown as a cube, with $8$ objects and $4$
homological levels. The two layers, top and bottom, correspond to the two different
resolutions of the highest crossing, labeled $3$. The chain map providing the
homotopy equivalence is the sum of the three (green) arrows each
connecting one layer of the left cube to a layer of the right cube.
The component maps are
\newcommand{\ps}[1]{\mathfig{0.1}{reidemeister_maps/R3_complex/R3a#1}}
\begin{align*}
-hf & = \ps{l_100} \xrightarrow{\Id} \ps{l_100} + \omega^2 \ps{l_110} \xrightarrow{\mathfig{0.075}{reidemeister_maps/R3_complex/c5}} \ps{r_101} \\
r & :
 \begin{pmatrix}
  \ps{l_101} \\ \\ \ps{l_011}
 \end{pmatrix}
 \To
 \begin{pmatrix}
  \ps{r_101} \\ \\ \ps{r_011}
 \end{pmatrix} =  \begin{pmatrix}
    -\omega^2 \mathfig{0.075}{reidemeister_maps/R3_complex/c6} & \omega^2 \mathfig{0.075}{reidemeister_maps/R3_complex/ms} \\
    - \Id & \mathfig{0.075}{reidemeister_maps/R3_complex/c1}
  \end{pmatrix}.
\end{align*}
\end{prop}
\begin{rem}
The names `$-hf$' and `$r$' shouldn't make any sense, unless you
know about the categorified Kauffman trick, and perhaps read a
future paper about the extension of Khovanov homology to disoriented tangles! If you do know the categorified
Kauffman trick, we'd be considering the cones over the morphisms
resolving the crossings labeled $3$.
\end{rem}
\begin{proof}
See \S\ref{ssec:reidemeister_appendix}.
\end{proof}

We won't need to know much about the details of this chain map,
however; what little we do is encapsulated in the following three
lemmas.
\begin{lem}[needed for MM6 and 10]
\label{lem:first-R3}%
The map from the bottom layer of the initial cube to the top layer
of the final cube is zero.
\end{lem}
\begin{lem}[needed for MM6, 8 and 10]
\label{lem:second-R3}%
The top layer of the initial cube is mapped identically to the top
layer of the final cube.
\end{lem}
\begin{lem}[needed for MM6]
\label{lem:third-R3}%
The leftmost and rightmost objects in the bottom layer are sent to
zero. That is, the map from the bottom layer to the bottom layer kills the highest and lowest homological height pieces.
Further, there is a single entry of that map, in the middle homological height, which is a multiple of the identity, that multiple is $-1$, and
every other nonzero entry has a disc component attached to a circle in either the source or target object (or both).
\end{lem}

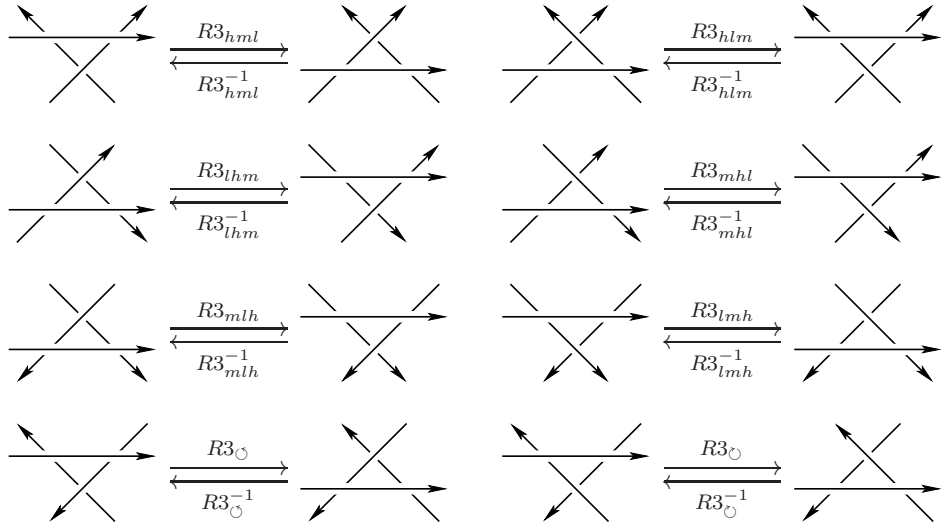
\begin{figure}[!ht]
\newcommand{\RIII}[1]{\xymatrix@C+20pt{%
  \mathfig{0.15}{reidemeister_maps/R3_variations/R3_#1_1} \ar@<0.5ex>[r]^{R3_{#1}} & \mathfig{0.15}{reidemeister_maps/R3_variations/R3_#1_2} \ar@<0.5ex>[l]^{R3_{#1}^{-1}}%
}}%
\begin{align*}
\RIII{hml} & & \RIII{hlm} & \\
\RIII{lhm} & & \RIII{mhl} & \\
\RIII{mlh} & & \RIII{lmh} & \\
\xymatrix@C+20pt{%
  \mathfig{0.15}{reidemeister_maps/R3_variations/R3_anti_1} \ar@<0.5ex>[r]^{R3_\circlearrowleft} & \mathfig{0.15}{reidemeister_maps/R3_variations/R3_anti_2} \ar@<0.5ex>[l]^{R3_\circlearrowleft^{-1}}%
} & &
\xymatrix@C+20pt{%
  \mathfig{0.15}{reidemeister_maps/R3_variations/R3_clock_1} \ar@<0.5ex>[r]^{R3_\circlearrowright} & \mathfig{0.15}{reidemeister_maps/R3_variations/R3_clock_2} \ar@<0.5ex>[l]^{R3_\circlearrowright^{-1}} %
}
 &
\end{align*}
\caption{The eight variations of the R3 move. These diagrams are taken from \cite{MR2186108}; they name the Reidemeister 3 moves differently, calling them $III_a$ through $III_h$, reading across the rows.}
\label{fig:R3-variations}
\end{figure}

Now there's not just one Reidemeister $3$ move; our version of Khovanov homology depends more explicitly on the orientations in
the original tangle than previous constructions, and as a consequence we need to do more work. There are eight R3 moves,
six `braidlike' and two `starlike'. We'll name the braidlike moves by walking counterclockwise around the
boundary, writing down the height of each outgoing strand. Thus in $R3_{hml}$
we see the `high' strand, the `middle' strand, then the `low
strand'. (We see the same sequence looking at the incoming strands.) The other braidlike moves are $R3_{hlm}$, $R3_{lhm}$, $R3_{mhl}$, $R3_{mlh}$ and $R3_{lmh}$.
There are then the two starlike R3 moves, which we'll call
$R3_\circlearrowleft$ and $R3_\circlearrowright$, depending on which
way we have to walk around the boundary in order to see the `outgoing low',
then `outgoing middle', then `outgoing high' strands. All eight Reidemeister 3 moves appear in Figure \ref{fig:R3-variations}.

When discussing these variations of the R3 move, we'll describe the left-hand of each pair of tangles as the `initial' tangle. In every case,
in the initial tangle the triangle lies to the left of the lowest strand, and in the final tangle it lies to the right. We also need to specify the
ordering of the crossings in these tangles. It turns out to be convenient to use a slightly unnatural ordering: in the initial tangle we number the
crossings as `middle', then `low', then `high', while in the final tangle we number them as `low', `middle', `high'. Notice this rule generalises the
ordering we used in describing $R3_{hml}$.

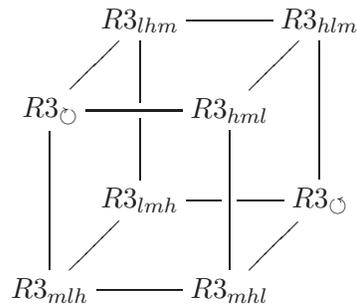
\begin{figure}[!ht]
\begin{equation*}
\xymatrix@!0@R+10pt@C+10pt{
                                                       & R3_{lhm} \ar@{-}[rr]\ar@{-}'[d][dd] &                                   & R3_{hlm} \ar@{-}[dd] \\
R3_\circlearrowright \ar@{-}[ur]\ar@{-}[rr]\ar@{-}[dd] &                                     & R3_{hml} \ar@{-}[ur]\ar@{-}[dd] \\
                                                       & R3_{lmh} \ar@{-}'[r][rr]            &                                   & R3_\circlearrowleft \\
R3_{mlh} \ar@{-}[rr]\ar@{-}[ur]                        &                                     & R3_{mhl} \ar@{-}[ur]
}
\end{equation*}
\caption{The cube of R3 moves.}
\label{fig:cube-of-R3s}
\end{figure}

The R3 moves fit together into a cube, shown in Figure \ref{fig:cube-of-R3s}.
The edges of this cube indicate pairs of $R3$ moves which are `related by R2 moves'. That is, for each edge there's a commutative diagram in the
category of tangles and tangle cobordisms. Here's one of the edges, connecting $R3_\circlearrowleft$ and $R3_{hlm}^{-1}$:
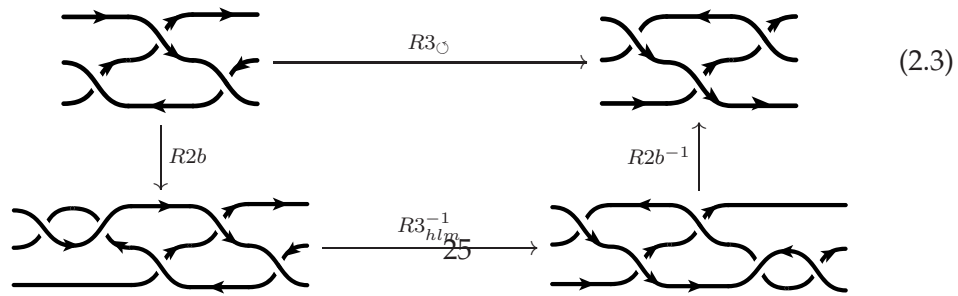
\begin{equation}
\label{eq:edge-example}
 \xymatrix@C+20mm{
    \mathfig{0.2}{knot_pieces/starlike_R3_l} \ar[d]^{R2b} \ar[r]^{R3_\circlearrowleft} & \mathfig{0.2}{knot_pieces/starlike_R3_r} \\
    \mathfig{0.3}{knot_pieces/starlike_R3_with_R2_l} \ar[r]^{R3_{hlm}^{-1}}     & \mathfig{0.3}{knot_pieces/starlike_R3_with_R2_r} \ar[u]^{R2b^{-1}} \\
 }
\end{equation}

We've already specified a chain map for $R3_{hml}$, in Proposition \ref{prop:R3}, and we can now specify chain maps for each of the others.
To do this, we pick some spanning tree for the cube. We'll now inductively define the chain map for a Reidemeister 3 variation
in terms of the already defined chain map for another variation adjacent in the spanning tree.
We simply write down the composition of the other three chain maps appearing in the commutative square corresponding to Equation \eqref{eq:edge-example}
for the appropriate edge.

We're never going to explicitly write down all the $R3$ maps; it would be incredibly tedious.
Instead, we'll just write down some lemmas (Lemmas \ref{lem:second-R3-general}, \ref{lem:first-R3-general} and \ref{lem:third-R3-general},
generalising Lemmas \ref{lem:second-R3}, \ref{lem:first-R3} and \ref{lem:third-R3} respectively), which encapsulate the facts we need for the movie move calculations.
We'll prove these statements by showing how they `propagate' along the edges of the cube in Figure \ref{fig:cube-of-R3s}.

Finally, you might worry about the choice of spanning tree. However, the sequence of movies moves corresponding
to a face of the cube is a cobordism isotopic to the identity, so functoriality will eventually assure us the choice didn't matter.

In order to state our more general lemmas, we'll need to describe various parts of the complexes appearing on either side of the variations of the R3 moves.
Thinking of such a complex as a cube, as in Proposition \ref{prop:R3}, we'll consider it as split into two layers, corresponding to the two resolutions
of the `highest' crossing (that is, the crossing between the `high' and `middle' strands). While we could describe the two layers as the `oriented' layer and the `disoriented layer', there's something
more useful; we'll describe them as the `orthogonal' ($\mathcal{O}$) and `parallel' ($\mathcal{P}$) layers, as shown in Figure \ref{fig:OP-layers}, depending on whether the strands in the resolution of the highest crossing are orthogonal or parallel
to the third strand not involved at the resolved crossing.

\begin{figure}[!ht]
\begin{equation*}
\mathfig{0.12}{reidemeister_maps/R3_variations/R3_hml_1} = \cone{\left(\mathcal{O} = \mathfig{0.12}{reidemeister_maps/R3_variations/R3_hml_1_O}\right)}{s}{\left(\mathcal{P}=\mathfig{0.12}{reidemeister_maps/R3_variations/R3_hml_1_P}\right)}
\end{equation*}
\caption{The complex for each tangle appearing in an R3 variation can be divided into two layers, the `orthogonal' ($\mathcal{O}$) and `parallel' ($\mathcal{P}$) layers.}
\label{fig:OP-layers}
\end{figure}

Notice that the differentials in the cube between the two layers point either from the $\mathcal{O}$ layer to the $\mathcal{P}$ layer, or from the $\mathcal{P}$ layer to the $\mathcal{O}$ layer.
This depends on which $R3$ variations we're looking at, in particular on the sign of the highest crossing, and whether its oriented resolution is orthogonal or parallel to the third strand.
Notice that order of layers alternates between $\OP$ and $\PO$ as we step across any edge in the cube in Figure \ref{fig:cube-of-R3s}.

\begin{figure}[!ht]
\begin{tabular}{c|ccc}
$R3$ variation           & highest crossing & order of layers & orthogonal layer is \\
\hline
$R3_{hml}$               & $+$ & $\OP$ & oriented \\
$R3_{hlm}$               & $+$ & $\PO$ & disoriented \\
$R3_{lhm}$               & $+$ & $\OP$ & oriented \\
$R3_{mhl}$               & $-$ & $\PO$ & oriented \\
$R3_{mlh}$               & $-$ & $\OP$ & disoriented \\
$R3_{lmh}$               & $-$ & $\PO$ & oriented \\
$R3_{\circlearrowleft}$  & $-$ & $\OP$ & disoriented \\
$R3_{\circlearrowright}$ & $+$ & $\PO$ & disoriented
\end{tabular}
\caption{The variations of the R3 move.}
\label{fig:table-of-R3s}
\end{figure}

We can then write each chain map $R3_\star$ (where $\star$ is one of $hml$, $hlm$, $lhm$, $mhl$, $mlh$, $lmh$, $\circlearrowleft$ or $\circlearrowright$) as the
sum of four components, $R3_\star = R3_\star^\OO + R3_\star^\OP + R3_\star^\PO + R3_\star^\PP$, where $R3_\star^{a \To b}$ are the components from the $a$
layer to the $b$ layer.

\begin{lem}
\label{lem:first-R3-general}%
If the layers of $\R3_\star$ are arranged as $\OP$, then the map from the parallel layer to the orthogonal layer, $R3_\star^\PO$, is zero.
Otherwise, if the layers are arranged as $\PO$, then the map $R3_\star^\OP$ is zero. (That is, the diagonal map pointing backwards in homological height is always zero.)
\end{lem}

\begin{lem}
\label{lem:second-R3-general}%
The map between the orthogonal layers, $R3_\star^\OO$, is the identity chain map, when $\star = hml, lhm, mhl$ or $lmh$.
When $\star = hlm, mlh, \circlearrowleft$ or $\circlearrowright$, the maps $R3_\star^\OO$ are nonzero multiples of a certain standard chain map; forgetting
disorientation data and coefficients, this map is the identity chain map. The disorientation seams are the minimal ones compatible with the boundary disorientation marks.
The coefficients
are all either $-1$ or $\omega^2$, and are determined by the rule that the coefficient $\kappa_\star$ of the chain map in the lowest homological height is
given by
\begin{align*}
\kappa_\star & =
\begin{cases}
\omega^2 & \text{if $\star = hlm$ or $\clockwise$} \\
-1 & \text{if $\star = mlh$ or $\counterclockwise$}.
\end{cases}
\end{align*}
The fact that
these are chain maps then determines the other coefficients; in particular, on the highest homological height the coefficient is $-\omega^2 \kappa_\star$.
\end{lem}
\begin{rem}
This dichotomy distinguishes whether the orthogonal layer of the cube comes from an oriented or disoriented resolution of the highest crossing. These data
are displayed in the last column of Figure \ref{fig:table-of-R3s}.
\end{rem}

As an example, the map $R3_{hlm}^\OO$ is
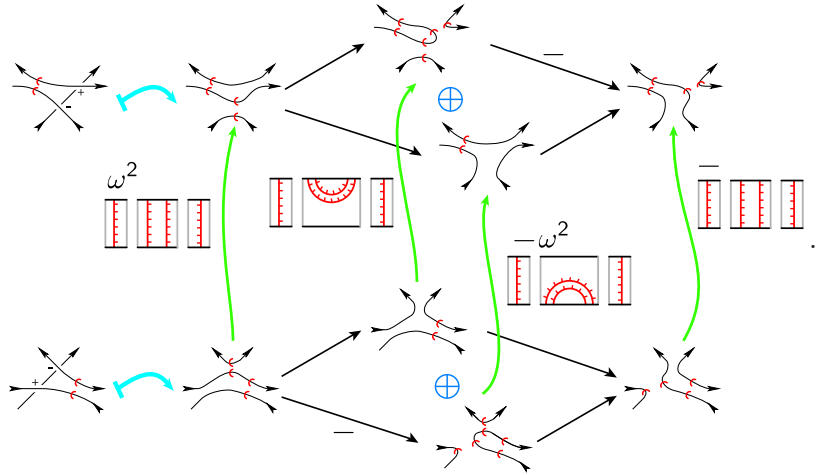
\begin{equation*}
\mathfig{0.8}{reidemeister_maps/R3_variations/O-O_chain_map_hlm}.
\end{equation*}
Notice here that the inverse map is obtained by taking the adjoint (reflection in the time direction) of each disoriented surface, and moving the
coefficient of $\omega^2$, but not the coefficient of $-1$, to the other component in homological height $0$.

As another example, the map $R3_{\circlearrowleft}^\OO$ is
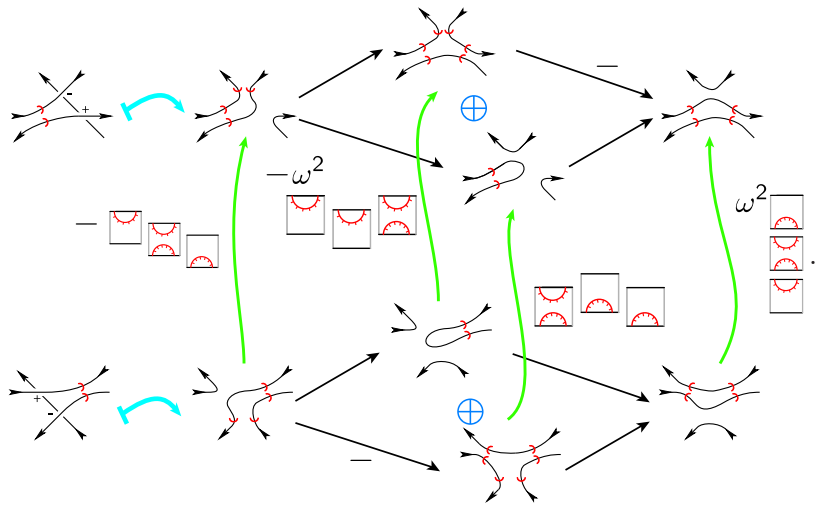
\begin{equation*}
\mathfig{0.8}{reidemeister_maps/R3_variations/O-O_chain_map_left}.
\end{equation*}
Again, the inverse map is obtained by taking the adjoint, and moving the coefficient of $\omega^2$ (but not the coefficient of $-1$), appearing at height $0$ over to the other map at that
height.

\begin{lem}
\label{lem:third-R3-general}%
The maps between the parallel layers, $R3_\star^\PP$, kill the highest
and lowest homological heights. Further, in the middle homological height
there are a pair of objects (one in the source complex, one in the target
complex) which have the same unoriented diagram, and the component of the
$R3_\star^\PP$ map between these is the unique disoriented surface with
minimal disorientation seams, and a coefficient of
\begin{align*}
p_\star & =
\begin{cases}
-1        & \text{if $\star = hml$ or $lmh$} \\
1         & \text{if $\star = hlm$ or $mlh$} \\
\omega^2  & \text{if $\star = lhm$ or $mhl$} \\
-\omega^2 & \text{if $\star = \clockwise$ or $\counterclockwise$.}
\end{cases}
\end{align*}
Every other entry of the map in the middle homological height is some multiple of a surface with a disc component attached to a circle in either the source or target object (or both).
\end{lem}

The proofs appear in \S\ref{ssec:R3_variations}.

At this point we can also give a description of the inverses of these chain maps.
\begin{cor}
\label{cor:R3-inverses}
Lemmas \ref{lem:first-R3-general}, \ref{lem:second-R3-general} and \ref{lem:third-R3-general} also hold without changes when describing the inverses of
the $R3$ chain maps.
\end{cor}
\begin{proof}
Consider the operation of rotating a tangle by $\pi$, and reversing all
orientations. Notice that this interchanges the source and target tangles
of each $R3$ variation.

Being a little more careful, and thinking about the source and target
tangles with their specified ordering of crossings, this operation
actually needs to be followed by switching the ordering of the low and
middle crossings.

Thus for each $R3$ variation, we produce a chain map pointing the
opposite direction, by rotating each component of the original chain map
by $\pi$, reversing all orientations and disorientations, and introducing
an extra sign in each component going between a pair of resolutions in
which for one or the other of the initial and final resolutions, but not
both, both the low and middle crossings have been resolved in the
disoriented way.

We now make two claims. Firstly, that this chain map really is the
inverse of the original map, and secondly, that this chain map is
correctly described by Lemmas \ref{lem:first-R3-general},
\ref{lem:second-R3-general} and \ref{lem:third-R3-general}.

First, we consider the $\OO$ parts of the map. It is readily seen
(trivial in the cases $hml, lhm, mhl$ or $lmh$, easy in the cases $hlm$
and $mlh$, and requiring an easy calculation involving disorientations in
the cases $\clockwise$ and $\anticlockwise$) that at the lowest
homological height, the composition of the original map and the candidate
inverse is the identity. This is enough to know that the candidate really
is the inverse.

Second, Lemma \ref{lem:first-R3-general} holds obviously, Lemma
\ref{lem:second-R3-general} holds because the signs introduced by reordering
occur at homological height $0$, so cannot affect the sign
$\kappa_\star$, and Lemma \ref{lem:third-R3-general} holds because the
reordering signs  occur at heights $\pm 1$, so cannot affect the sign
$p_\star$.

Notice that the inverses of the example $\OO$ maps given above agree with
the description here.
\end{proof}

The third task of this section is to describe an alternative chain
map for each Reidemeister 3 move. This alternative will be chain
homotopic to the one described above, but not identical.

The mirror image (in the direction perpendicular to the plane) of a tangle is simply the obvious topological operation.
At the level of the corresponding Khovanov complexes, this corresponds to
negating the homological height of each step
of the complex, and replacing each differential with its time
reverse, by switching source and target. That is, the mirror image of
a complex $\left(C^\bullet, d\right)$ is $\left(\overline{C}^\bullet,
\overline{d}\right)$, with $\overline{C}^i = C^{-i}$, and
$(\overline{d}_i : \overline{C}^i \To \overline{C}^{i+1}) =
(d_{-i-1} : C^{-i-1} \To C^{-i})^*$, where the $*$ here means time
reversal, or `adjoint'. By the mirror image of a chain map
$f^\bullet$, we mean $\overline{f}^\bullet$, with $\overline{f}^i$
= $f^{-i}$; that is, exactly the same components, but each in negated
homological height.

We can think of the alternative chain
map in two different ways. First, and secretly, we think of it as
coming from performing the Kauffman trick on the lowest crossing,
rather than the highest crossing as above. Second, we can simply
think of it, and define it, as the mirror image of
one of the chain maps above. Actually, more precisely, we need to modify this mirror image in two ways.
First, in all cases, we must pre- and post-compose with crossing reordering maps,
to ensure that we start and finish at the same ordered tangles as the usual chain maps.
Second, only for the starlike R3 variations, we need to multiply the mirror image chain map by $-\omega^2$. (This will ensure that
the mirror image chain map really is homotopic to the usual one. Recall of course that in the disoriented theory, $-\omega^2 = 1$!)

Notice that taking mirror image exchanges pairs of R3 move variations, switching the labels `h' and `l', and interchanging $\circlearrowleft$ and $\circlearrowright$.
Thus $R3_{hml}$ and $R3_{lmh}$, which are antipodal in the cube of R3 variations in Figure \ref{fig:cube-of-R3s}, are exchanged, as are $R3_\circlearrowleft$ and $R3_\circlearrowright$.
The other pairs are $R3_{lhm}$ and $R3_{hlm}^{-1}$, and $R3_{mlh}$ and $R3_{mhl}^{-1}$, which are each adjacent in the cube.

To distinguish the chain maps defined in this way from the ones described above, we'll write a bar over the top. Thus $\overline{R3}_{hml}$ is defined
by taking the chain map for $R3_{lmh}$, and applying the mirror image operation described in the paragraph above, and reordering crossings in the source and
target tangles appropriately.

Passing to the mirror image move reverses the `order of layers' appearing in
Figure \ref{fig:table-of-R3s}. It's easy to see that the mirror image of a chain map for one vertex does not give the chain map for the opposite vertex
which has been described above. This is essentially because the lemmas above are written in terms of the orthogonal and parallel layers with
respect to the highest crossing, which are not preserved by mirror image. For example, look at the pair $R3_{hml}$ and $R3_{lmh}$, and in particular
the completely oriented resolution. The completely oriented resolution is in the orthogonal layer
for $R3_{hml}$, so the chain map above acts as the identity here, by Lemma \ref{lem:first-R3-general} (or indeed, the original special case Lemma \ref{lem:first-R3}).
However the completely oriented resolution is killed by the usual chain map for $R3_{lmh}$, being in the parallel layer, using Lemma \ref{lem:third-R3-general}.
Thus we see that the chain map for $R3_{hml}$ coming from the mirror image of the chain map for $R3_{lmh}$ is in fact different from the usual one.

On the other hand, these maps turn out to be homotopic to the usual maps, even though we have seen they are not equal on the nose.
The argument relies on two results which live more naturally later in the paper, namely Corollary \ref{cor:reidemeister-maps-unique}, appearing in the next section,
and Lemma \ref{lem:loopless} appearing in \S \ref{ssec:movie_moves}, so the reader may prefer to postpone deciphering this argument until having reached those statements!
Corollary \ref{cor:reidemeister-maps-unique}, appearing in the next section, assures us that the relevant space of chain maps,
up to homotopy, is $1$ dimensional. Thus we know that each mirror image map must be homotopic to some multiple of the usual map, and we only need to show
that multiple is always $1$. To do this, we look at a particular resolution, namely the unique resolution which is in an extreme homological height of the
$\mathcal{O}$ layer for both the usual map and the mirror image map. Lemma \ref{lem:second-R3-general} then describes how this resolution is mapped to the
corresponding resolution of the target tangle, and it suffices, by Lemma \ref{lem:loopless} to check that both the usual map and the mirror map act in the same way, without coefficients, on this resolution.
That check follows directly from Lemma \ref{lem:second-R3-general}, along with the relevant crossing reordering calculations. Recall also the coefficient of $-\omega^2$ which we smuggled into the definition of the mirror image maps for the starlike moves, precisely to allow the present result.

An important point we  need to make is that the three Lemmas \ref{lem:first-R3-general}, \ref{lem:second-R3-general} and \ref{lem:third-R3-general}
still apply to the mirror image maps, replacing as needed each reference to an R3 variation $R3_\star$ with $\overline{R3_{\star'}}$, where $\star'$ is the mirror variation,
and understanding `orthogonal' and `parallel' layers as referring to the layers
given by resolving the lowest, rather than the highest, crossing.

Note in particular in regard to Lemma \ref{lem:second-R3-general}, that while $R3_\star^\OO$ is the identity chain map, when $\star = hml, lhm, mhl$ or $lmh$,
when we look at $\overline{R3}_\star^\OO$, it is $\star = hml, lmh, mlh$ and $hlm$ that give the identity. This will be important in the discussion of movie move 6.

\section{Checking movie moves}
\label{sec:moviemoves}
\subsection{Duality, and dimensions of spaces of chain
maps}
\label{sec:duality}

Most nice (or at least, interesting to topologists) monoidal
categories have duals. There are many formulations of this; see for
example \cite{MR1686423} for `pivotal categories', etc. The category
$\C$ should have an involution $*$ on objects, called the dual, and
isomorphisms between hom-sets of the form
$$\Hom{\C}{U \tensor V}{W} \Iso \Hom{\C}{U}{W \tensor V^*}$$
(along with the three other obvious variations of this), satisfying
some axioms (corresponding diagrammatically to `straightening an S-bend'). 

There's no shortage of examples. Categories of diagrams up to
isotopy \cite{MR1113284} are generally tautologically equipped with
duals, given by $\pi$ rotations, and the natural isomorphisms
between hom-sets are just planar isotopies. Categories of
representations of quantum groups have duals, provided by the
antipode in the Hopf algebra structure of the quantum group.
Bimodules over a von Neumann algebra have duals; there the
isomorphism between hom-sets is called ``Frobenius reciprocity''
\cite{MR1424954}.

We'll prove a result along these lines here. To fit with the above
pattern, briefly consider the 2-category whose objects are
(oriented) points on a line, whose 1-morphisms are tangles between
these points, and whose 2-morphisms are chain maps up to homotopy
between the Khovanov complexes associated to the tangles. There's a
duality functor, at least at the level of 0- and 1-morphisms, given
by reflection. We'll prove that there are isomorphisms of the type
described above.

In our case there is more structure than in the above examples,
since we're actually in a 3- or 4-category rather than a 2-category.
(3-category if we're thinking in terms of tangle projections living in $B^2$;
4-category if we're thinking in terms of unprojected tangles living in $B^3$, with cobordisms in $B^4$.)
More specifically, we can glue tangles $P$ and $Q$ together anywhere
along their boundaries --- we're not limited to
tensoring on the right or tensoring on the left.
We'll denote any gluing of tangles $P$ and $Q$ by $P \juxta Q$ (or, equivalently,
by $Q\juxta P$).

\begin{prop}
\label{prop:duality}%
Given oriented tangles $P$, $Q$ and $R$, there is an isomorphism
between the spaces of chain maps up to homotopy
$$F : \Hom{Kh}{\Kh{P \juxta Q}}{\Kh{R}} \IsoTo \Hom{Kh}{\Kh{P}}{\Kh{R \juxta \refl{Q}}}.$$
($\refl{Q}$ denotes the reflection of $Q$.)

Diagrammatically, this statement claims that there's an isomorphism
between the spaces of chain maps we can fill inside the following
two cylinders.
$$\mathfig{0.175}{duality/cylinder1} \Iso
\mathfig{0.175}{duality/cylinder2}$$

These isomorphisms are natural in the sense that they are
compatible with pre-composition with a
morphism into $P$, and with post-composition with a morphism out of
$R$.

We can actually make a stronger statement, which includes grading
shifts. If the tangle $Q$ has $m$ boundary points attached to $P$ to
in $P \juxta Q$, and $n$ boundary points attached to $R$ in $R
\juxta \refl{Q}$, the isomorphism is in fact
 $$\Hom{Kh}{\Kh{P\juxta Q}}{\Kh{R}} \Iso \Hom{Kh}{\Kh{P}}{\Kh{R \juxta\refl{Q}}}\grading{\frac{m-n}{2}}.$$
\end{prop}
\begin{rem}
For now, we're just claiming that there is some isomorphism; in
particular, all we'll need for now is that the dimensions of the
morphisms spaces are the same.

In a future paper, we'll explain a coherence result for these
isomorphisms. Essentially this result is the difference between
`functoriality in $B^3$' and `functoriality in $S^3$'. There are pairs of
cobordisms in $B^3$ which are not isotopic in $B^3$, but become isotopic
in $S^3$. The coherence result for the maps described in the proposition
above requires us to show that such pairs give homotopic chain maps, and
this remains beyond the scope of the current paper.
\end{rem}
\begin{proof}
\newcommand{\rQ}{\mathfig{0.05}{duality/positive_crossing}}
\newcommand{\Q}{\mathfig{0.05}{duality/negative_crossing}}
\newcommand{\st}{\mathfig{0.0275}{duality/strand}}
\newcommand{\rst}{\mathfig{0.0275}{duality/other_strand}}

We'll prove the result for a short list of (very!) small tangles
$Q$, which easily imply the rest. Namely $Q = \Q, \rQ, \st$ and $\rst$, and the other oriented versions of these tangles.
We can then build the
isomorphism for an arbitrary $Q$ by composing isomorphisms for the
constituent pieces of the tangle $Q$.

We'll begin with $Q=\Q$, a negative crossing oriented to the right.
(The case for a positive crossing is exactly analogous.) Given a
chain map $f \in \Hom{Kh}{\Kh{P \juxta \Q}}{\Kh{R}}$, we'll produce
the chain map
 $$F(f) = (f \juxta \Id_{\rQ}) \compose
 (\Id_P \juxta R2)
 \in \Hom{Kh}{\Kh{P}}{\Kh{R \juxta \rQ}}.$$
We propose that the inverse of this construction is given by
$$\Hom{Kh}{\Kh{P}}{\Kh{R \juxta \rQ}} \ni g \mapsto
F^{-1}(g) = (\Id_R \juxta R2^{-1}) \compose (g \juxta \Id_{\Q}).$$

The composition $F^{-1} \compose F$ applied to a chain map $f$ is
$$(\Id_R \juxta R2^{-1}) \compose (((f \juxta
 \Id_{\Q}) \compose
 (\Id_P \juxta R2)) \juxta \Id_{\Q}) =
 \mathfig{0.2}{duality/invertibility_stack}.$$
To see that this just $f$, we can do some tensor category arithmetic;%
\begin{align*}
 F^{-1}(F(f)) & = (f \juxta \Id_{\mathfig{0.06}{duality/two_strands}}) \compose (\Id_P \juxta (\Id_\Q \juxta R2^{-1} \compose R2 \juxta \Id_\Q)) \\
              & = \mathfig{0.2}{duality/invertibility_stack2} \\
              & = f.
\end{align*}
The critical step in this calculation came at the end, in claiming
that $(\Id_\Q \juxta R2^{-1}) \compose (R2 \juxta \Id_\Q) = \Id_\Q$.
This is exactly checking MM9, the ninth movie move. Although it
strains the logical order of the paper somewhat, we'll postpone that
calculation until \S \ref{sssec:MM9}, where we do all the other
movie moves, being careful to point out that we don't use any of the
results of this section while checking MM9.

A very similar argument shows $F(F^{-1}(g))$ is also just $g$.

The case $Q=\rQ$ is very similar.

Next, we deal with the case that the tangle $Q$ is just an arc,
$\rst$. This time, the map $F$ is given by $$F(f) = (f \juxta
\Id_\st) \compose (\Id_P \juxta
\rotatemathfig{0.05}{90}{cobordisms/cap_bdy_right}),$$ with inverse
$$F^{-1}(g) = (\Id_R \juxta
\mathfig{0.05}{cobordisms/saddle_unoriented2}) \compose (g \juxta
\Id_\rst).$$

The argument that $F$ and $F^{-1}$ are inverses is even easier than
before; some formal tensor category arithmetic and cobordism
arithmetic is all we need. For example,
\begin{align*}
F(F^{-1}(g)) = \mathfig{0.2}{duality/invertibility_stack3} = g.
\end{align*}

The other three cases where $Q$ is an arc are very similar.
\end{proof}

We now get an easy corollary, which you should think of as a nice
analogue of Bar-Natan's result about simple tangles in
\cite{MR2174270}.
\begin{cor}
\label{cor:simple_tangles}%
Let $T_1$ and $T_2$ be tangles with $k$ endpoints such that
$\refl{T_1}T_2$ is an unlink with $m$ components. Then the space of
chain maps modulo chain homotopy from $\Kh{T_1}$ to $\Kh{T_2}$ in
grading $m-k$ is 1-dimensional, and all chain maps of grading higher
than $m-k$ are chain homotopic to zero.
\end{cor}
\begin{proof}
By Proposition \ref{prop:duality}
\begin{eqnarray*}
    \Hom{Kh}{T_1}{T_2} &\Iso& \Hom{Kh}{\emptyset}{\refl{T_1}T_2}\grading{-k} \\
        &\Iso& \Kh{\refl{T_1}T_2}\grading{-k} \\
        &\Iso& (\R\grading{-1}\directSum\R\grading{+1})^{\tensor m}\grading{-k}
\end{eqnarray*}
\end{proof}

The next corollary is well known in the field, but perhaps worth
stating again.
\begin{cor}
\label{cor:reidemeister-maps-unique}%
The chain maps defined for the three Reidemeister moves in \S
\ref{ssec:disoriented_khovanov_homology} are, up to chain homotopy and scalar multiples,
the unique chain maps between the complexes in the appropriate grading.
\end{cor}

\subsection{Movie moves}
\label{ssec:movie_moves}

In this section, we'll complete the proofs of Theorems
\ref{thm:isofunctoriality} and \ref{thm:cobfunctoriality}, by
checking that changing the presentation of a cobordism by a movie
move does not change the associated chain map.

We'll first prove some preparatory lemmas, which will significantly
reduce the computational burden.

\begin{defn}
\label{defn:homotopy_isolation}%
Say $C^\bullet$ is a complex in some additive category, and $A$ is a
direct summand of some $C^i$. We say $A$ is homotopically isolated
if for any homotopy $h:C^\bullet \To C^{\bullet-1}$, the restriction
of $dh + hd$ to $A$ is zero.

If we're in a graded category then $A$ is homotopically isolated if
$dh + hd$ is zero for every grading $0$ homotopy $h$.
\end{defn}

\begin{lem}
\label{lem:loopless}%
Say $C^\bullet$ is the complex associated to some tangle diagram (so a
complex in the category of abstract disoriented cobordisms), and say
$A$ is a smoothing appearing as a direct summand of some step of the
complex. If $A$ does not contain any loops, and is not connected by
differentials to diagrams containing loops, then $A$ is
homotopically isolated.
\end{lem}
\begin{proof}
This is easy from the definition of the invariant in Equation
\ref{eq:invariant}. If two smoothings $B$ and $C$ are connected by a
differential $d: B \To C$, $C$ appears with a grading shift one more
than that of $B$. Thus a homotopy $h: C \To B$ would have to have
`bare' grading $+1$, but there are no positive grading morphisms
between loopless diagrams, by Euler characteristic considerations.
\end{proof}

\begin{lem}
\label{lem:homotopy_isolation}%
In each of movie movies 6 through 8, and in movie moves 11, 13 and
15, every smoothing in the complex associated to the initial frame
is homotopically isolated.
\end{lem}
\begin{proof}
This is trivial; no loops occur anywhere in the complex associated
to the initial frame.
\end{proof}

We don't need to say anything about homotopy isolation in MM9,
because we won't be using any of these simplifying lemmas in that
case---instead, the complete calculations are necessary for the sake
of Proposition \ref{prop:duality} on duality for Khovanov homology.

We can't say anything about homotopy isolation in MM12 and MM14,
because, when reading backwards in time, there aren't any isolated
objects! We will also use homotopy isolation in MM10, but identifying a
different smoothing in each of the many variations; the
details are in \S \ref{sssec:MM10}.

\begin{lem}
\label{lem:homotopy_isolation_factor}%
Suppose $f$ and $g$ are chain maps between the complexes $\Kh{T_1}$
and $\Kh{T_2}$, and we know $f \htpy \alpha g$ for some $\alpha$. If
$f$ and $g$ agree on some homotopically isolated object in the
complex $\Kh{T_1}$, say $O$, then  in fact $f \htpy g$ are
actually homotopic.
\end{lem}
\begin{proof}
On $O$, $f - \alpha g = dh+hd = 0$, so $f = \alpha g = g$. If $g$ is
just $0$, then $f$ is zero too, so $f$ and $g$ are trivially
homotopic. Otherwise $\alpha$ must be $1$, so $f$ and $g$ are
homotopic.
\end{proof}

Finally, we observe that Corollary \ref{cor:simple_tangles} applies
to every movie move. The join of the initial and final tangle is
always just an unlink, so the relevant space of chain maps modulo
homotopy is always one dimensional. Combined with the lemmas above,
we see that every movie move must come out right up to a multiple (in $\Integer[\frac{1}{2},\omega]$), and
to detect this multiple we can simply look at the restriction of the map
to a single homotopically isolated object. (Remembering, of course,
that MM9, MM12, and MM14 take a little more work; MM9 because there
we don't have access to any of the results on duality,
in particular Corollary \ref{cor:simple_tangles}, and MM12 and MM14
because we can't find homotopically isolated objects in the reverse
time direction.) Movie moves MM1 through MM10 describe isotopies, not general cobordisms, so there any multiple would actually have to be a unit.

In the calculations for MM6, MM8, MM9 and MM14, we'll explicitly keep
track of the ordering of the crossings. In all of the other
calculations, it turns out the ordering of crossings is irrelevant;
using the tricks described above, we only need to look at the action
of the chain maps on part of the complex, and in most cases any
crossing reordering maps automatically act on the objects we're
interested in by $+1$, simply because there's at most one crossing
which has been resolved disorientedly.

\subsubsection{MM1-5}
The first five movie moves are trivial; they simply say that a
Reidemeister move followed by its inverse is the identity.

\subsubsection{MM6-10}
Movie moves 6 through 10 involve no Morse moves, and so are
reversible. We only need to check one time direction.

In the following calculations (and those for MM11-15), red and
purple bands appearing in diagrams in complexes are simply a hint to
the reader, marking where crossings appeared in the original tangle.
(We hope they don't obscure too much for reader looking at a black and white printout.)

\subsubsection*{MM6}
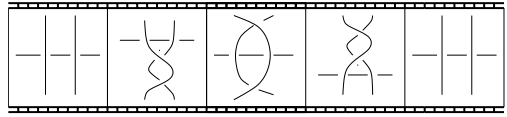
\label{sssec:MM6}
$$\mathfig{0.5}{movie_moves/MM6}$$

\begin{figure}[ht]
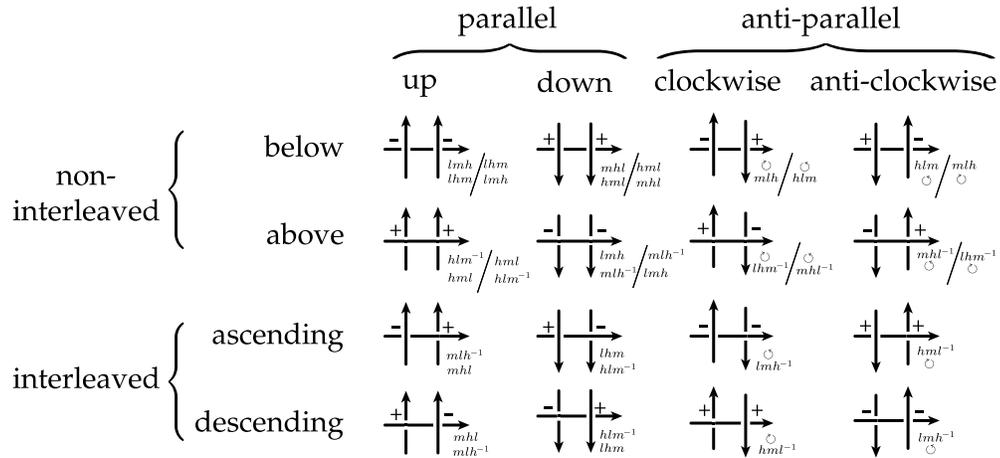

$$\mathfig{0.99}{movie_moves/MM6/table_of_variations}$$
\caption{16 variations for the initial frame of MM6.}\label{fig:MM6-variations}
\end{figure}

There are 24 variations of MM6. To see this we'll first of all make use of rotational symmetry to require that the 'horizontal' strand (the one
not involved in either R2 move) points from left to right. There are then sixteen possibilities for the initial frame of the movie move; these come from four choices of height orderings and four choices of orientations.
The horizontal strand can either lie entirely above or entirely below the two vertical strands ('non-interleaved'), or it may pass under one and over the
other ('interleaved', 'ascending' or 'descending').
The two vertical strands may be either parallel or anti-parallel. When they are parallel, they may point up or down, and when they
are anti-parallel they may have a clockwise or anti-clockwise orientation. All of these variations are displayed in Figure \ref{fig:MM6-variations}.

Note that the interleaved variations were not treated at all in versions of this paper before `\textbf{v2}' on the arXiv (see \S \ref{ssec:changelog}), or in Caprau's paper \cite{caprau} on the disoriented version of Khovanov homology.

Further, the eight variations in which the strands are 'non-interleaved' (the first two rows of Figure \ref{fig:MM6-variations}) each have two sub-variations,
which we don't see until the second frame of the movie. Of the two vertical strands, either one can pass above the other during the R2 moves; in Figure \ref{fig:MM6-variations}, the 'left passing above the right' sub-variation is listed to the left of the slash. In the 'interleaved' variations, there is no choice here.

We will thus treat four major cases,
\begin{itemize}
\item non-interleaved, parallel variations,
\item non-interleaved, anti-parallel variations,
\item interleaved, parallel variations and
\item interleaved, anti-parallel variations.
\end{itemize}

\paragraph{Non-interleaved parallel variations}
There are four possible initial frames which are 'non-interleaved' and have parallel vertical strands. Each of these
initial frames has two possible sub-variations, depending on the relative heights of the vertical strands during the R2 moves.
For each of the four initial frames, we will treat uniformly the sub-variations in which the upper $R2$-induced
crossing is negative and the lower one is positive, and then indicate how to treat the other four sub-variations.

Recall that our lemmas encapsulating the details of the $R3$ variations require that we
separate the initial and final complexes into layers $\mathcal{O}$ and $\mathcal{P}$
by resolving a crossing. Maneuvering through the pair of $R3$s in
this movie move is most efficiently managed by resolving the
$R2$-induced crossings: the upper one for the first $R3$, and the
lower one for the second $R3$. Notice that since the upper crossing
is negative, the first $R3$ will have homological ordering $\OP$,
while the second $R3$ will have ordering $\PO$. It's also worth mentioning
that the horizontal strand could be above or below the vertical
ones, meaning that these two crossing could be either the high or
low crossings in their respective $R3$ moves. However, Lemmas
\ref{lem:first-R3-general}, \ref{lem:second-R3-general}, and
\ref{lem:third-R3-general} work regardless\footnote{Recall the paragraphs following the statements of these Lemmas.} of whether the resolved
crossing is high or low, so we needn't treat them separately.

Our `bundle' of maps for this subcase, then, will look like this:

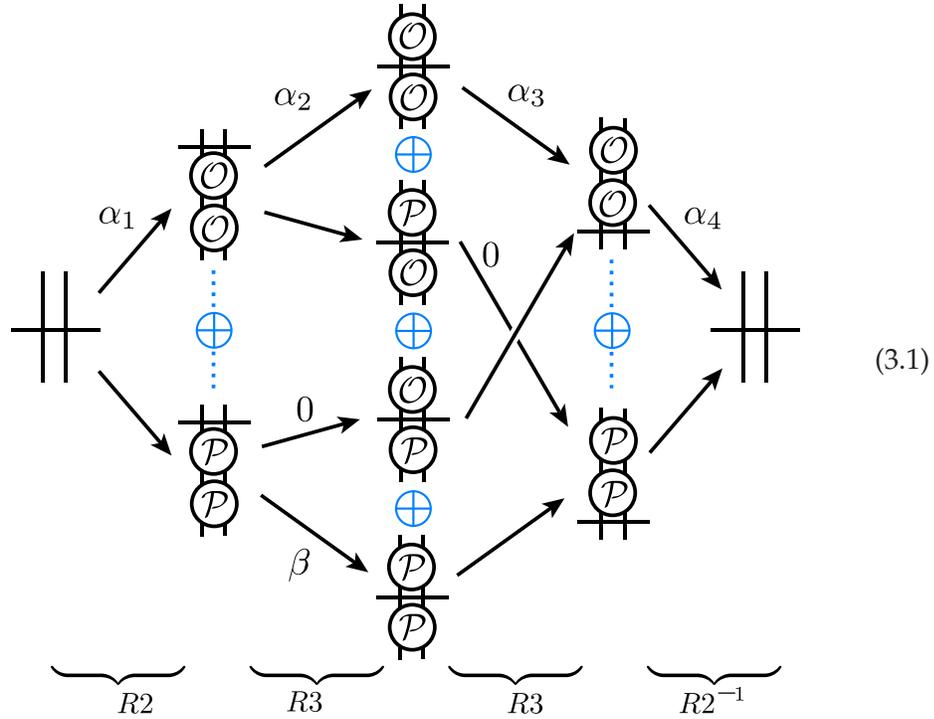
\begin{equation}
\label{eq:MM6_compositions}%
\mathfig{0.8}{movie_moves/MM6/MM6_2}
\end{equation}

In this diagram, the $\mathcal{O}$s and $\mathcal{P}$s describe
whether the indicated crossing resolution has strands orthogonal
or parallel to the horizontal strand. For example
$\mathfig{0.06}{movie_moves/MM6/MM6_4}$ is our notation for
$\mathfig{0.06}{movie_moves/MM6/MM6_5}$. Also, we've
cheated slightly with this diagram: the fourth column should
contain two additional summands, those with mixed
$\mathcal{O}$s and $\mathcal{P}$s. However, while there are
non-zero maps into these summands, the $R2^{-1}$ maps out are always
zero. Thus we needn't excessively complicate things with their
presence.

We're left with a sum of four compositions. The two middle
compositions are both zero, as each contains a leg (labelled with
`$0$') that's zero by Lemma \ref{lem:first-R3-general}. The top
composition ($\alpha_{i}$'s) is just the identity: $\alpha_1$ and $\alpha_4$ are components of $R2a$ moves,
and $\alpha_2$ and $\alpha_3$ are each the identity, by Lemma
\ref{lem:second-R3-general} (each map is a component of the $\OO$ map; when the horizontal strand lies below,
the R3 moves are $lmh,lhm,mhl$ and $hml$, which are exactly the four for which the $\OO$ part of the R3 map is the identity, and when the
horizontal strand lies above, the R3 moves are $hml,hlm,lmh$ and $mlh$, which are exactly the four for which the $\OO$ part of the mirror image
R3 map is the identity). The bottom composition is slightly more
mysterious, but we see that the map $\beta$ sends homologically
extreme smoothings to zero by Lemma \ref{lem:third-R3-general}.
Thus, if we choose an extreme smoothing to begin with, for example
the doubly oriented one, it will necessarily map to an extreme smoothing in
$\mathfig{0.06}{movie_moves/MM6/MM6_3}$, and thence to zero.
Further, as mentioned before, any initial smoothing here is
homotopically isolated, so the computation with this particular
smoothing suffices.

The other four sub-variations, in which the signs of the
$R2$-induced crossings are reversed, are proven analogously: note
that Equation \eqref{eq:MM6_compositions} will then have all
$\mathcal{O}$s and $\mathcal{P}$s swapped.

\paragraph{Non-interleaved anti-parallel variations}
Let's consider first those cases in which the vertical strands are oriented in the anti-clockwise direction,
so the left vertical strand is oriented downward, and
the right upward. Again we'll be referring to Equation
\eqref{eq:MM6_compositions}. Consider
the smoothing $\mathfig{0.1}{movie_moves/MM6/MM6f1_trace}$. Since we are looking at non-interleaved anti-parallel variations,
the two signs of the initial crossings differ, and so this resolution has homologically extreme height. In particular, when the horizontal strand
is below the vertical strands, this resolution has height $+1$, and when the horizontal strand is above the vertical strands, it has height $-1$.

The composition $\alpha_4 \compose \alpha_3 \compose \alpha_2 \compose \alpha_1$ then looks like

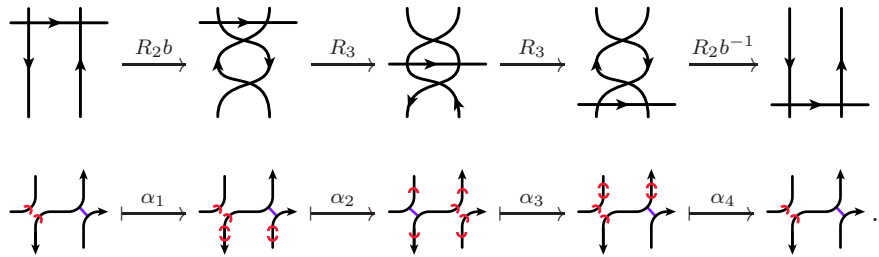
\begin{align*}
\newcommand{\pd}[1]{\mathfig{0.1}{movie_moves/MM6/#1}}
\newcommand{\pa}[1]{\mathfig{0.05}{movie_moves/MM6/#1}}
\xymatrix@R-7.5mm{
    \pd{MM6f5} \ar[r]^{R_2b} & \pd{MM6f4} \ar[r]^{R_3} & \pd{MM6f3} \ar[r]^{R_3} & \pd{MM6f2} \ar[r]^{{R_2b}^{-1}} & \pd{MM6f1} \\
    \\
    \pd{MM6f1_trace} \ar@{|->}[r]^{\alpha_1} & \pd{MM6f4_trace} \ar@{|->}[r]^{\alpha_2} & \pd{MM6f3_trace} \ar@{|->}[r]^{\alpha_3} & \pd{MM6f2_trace} \ar@{|->}[r]^{\alpha_4} & \pd{MM6f1_trace}. \\
}
\end{align*}


We now need to describe the maps $\alpha_i$, using our definitions of the
R2 chain maps from Figure \ref{fig:R2b_maps} for $\alpha_1$ and $\alpha_4$, and Lemma \ref{lem:second-R3-general} for $\alpha_2$ and $\alpha_3$.
We use the usual chain maps when the horizontal strand lies behind the others, and the mirror image chain maps when it is in front.
This description comes in three steps; first the underlying surfaces, ignoring disorientation data, then any associated coefficients, and finally
the arrangement of disorientation seams. The underlying surface for each map is simply the cylinder over the initial (and final) resolution.
Figure \ref{fig:R2b_maps} shows that $\alpha_1$ carries no coefficient, while $\alpha_4$ carries a coefficient of $\omega^2$.

For the R3 coefficients, notice that the R3 moves occurring in this configuration are one of the following pairs:
$hlm$/$\circlearrowleft$, $mlh$/$\circlearrowright$, $\overline{R3_{mhl}}^{-1}$/$\overline{R3_\circlearrowleft}$, or $\overline{R3_{lhm}}^{-1}$/$\overline{R3_\circlearrowright}$.
Our computation involves an extreme resolution on the $\mathcal{O}$ layer in both R3 moves; let $\sigma_2$ and $\sigma_3$ be the coefficients on the appropriate height (either high or low) components of the $\OO$ part of the corresponding R3 moves.
Then, according to Lemma \ref{lem:second-R3-general} and Corollary \ref{cor:R3-inverses}, it is always the case that one of the $\sigma_i$'s is $-1$, while the other is $\omega^2$.

In each of the four cases, a (cancelling) pair of reordering signs is needed. 

Finally, we add disorientation seams:
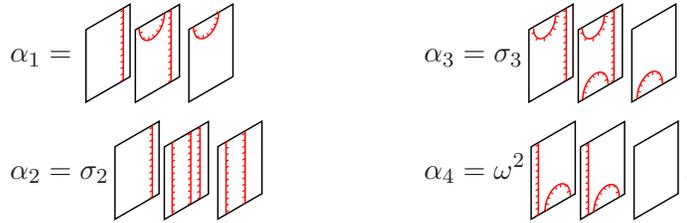
\begin{align*}
 \alpha_1      & = \mathfig{0.15}{movie_moves/MM6/MM6_example1_cob1} &
 \alpha_3      & = \sigma_3 \mathfig{0.15}{movie_moves/MM6/MM6_example1_cob3}
 \displaybreak[1] \\
 \alpha_2      & = \sigma_2 \mathfig{0.15}{movie_moves/MM6/MM6_example1_cob2} &
 \alpha_4      & = \omega^2 \mathfig{0.15}{movie_moves/MM6/MM6_example1_cob4}.
\end{align*}
Notice that in $\alpha_2$, the left-most seams on the second two sheets are vertical, as the associated crossing is not involved in the $R3$ move; the other seams are the unique minimal ones connecting the remaining eight disorientation marks. Similarly, in $\alpha_3$ the right-most seams on the first two sheets are vertical, leaving the others to be determined by minimality.

Thus, our composition $\alpha_4 \circ \alpha_3 \circ \alpha_2 \circ \alpha_1$ is just
$$\omega^2 \sigma_2 \sigma_3 \ \mathfig{0.2}{movie_moves/MM6/MM6_example1_cobf1} \  = \  - \ \mathfig{0.2}{movie_moves/MM6/MM6_example1_cobf2} = -\omega^2 \Id.$$

Of course, starting with an extreme object also guarantees this $\alpha$ composition is the only one we need to worry about, as $\beta=0$ from Lemma \ref{lem:third-R3-general}.

The argument for the case in which the left vertical strand is oriented upward, and the right downward, is essentially the same.

\paragraph{Interleaved variations}
There are eight variations, and essentially two distinct computations will cover them all. Let's start with $hml^{-1}/\circlearrowright$, $\circlearrowleft/lmh^{-1}$, $mlh^{-1}/mhl$, and $lhm/hlm^{-1}$; we'll show the calculation for the first, and explain the necessary alterations for the other three versions.

{
\newcommand{\pd}[1]{\mathfig{0.1}{movie_moves/MM6/#1}}
\newcommand{\pa}[1]{\mathfig{0.05}{movie_moves/MM6/#1}}

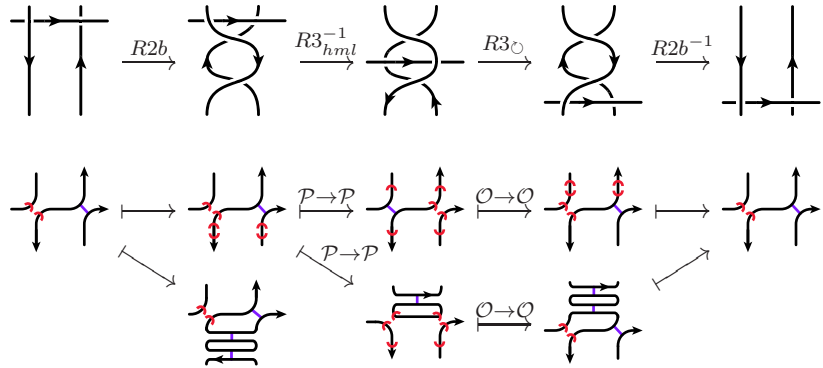
\begin{align*}%
\xymatrix@R-7.5mm@C-1.5mm{
    \pd{MM6intaf5} \ar[r]^{R2b} & \pd{MM6intaf4} \ar[r]^{R3_{hml}^{-1}} & \pd{MM6intaf3} \ar[r]^{R3_{\circlearrowright}} & \pd{MM6intaf2} \ar[r]^{{R2b}^{-1}} & \pd{MM6intaf1} \\
    \\
    \pd{MM6f1_trace} \ar@{|->}[dr] \ar@{|->}[r] & \pd{MM6f4_trace} \ar@{|->}[dr]^{\PP} \ar@{|->}[r]^{\PP} & \pd{MM6f3_trace} \ar@{|->}[r]^{\OO} & \pd{MM6f2_trace} \ar@{|->}[r] & \pd{MM6f1_trace} \\
     & \pd{MM6intaf4r2_trace} & \pd{MM6intaf3r2_trace} \ar@{|->}[r]^{\OO} & \pd{MM6intaf2r2_trace} \ar@{|->}[ur] &
}
\end{align*}
}

Notice that our first $R3$ map is ordered $\OP$ and the second $\PO$, each with the high crossing resolved, and that the maps for these moves are labeled by their source and target layers; in particular, the initial $\mathcal{O}$ layer for the second move and the final $\mathcal{P}$ layer for the first move coincide.

Lemma \ref{lem:first-R3-general} tells us there are only three compositions we need to keep track of here. The first map into the second row has an extreme target in the initial $\mathcal{P}$ layer of $R3_{hml}^{-1}$, which thereafter maps to zero by Lemma \ref{lem:third-R3-general}. The composition including the rest of the second row contains a sphere; this is because, disregarding coefficients and disorientation seams, the first and third maps are cylinders from the $R2b$ chain map definitions and Lemma \ref{lem:second-R3-general}, the second map contains a cup by Lemma \ref{lem:third-R3-general}, and the fourth map, an $R2b$ untuck, contains a cap. Thus we're left with the first row, and a brief look at Lemmas \ref{lem:second-R3-general} and \ref{lem:third-R3-general}, and a check that there are no signs from crossing reorderings, confirms that this composition looks like

$$- \ \mathfig{0.2}{movie_moves/MM6/MM6_example1_cobf1} \  = \  - \ \mathfig{0.2}{movie_moves/MM6/MM6_example1_cobf2} = -\omega^2 \Id.$$

The calculations for the $\circlearrowleft/lmh^{-1}$, $mlh^{-1}/mhl$, and $lhm/hlm^{-1}$ variations are very similar. For $\circlearrowleft/lmh^{-1}$, the initial object will have a disoriented left crossing and an oriented right crossing, and we'll resolve each $R3$ move into layers using the low crossing. Thus we'll need to compute using the mirror image maps, which will introduce an extra factor of $-\omega^2$. The $mlh^{-1}/mhl$ and $lhm/hlm^{-1}$ variations are even easier: we start with the doubly oriented object in each case, and resolve into layers using the high crossings or the low crossings, respectively. Crossing reordering maps are trivial in all three of these additional variations, and the overall coefficient for each is just 1.

The computations for $\circlearrowright/hml^{-1}$, $lmh^{-1}/\circlearrowleft$, $mhl/mlh^{-1}$, and $hlm^{-1}/lhm$ are somewhat different; again, we'll explicitly show the first.

{
\newcommand{\pd}[1]{\mathfig{0.1}{movie_moves/MM6/#1}}
\newcommand{\pa}[1]{\mathfig{0.05}{movie_moves/MM6/#1}}

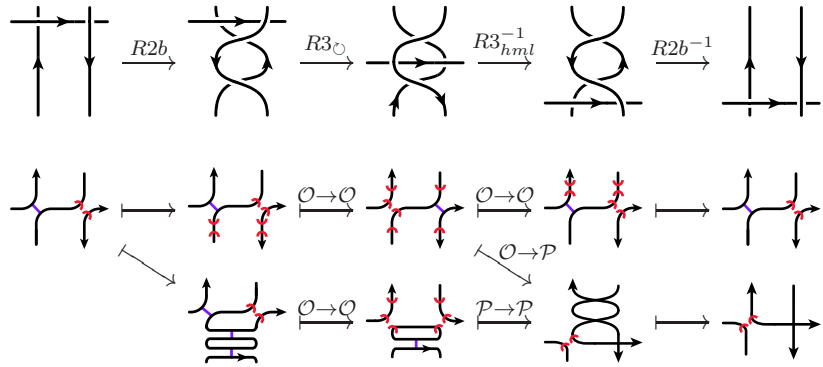
\begin{align*}%
\xymatrix@R-7.5mm@C-1.5mm{
    \pd{MM6intbf5} \ar[r]^{R2b} & \pd{MM6intbf4} \ar[r]^{R3_{\circlearrowright}} & \pd{MM6intbf3} \ar[r]^{R3_{hml}^{-1}} & \pd{MM6intbf2} \ar[r]^{{R2b}^{-1}} & \pd{MM6intbf1} \\
    \\
    \pd{MM6intbf1_trace} \ar@{|->}[r] \ar@{|->}[dr] & \pd{MM6intbf4_trace} \ar@{|->}[r]^{\OO} & \pd{MM6intbf3_trace} \ar@{|->}[r]^{\OO} \ar@{|->}[dr]^{\OP} & \pd{MM6intbf2_trace} \ar@{|->}[r] & \pd{MM6intbf1_trace} \\
     & \pd{MM6intbf4r2_trace} \ar@{|->}[r]^{\OO} & \pd{MM6intbf3r2_trace} \ar@{|->}[r]^{\PP} & \pd{MM6intbf2r2_trace} \ar@{|->}[r] & \pd{MM6intbf1r2_trace}
}
\end{align*}
}

Now our first $R3$ map is ordered $\PO$ with the high crossing resolved, and the second is ordered $\OP$ with the low crossing resolved. Again, we'll keep track of the layers to which objects belong by referring to the labels on the maps.

By Lemma \ref{lem:first-R3-general}, we have three compositions to consider. Two of them factor through the second row, and thus map to a complex with the left crossing disoriented; since our map is a multiple of the identity, these compositions must sum to zero. So we're left with the first row. Using Lemma \ref{lem:second-R3-general} (and its mirror image variant for the second $R3$), the $R2b$ map definitions, and the fact that crossing reorderings give a minus sign here, it's straightforward to verify this composition is given by

$$- \ \mathfig{0.2}{movie_moves/MM6/MM6_example2_cobf1} \  = \  - \ \mathfig{0.2}{movie_moves/MM6/MM6_example2_cobf2} = -\omega^2 \Id.$$

There are a few modifications necessary for $lmh^{-1}/\circlearrowleft$, $mhl/mlh^{-1}$, and for $hlm^{-1}/lhm$. In the $lmh^{-1}/\circlearrowleft$ case, we start with the object with oriented left crossing and disoriented right crossing, and resolve the first $R3$ on low and the second on high. A crossing reordering sign gives us an overall coefficient of $-\omega^2$. For each of $mhl/mlh^{-1}$ and $hlm^{-1}/lhm$ our initial object will be the doubly oriented one, so crossing reordering maps act trivially. In the calculations, $hlm^{-1}$ and $mlh^{-1}$ should be resolved on low, while $mhl$ and $lhm$ should be resolved on high. An overall coefficient of 1 will result in each of these cases.

\subsubsection*{MM7}
$$\mathfig{0.4}{movie_moves/MM7}$$
We need to consider four variations of MM7, depending on the
orientation of the strand, and whether the `first' crossing is
positive or negative. It's easy to check that reversing orientations
in the two subsequent calculations doesn't change the result.

First we deal with a positive crossing:
\begin{align*}
\newcommand{\pd}[1]{\mathfig{0.1}{movie_moves/MM7/#1}}
\newcommand{\pa}[1]{\mathfig{0.05}{movie_moves/MM7/#1}}
\xymatrix@R-7.5mm@C+7mm{
    \pd{MM7f1} \ar[rr]^{R_1a} & & \pd{MM7f2} \ar[r]^{R_1b} & \pd{MM7f3} \ar[r]^{{R_2b}^{-1}} & \pd{MM7f4} \\
    \\
    \pd{MM7f1} \ar@{|->}[rr]^{\frac{1}{2}\left(\mathfig{0.05}{cobordisms/curtain_and_handle_creation}
        - \omega^{-2} \mathfig{0.05}{cobordisms/curtain_handle_and_disc_creation}\right)} & & %
    \pd{MM7f2_trace} \ar@{|->}[r]^{\mathfig{0.075}{cobordisms/curtain_and_cylinder_and_disc_creation}} & %
    \pd{MM7f3_trace} \ar@{|->}[r]^{- \omega^2 \mathfig{0.065}{movie_moves/MM7/untuck_map}} & %
    \pd{MM7f4} %
}.
\end{align*}

Composing, we see that the second term of the first map gives zero
when composed with the later maps. Cancelling the factor of
$\frac{1}{2}$ with the torus, we get $-\omega^2$ times the identity.

For the negative crossing, we have
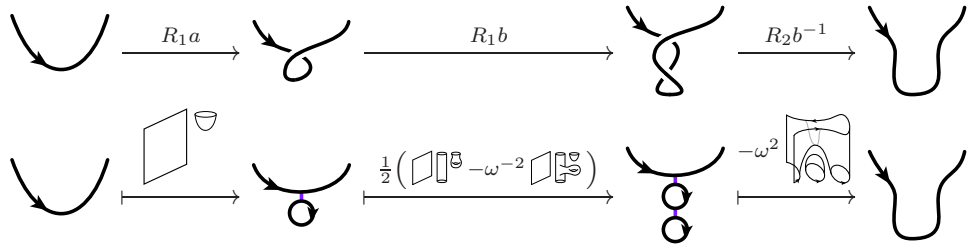
\begin{align*}
\newcommand{\pd}[1]{\mathfig{0.1}{movie_moves/MM7/#1}}
\newcommand{\pa}[1]{\mathfig{0.05}{movie_moves/MM7/#1}}
\xymatrix@R-7.5mm@C+7mm{
    \pd{MM7f1} \ar[r]^{R_1a} & \pd{MM7f2n} \ar[rr]^{R_1b} & & \pd{MM7f3n} \ar[r]^{{R_2b}^{-1}} & \pd{MM7f4} \\
    \\
    \pd{MM7f1} \ar@{|->}[r]^{\mathfig{0.075}{cobordisms/curtain_and_disc_creation}} & %
    \pd{MM7f2_trace} \ar@{|->}[rr]^{\frac{1}{2}\left(\mathfig{0.05}{movie_moves/MM7/curtain_and_cylinder_and_handle_creation}
        - \omega^{-2} \mathfig{0.05}{movie_moves/MM7/curtain_and_cylinder_with_handle_and_disc_creation}\right)} & & %
    \pd{MM7f3_trace} \ar@{|->}[r]^{- \omega^2 \mathfig{0.065}{movie_moves/MM7/untuck_map}} & %
    \pd{MM7f4} %
}
\end{align*}
and the composition is just the identity.

\subsubsection*{MM8}
$$\mathfig{0.5}{movie_moves/MM8}$$

This is the only movie move involving all three Reidemeister moves. There
are quite a few variations. By a rotation of the whole diagram, we can
assume the R1 move happens on the horizontal strand, beginning on the
right. Moreover, we can assume that the horizontal strand is oriented
right to left (otherwise, we can obtain this condition by a $\pi$
rotation of its time reversal).

There are then sixteen variations, depending on whether the vertical strand
lies above or below the horizontal strand, its orientation, the sign of
the crossing introduced by the first Reidemeister move in the first
frame, and finally whether the first Reidemeister move introduces a twist on the left or right side. The following diagram shows all the maps involved, independent of crossing sign choices and thus without disorientation marks (we will add them later):

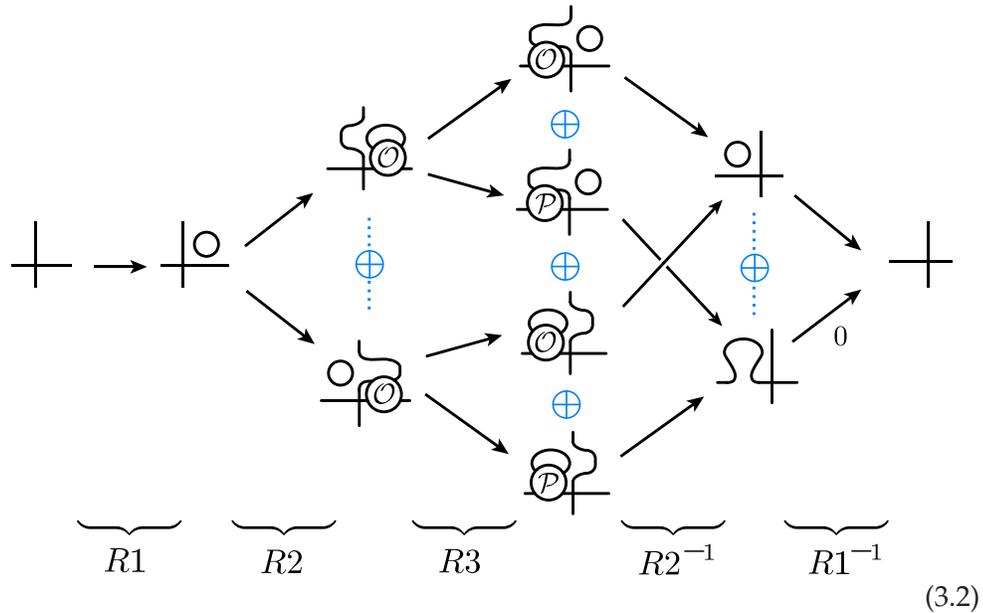
\begin{equation}
\label{eq:MM8_compositions}%
\mathfig{0.95}{movie_moves/MM8/MM8_compositions}
\end{equation}

Note that the crossing introduced by the $R1$ move is always either the low or high crossing in the $R3$ move,
so we will denote its resolution with either $\mathcal{O}$ or $\mathcal{P}$ as we did in the computation for MM6.
We can also observe that any map factoring through the resolution
$\mathfig{0.1}{movie_moves/MM8/killstate}$
must be zero, since this object maps to zero under $R1$ (see \S \ref{sec:R1-maps}). Thus we need only concern ourselves with the other two compositions in Equation \eqref{eq:MM8_compositions}.

Let's first treat the positive twist.  We'll show
calculations for the case in which the vertical strand is oriented downward (but ignore
whether the twist appears on the left or right side of the horizontal strand; this barely changes any of the
calculations). Also, our computation will work regardless of whether the vertical strand is above or below the horizontal strand.

{
\newcommand{\pd}[1]{\mathfig{0.1}{movie_moves/MM8/#1}}%
\newcommand{\pa}[1]{\mathfig{0.05}{movie_moves/MM8/#1}}%
Beginning with a downward-oriented vertical strand, the two relevant compositions are
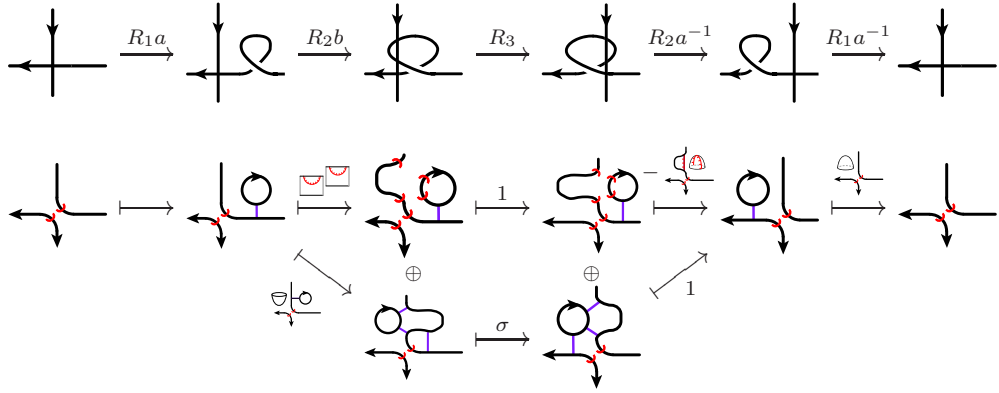
\begin{align*}%
\xymatrix@R-7.5mm@C-1.5mm{
    \pd{MM8v1f1} \ar[r]^{R_1a} & \pd{MM8v1f2} \ar[r]^{R_2b} & \pd{MM8v1f3} \ar[r]^{R_3} & \pd{MM8v1f4} \ar[r]^{{R_2a}^{-1}} & \pd{MM8v1f5} \ar[r]^{{R_1a}^{-1}} & \pd{MM8v1f1} \\
    \\
    \pd{MM8v1f1_trace} \ar@{|->}[r]^{} & %
    \pd{MM8v1f2_trace} \ar@{|->}[r]^{\mathfig{0.05}{reidemeister_maps/R2b_tuck_map_sheets}} \ar@{|->}[dr]_{\pa{MM8v1_map2}} & %
    \pd{MM8v1f3_trace} \ar@{|->}[r]^{1} \ar@{}[d]|{\directSum} & %
    \pd{MM8v1f4_trace} \ar@{|->}[r]^{- \pa{MM8v1_map4}} \ar@{}[d]|{\directSum} & %
    \pd{MM8v1f5_trace} \ar@{|->}[r]^{\pa{MM8v1_map5}} & %
    \pd{MM8v1f1_trace} \\
     & &
     \pd{MM8v1f3_trace2} \ar@{|->}[r]^{\sigma} &
     \pd{MM8v1f4_trace2} \ar@{|->}[ur]_{1} & &
}
\end{align*}
where the initial $R_1a$ map is
$\pd{MM8v1f1_trace} \xrightarrow{\frac{1}{2}\left(\pa{MM8v1_map1a} - \omega^{-2} \pa{MM8v1_map1b}\right)} \pd{MM8v1f2_trace}$.
}
Here the R3 move is either $\overline{R3_{mlh}}$ or $R3_{hml}^{-1}$, depending on whether the vertical strand is in front or behind.

These chain map components come from \S \ref{sec:R1-maps}. In particular,
we use Lemma \ref{lem:second-R3-general} (or its `mirror image' analogue,
depending on whether the vertical or horizontal strand is on top), to see
what the $R3$ maps do. In the first row our object lies in the $\mathcal{O}$ layer, and thus maps via the identity. (Note that for $\overline{R3_{mlh}}$, when looking
for a description of the $\OO$ layer in Lemma \ref{lem:second-R3-general}, we actually need to look at the case corresponding to
$R3_{mhl}$, since $\overline{R3_{mlh}}$ is defined in terms of  $R3_{mhl}$.)
All we need to know about the $R3$ map in the second row, labelled by $\sigma$, is that, ignoring disorientation seams
it is the cylinder cobordism. This tells us that the
lower row has a spherical component and can thus be ignored. As such, the
composition simplifies to

{
\newcommand{\pd}[1]{\mathfig{0.1}{movie_moves/MM8/#1}}%
\begin{equation*}
\begin{split}%
 \pd{MM8v1_map5} \compose -\pd{MM8v1_map4} \compose \Id \compose \mathfig{0.1}{reidemeister_maps/R2b_tuck_map_sheets} \compose
    \frac{1}{2}\left(\pd{MM8v1_map1a} - \omega^{-2} \pd{MM8v1_map1b}\right) \\
  \shoveright{= -\omega^{2} \  \frac{1}{2} \ \mathfig{0.075}{cobordisms/torus} \pd{MM8v1f1_trace}} \\
  = -\omega^{2} \pd{MM8v1f1_trace}.
\end{split}
\end{equation*}
}

When the first Reidemeister move introduces a negative crossing, we
see instead
{
\newcommand{\pd}[1]{\mathfig{0.1}{movie_moves/MM8/#1}}%
\newcommand{\pa}[1]{\mathfig{0.05}{movie_moves/MM8/#1}}%
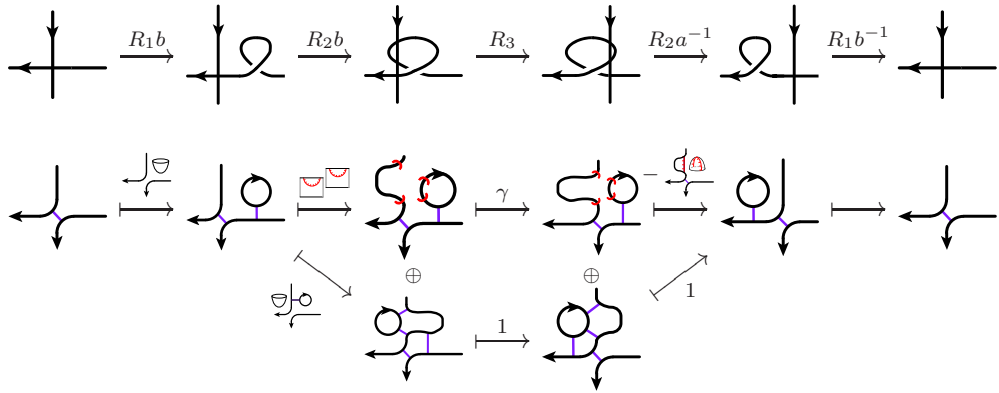
\begin{align*}%
\xymatrix@R-7.5mm@C-1.5mm{
    \pd{MM8v2f1} \ar[r]^{R_1b} & \pd{MM8v2f2} \ar[r]^{R_2b} & \pd{MM8v2f3} \ar[r]^{R_3} & \pd{MM8v2f4} \ar[r]^{{R_2a}^{-1}} & \pd{MM8v2f5} \ar[r]^{{R_1b}^{-1}} & \pd{MM8v2f1} \\
    \\
    \pd{MM8v2f1_trace} \ar@{|->}[r]^{\pa{MM8v2_map1}} & %
    \pd{MM8v2f2_trace} \ar@{|->}[r]^{\mathfig{0.05}{reidemeister_maps/R2b_tuck_map_sheets}} \ar@{|->}[dr]_{\pa{MM8v2_map2}} & %
    \pd{MM8v2f3_trace} \ar@{|->}[r]^{\gamma} \ar@{}[d]|{\directSum} & %
    \pd{MM8v2f4_trace} \ar@{|->}[r]^{- \pa{MM8v2_map4}} \ar@{}[d]|{\directSum} & %
    \pd{MM8v2f5_trace} \ar@{|->}[r] & %
    \pd{MM8v2f1_trace} \\
     & &
     \pd{MM8v2f3_trace2} \ar@{|->}[r]^{1} &
     \pd{MM8v2f4_trace2} \ar@{|->}[ur]_{1} & &
}
\end{align*}
where the final $R_1b^{-1}$ map is $\pd{MM8v2f5_trace} \xrightarrow{\frac{1}{2}\left(\pa{MM8v2_map5a} - \omega^{-2} \pa{MM8v2_map5b}\right)} \pd{MM8v2f1_trace}$.
}
Here the R3 move is either $\overline{R3_{lmh}^{-1}}$ or $R3_{mhl}^{-1}$, depending on whether the vertical strand is in front or behind.

This time the first row gives zero ($\gamma$ is some disoriented cylinder, so the composition contains a sphere), and we obtain
{
\newcommand{\pd}[1]{\mathfig{0.1}{movie_moves/MM8/#1}}%
\begin{equation*}
\begin{split}%
 \frac{1}{2}\left(\pd{MM8v2_map5a} - \omega^{-2} \pd{MM8v2_map5b}\right)
 \compose \Id \compose \Id \compose \pd{MM8v2_map2} \compose \pd{MM8v2_map1} \\
 \shoveright{ = \frac{1}{2} \mathfig{0.075}{cobordisms/torus} \pd{MM8v2f1_trace}} \\
 = \pd{MM8v2f1_trace}.
\end{split}
\end{equation*}
}

Changing the orientation of the vertical strand (for either a positive or negative twist) alters the computations only slightly (in particular the coefficient appearing on the R3 map is still always $+1$), and we obtain the same result: the coefficient, $1$ or $-\omega^2$, just depends on the sign of the crossing introduced by the first Reidemeister move.

\subsubsection*{MM9}
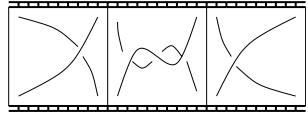
\label{sssec:MM9}
$$\mathfig{0.3}{movie_moves/MM9}$$
For MM9 we have to be particularly careful; the proof of
Proposition \ref{prop:duality} relied on this movie
move, so while checking MM9 we don't have access to any results
about the space of chain maps being one dimensional. Thus we'll
fully calculate the map, checking it's the identity on every object
in the complex associated to the initial tangle.

There are four variations of MM9; we can fix the orientation of one
strand, then have to deal with either orientation of the other
strand, and either sign for the crossing.

We'll do the calculations for both types of crossings, in a given
orientation. It's easy to see that changing an orientation
essentially interchanges these cases.

With a positive crossing, we have
{
\newcommand{\pd}[1]{\mathfig{0.1}{movie_moves/MM9/positive/#1}}
\newcommand{\pa}[1]{\mathfig{0.05}{movie_moves/MM9/positive/#1}}
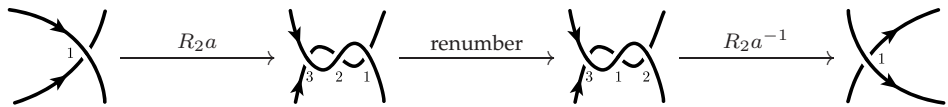
\begin{align*}
\xymatrix@R-7.5mm@C+12mm{
    \pd{MM9f1} \ar[r]^{R_2a} & \pd{MM9f2} \ar[r]^{\text{renumber}} & \pd{MM9f2renumbered} \ar[r]^{{R_2a}^{-1}} & \pd{MM9f3} \\
}
\intertext{and the components of the chain map are given by:}
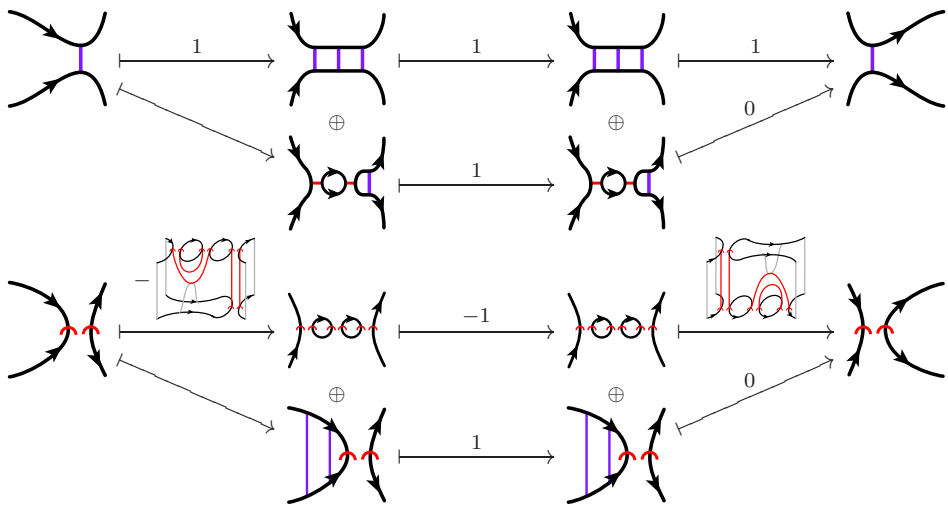
\xymatrix@R-7.5mm@C+12mm{
    \pd{MM9f1_trace} \ar@{|->}[r]^{1} \ar@{|->}[dr] & %
    \pd{MM9f2_trace1} \ar@{|->}[r]^{1} \ar@{}[d]|{\directSum} & %
    \pd{MM9f2_trace1} \ar@{|->}[r]^{1} \ar@{}[d]|{\directSum} & %
    \pd{MM9f3_trace} \\ %
    & \pd{MM9f2_trace2} \ar@{|->}[r]^{1} & %
      \pd{MM9f2_trace2} \ar@{|->}[ur]^{0} & \\\\%
    \pd{MM9f1_dtrace} \ar@{|->}[r]^{-\pd{tuck}} \ar@{|->}[dr] & %
    \pd{MM9f2_dtrace1} \ar@{|->}[r]^{-1} \ar@{}[d]|{\directSum} & %
    \pd{MM9f2_dtrace1} \ar@{|->}[r]^{\pd{untuck}} \ar@{}[d]|{\directSum} & %
    \pd{MM9f3_dtrace} \\ %
    & \pd{MM9f2_dtrace2} \ar@{|->}[r]^{1} & %
      \pd{MM9f2_dtrace2} \ar@{|->}[ur]^{0} & %
}
\end{align*}
}
and the composition is just the identity.

With a negative crossing, we have
{
\newcommand{\pd}[1]{\mathfig{0.1}{movie_moves/MM9/negative/#1}}
\newcommand{\pa}[1]{\mathfig{0.05}{movie_moves/MM9/negative/#1}}
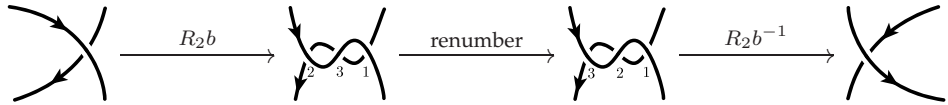
\begin{align*}
\xymatrix@R-7.5mm@C+12mm{
    \pd{MM9f1} \ar[r]^{R_2b} & \pd{MM9f2} \ar[r]^{\text{renumber}} & \pd{MM9f2renumbered} \ar[r]^{{R_2b}^{-1}} & \pd{MM9f3} \\
}
\intertext{with the components of the chain map being given by}
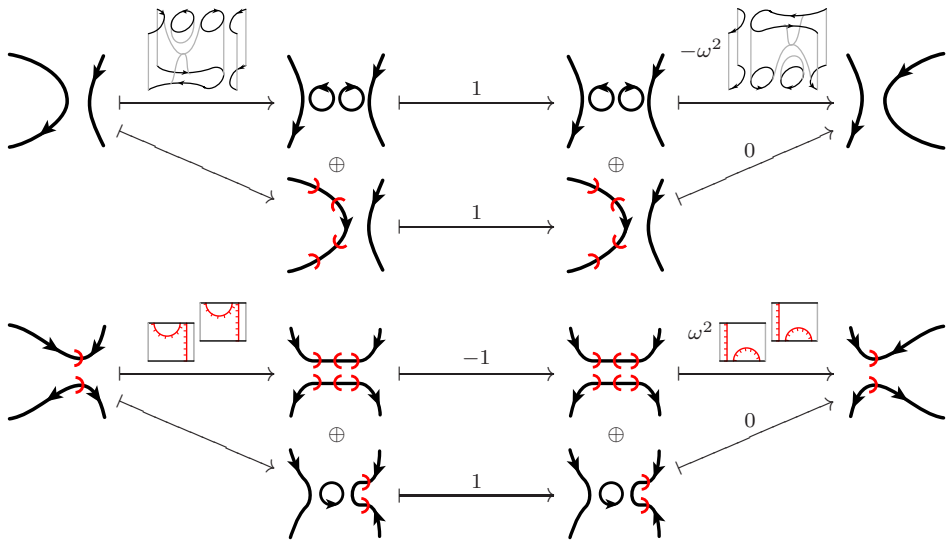
\xymatrix@R-7.5mm@C+12mm{
    \pd{MM9f1_trace} \ar@{|->}[r]^{\pd{tuck}} \ar@{|->}[dr] & %
    \pd{MM9f2_trace1} \ar@{|->}[r]^{1} \ar@{}[d]|{\directSum} & %
    \pd{MM9f2_trace1} \ar@{|->}[r]^{- \omega^2 \pd{untuck}} \ar@{}[d]|{\directSum} & %
    \pd{MM9f3_trace} \\ %
    & \pd{MM9f2_trace2} \ar@{|->}[r]^{1} & %
      \pd{MM9f2_trace2} \ar@{|->}[ur]^{0} & \\\\%
    \pd{MM9f1_dtrace} \ar@{|->}[r]^{\pd{tuck_sheets}} \ar@{|->}[dr] & %
    \pd{MM9f2_dtrace1} \ar@{|->}[r]^{-1} \ar@{}[d]|{\directSum} & %
    \pd{MM9f2_dtrace1} \ar@{|->}[r]^{\omega^2 \pd{untuck_sheets}} \ar@{}[d]|{\directSum} & %
    \pd{MM9f3_dtrace} \\ %
    & \pd{MM9f2_dtrace2} \ar@{|->}[r]^{1} & %
       \pd{MM9f2_dtrace2} \ar@{|->}[ur]^{0} & %
}
\end{align*}}
and the composition is $-\omega^2$ times the identity.

\subsubsection*{MM10}
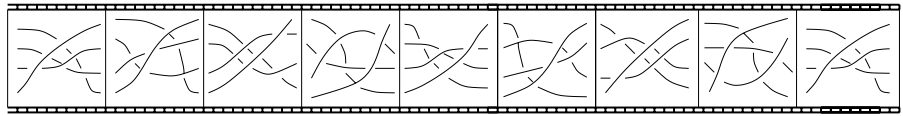
\label{sssec:MM10}
\begin{align*}
\mathfig{0.9}{movie_moves/MM10}
\end{align*}

This is the tetrahedron move, and is surprisingly easy. On the
other hand, there are a great many variations which we need to
treat.

Firstly, let's consider the case in which all strands are oriented
to the right. Here, all the crossings are positive, and if we
consider the object in the initial complex with homological height
zero (ie, we've smoothed every crossing in the oriented way), we
see that it is homotopically isolated.

Notice further that each of the eight $R3$ moves in this movie is
of type $R3_{lmh}$, and so has homological ordering $\PO$ when
resolving the highest (or lowest) crossing. The oriented smoothing
lives in the $\mathcal{O}$ layer, and thus maps via the identity
to the oriented smoothing in the next frame by Lemma
\ref{lem:second-R3-general}. Also, by Lemma
\ref{lem:first-R3-general}, there is no map to the $\mathcal{P}$
layer. It's easy to see that the same happens at each of the seven other
$R3$ moves, so we're just left with a string of identity maps: {
\newcommand{\pd}[1]{\mathfig{0.075}{movie_moves/MM10/#1}}
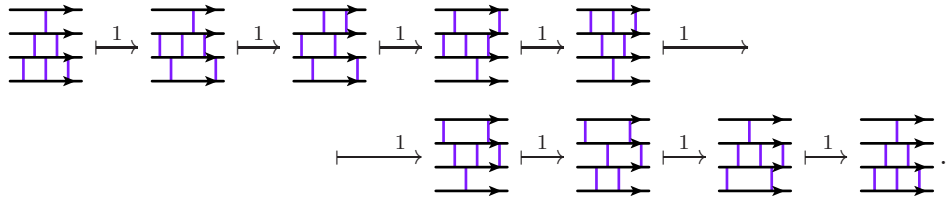
\begin{gather*}
\xymatrix@R-7.5mm@C-3mm{
    \pd{MM10f1_trace} \ar@{|->}[r]^{1} & %
    \pd{MM10f2_trace} \ar@{|->}[r]^{1} & %
    \pd{MM10f3_trace} \ar@{|->}[r]^{1} & %
    \pd{MM10f4_trace} \ar@{|->}[r]^{1} & %
    \pd{MM10f5_trace} \ar@{|->}[r]^{1} & \\ & & \ar@{|->}[r]^{1} & %
    \pd{MM10f6_trace} \ar@{|->}[r]^{1} & %
    \pd{MM10f7_trace} \ar@{|->}[r]^{1} & %
    \pd{MM10f8_trace} \ar@{|->}[r]^{1} & %
    \pd{MM10f1_trace}. \\ %
}
\end{gather*}
}%
Thus this movie induces the identity chain map.

\begin{figure}[ht]
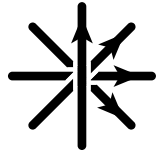

$$\mathfig{0.15}{movie_moves/MM10/non_generic_projection}$$
\caption{A non-generic projection corresponding to a MM10
2-cell.}\label{fig:MM10-projection}
\end{figure}

Beyond this, there are a frightening forty-eight variations. In the space
of tangle diagrams, MM10 corresponds to a codimension 2 stratum,
appearing as a non-generic projection in which four strands cross
at a point. (See Figure \ref{fig:MM10-projection}). Rotating the
projection to put the highest strand in a standard position, there
are then $3!$ height orderings we need to consider for the other
strands, and $2^3$ orientations.

It turns out, however, that every variation of MM10 is actually
equivalent, modulo MM6.

The idea, essentially, is to add an extra crossing to MM10. We can
do this at any adjacent pair of boundary points; for concreteness,
let's imagine adding an extra crossing at the top right, with
opposite sign to the crossing that already appears in the top
right in the first and last frames. There's now a pair of strands
carrying two crossings. We can now consider two different
variations of MM10, each of which involves only one of those two
crossings, and see that these two MM10 moves differ by some MM6
moves (and some `distant Reidemeister moves commute' moves).

More generally, we can stratify the space of smooth tangles so
that
in the dual cell complex (where a $k$-cell corresponds to a codimension $k$ stratum)%
\begin{itemize}
\item 0-cells correspond to tangles whose projection to $B^2$ is a
generic immersion. \item 1-cells correspond to Reidemeister moves.
\item 2-cells correspond to movie moves and pairs of distant
Reidemeister moves. \item 3-cells correspond to redundancies
amongst movie moves.
\end{itemize}%

\begin{figure}[ht]
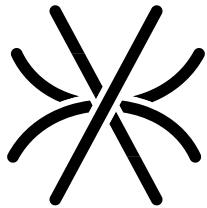

$$\mathfig{0.2}{movie_moves/MM10/3-cell}$$
\caption{A non-generic projection corresponding to a 3-cell
involving MM10 and MM6.}\label{fig:3-cell}
\end{figure}

If we consider a 3-cell dual to the non-generic projection shown
in Figure \ref{fig:3-cell}, we find that the 2-cells on its
boundary consist of two MM10 2-cells, four MM6 2-cells, and six
distant R-move 2-cells: see Figure \ref{fig:3-cell diagram}. Thus
invariance for one of the two MM10's, plus invariance for all
MM6's and pairs of distant Reidemeister moves, implies invariance
for the other MM10.

\begin{figure}[ht]
$$\mathfig{1.0}{movie_moves/MM10/3-cell_diagram}$$
\caption{The 3-cell for the singularity in Figure
\ref{fig:3-cell}, rotated 90 degrees. The 0-cells here are the
generic tangle projections neighboring this singularity, achieved
by untucking the curved strands ($z$ direction) and translating
the crossing ($x$ and $y$ directions). The 2-cells marked with an
asterisk correspond to distant Reidemeister
moves.}\label{fig:3-cell diagram}
\end{figure}

This argument shows that a certain pair of variations of MM10 are
equivalent. Thinking about the non-generic projection
corresponding to MM10 in Figure \ref{fig:MM10-projection}, the two
variations are related simply by rotating one strand past an
adjacent one. It's relatively straightforward to see that these
pairs suffice to connect any two variations.

\subsubsection{MM11-15}

\subsubsection*{MM11}
$$\mathfig{0.2}{movie_moves/MM11}$$
This is trivial in either time direction; the complexes involved
only have a single object, and the relevant pairs of cobordisms are isotopic.

\subsubsection*{MM12}
$$\mathfig{0.2}{movie_moves/MM12}$$
We can't use a homotopy isolation argument for MM12, but it's easy
enough to look at all components of the map.

{
\newcommand{\pd}[1]{\mathfig{0.1}{movie_moves/MM12/#1}}
\newcommand{\pds}[1]{\mathfig{0.06}{movie_moves/MM12/#1}}
\newcommand{\birth}{\rotatemathfig{0.035}{90}{cobordisms/cap_bdy_right}}
\newcommand{\death}{\rotatemathfig{0.035}{90}{cobordisms/cap_bdy_left}}
We need to deal with MM12 in two mirror images. In the first mirror
image, there is a positive crossing. Reading down, we have on the
left
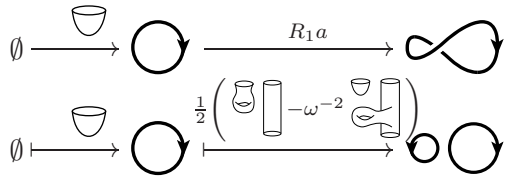
\begin{align*}
\xymatrix@R-5mm@C+3mm{
    \emptyset \ar[r]^{\birth} & \pds{MM12f2} \ar[rr]^{R_1a} & & \pd{MM12f3a} \\
     \emptyset
     \ar@{|->}[r]^{\birth} & %
     \pds{MM12f2} \ar@{|->}[rr]^{\frac{1}{2}\left(\reflectmathfig{0.05}{movie_moves/MM12/cylinder_and_handle_creation} - \omega^{-2} \reflectmathfig{0.05}{movie_moves/MM12/cylinder_with_handle_and_disc_creation}\right)} & & %
     \pd{MM12f3a_trace} \\ %
}
\end{align*}
while on the right we have
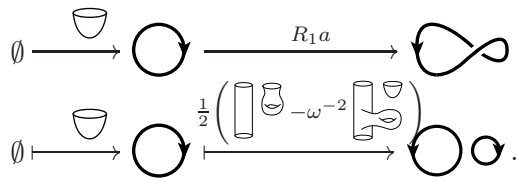
\begin{align*}
\xymatrix@R-5mm@C+3mm{
    \emptyset
    \ar[r]^{\birth} & \pds{MM12f2} \ar[rr]^{R_1a} & & \pd{MM12f3b} \\
     \emptyset \ar@{|->}[r]^{\birth} & %
     \pds{MM12f2} \ar@{|->}[rr]^{\frac{1}{2}\left(\mathfig{0.05}{movie_moves/MM12/cylinder_and_handle_creation} - \omega^{-2} \mathfig{0.05}{movie_moves/MM12/cylinder_with_handle_and_disc_creation}\right)} & & %
     \pd{MM12f3b_trace}. \\ %
}
\end{align*}
Composing, we see the morphisms agree in the disoriented theory,
when $\omega^2 = -1$, but differ by a sign in the unoriented theory.

Reading up, we have on the left
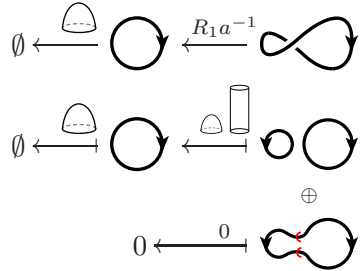
\begin{align*}
\xymatrix@R-5mm{
    \emptyset & \pds{MM12f2} \ar[l]_{\death} & \pd{MM12f3a} \ar[l]_{{R_1a}^{-1}}  \\
     \emptyset & %
     \pds{MM12f2} \ar@{|->}[l]_{\death} & %
     \pd{MM12f3a_trace} \ar@{}[d]|{\directSum} \ar@{|->}[l]_{\reflectmathfig{0.05}{movie_moves/MM12/cylinder_and_disc_annihilation}} \\ %
     & 0 & \pd{MM12f3a_dtrace} \ar@{|->}[l]_0 \\
}
\end{align*}
while on the right we have
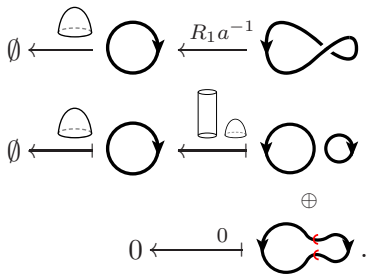
\begin{align*}
\xymatrix@R-5mm{
    \emptyset & \pds{MM12f2} \ar[l]_{\death} & \pd{MM12f3b} \ar[l]_{{R_1a}^{-1}}  \\
     \emptyset & %
     \pds{MM12f2} \ar@{|->}[l]_{\death} & %
     \pd{MM12f3b_trace} \ar@{}[d]|{\directSum} \ar@{|->}[l]_{\mathfig{0.05}{movie_moves/MM12/cylinder_and_disc_annihilation}} \\ %
     & 0 & \pd{MM12f3b_dtrace}. \ar@{|->}[l]_0 \\
}
\end{align*}
}%

These chain maps agree exactly, in both the unoriented and the
disoriented theory.

The mirror image is much the same, although it's the forward in time
maps that agree exactly, and the backwards in time maps that agree
up to a factor of $-\omega^2$.

\subsubsection*{MM13}
{
\newcommand{\pd}[1]{\mathfig{0.1}{movie_moves/MM13/#1}}
\newcommand{\pds}[1]{\mathfig{0.06}{movie_moves/MM13/#1}}
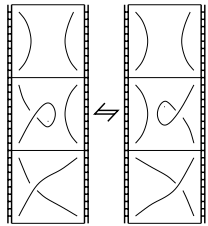
$$\mathfig{0.2}{movie_moves/MM13}$$
This time there are no orientation variations; we can take both
strands in the initial frame to be oriented upwards. We need to
compare the two clips read both up and down, and also consider the
mirror image.

Reading down we have on the left
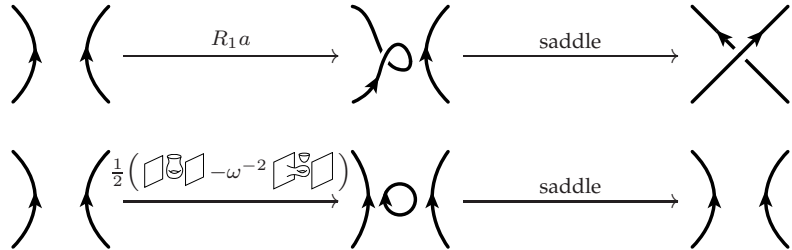
\begin{align*}
\xymatrix@R-5mm@C+20mm{
    \pd{MM13f1} \ar[r]^{R_1a} & \pd{MM13f2a} \ar[r]^{\text{saddle}} & \pd{MM13f3} \\
    \pd{MM13f1} \ar[r]^{\frac{1}{2}\left(\pds{handle_creation} - \omega^{-2} \pds{disc_creation_a}\right)} & \pd{MM13f2a_trace} \ar[r]^{\text{saddle}} & \pd{MM13f1}
}
\end{align*}
and on the right
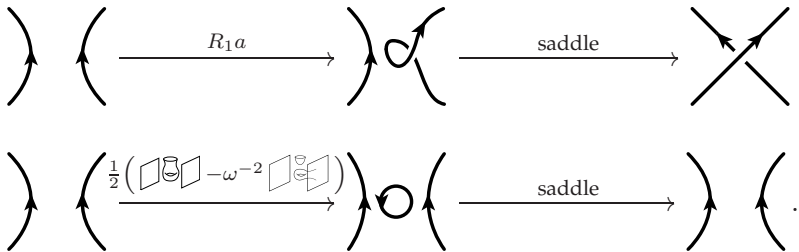
\begin{align*}
\xymatrix@R-5mm@C+20mm{
    \pd{MM13f1} \ar[r]^{R_1a} & \pd{MM13f2b} \ar[r]^{\text{saddle}} & \pd{MM13f3} \\
    \pd{MM13f1} \ar[r]^{\frac{1}{2}\left(\pds{handle_creation} - \omega^{-2} \pds{disc_creation_b}\right)} & \pd{MM13f2b_trace} \ar[r]^{\text{saddle}} & \pd{MM13f1}.
}
\end{align*}
These maps differ by a sign of $-\omega^2$. Reading up, both maps
are the identity on the oriented smoothing, and zero on the
disoriented smoothing, and hence agree on the nose.

In the mirror image, we see the opposite pattern (since the
`interesting' morphism in the R1a and R1b maps appears in
opposite directions). There are no other orientations to deal with.
}

\subsubsection*{MM14}
$$\mathfig{0.2}{movie_moves/MM14}$$
Fixing the orientation of the strand to be from bottom to top, the
loop can either be clockwise or counterclockwise, and lie either
below or above the strand. We'll first deal with the case in which
is loop is oriented counterclockwise, and lies below the strand.

We can't use a homotopy isolation argument for MM12, but it's easy
enough to look
 at all components of the map. Because we're looking
at all components, we actually need to pay attention to the ordering
of crossings; for compatibility with the Reidemeister maps described
in \S \ref{sssec:R2_maps}, we'll number the crossings from the
bottom up, so that the negative crossing comes first.

On either the left or right sides of MM14, we have an R2 map.
Looking at Figures \ref{fig:R2a_maps} and \ref{fig:R2b_maps}, we see
that on both sides we obtain the map
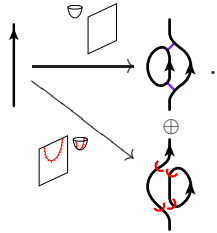
\begin{align*}
\xymatrix@R-7.5mm@C+5mm{
 \mathfig{0.01}{movie_moves/MM14/MM14f1}
    \ar[r]^{\mathfig{0.05}{movie_moves/MM14/curtain_and_disc_creation}}
    \ar[dr]_{\mathfig{0.05}{movie_moves/MM14/curtain_and_disc_creation_disoriented}} &
 \mathfig{0.05}{movie_moves/MM14/MM14f3cx1} \ar@{}[d]|{\directSum} \\
 &
 \mathfig{0.05}{movie_moves/MM14/MM14f3cx2}
}.
\end{align*}

Backwards in time, we obtain different maps. On the left we see
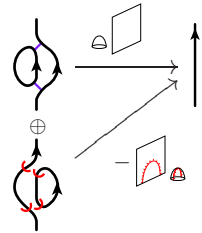
\begin{align*}
\xymatrix@R-7.5mm@C+5mm{
 \mathfig{0.05}{movie_moves/MM14/MM14f3cx1} \ar@{}[d]|{\directSum}
     \ar[r]^{\mathfig{0.05}{movie_moves/MM14/curtain_and_disc_annihilation}}
 &
 \mathfig{0.01}{movie_moves/MM14/MM14f1}
 \\
 \mathfig{0.05}{movie_moves/MM14/MM14f3cx2}
    \ar[ur]_{-\mathfig{0.05}{movie_moves/MM14/curtain_and_disc_annihilation_disoriented}}
 &
}
\end{align*}
and on the right
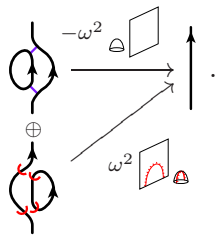
\begin{align*}
\xymatrix@R-7.5mm@C+5mm{
 \mathfig{0.05}{movie_moves/MM14/MM14f3cx1} \ar@{}[d]|{\directSum}
     \ar[r]^{- \omega^2 \mathfig{0.05}{movie_moves/MM14/curtain_and_disc_annihilation}}
 &
 \mathfig{0.01}{movie_moves/MM14/MM14f1}
 \\
 \mathfig{0.05}{movie_moves/MM14/MM14f3cx2}
    \ar[ur]_{\omega^2 \mathfig{0.05}{movie_moves/MM14/curtain_and_disc_annihilation_disoriented}}
 &
}.
\end{align*}

Thus we see that forwards in time the maps agree, but backwards in
time they only agree in the disoriented theory.

Reversing the relative heights of the loop and the strand doesn't
change the calculation; similarly reversing the orientation of one
strand has no effect.

\subsubsection*{MM15}
We now consider both time directions in MM15.
$$\mathfig{0.2}{movie_moves/MM15}$$

We need to deal with 4 variations; assuming the middle strand is
oriented left to right, we can orient the highest strand either to
the left or to the right (forcing the lowest strand to be oriented
oppositely), and we can tuck the middle strand either under or over
the other strands.

We'll start by choosing orientations so the upper two strands are
oriented to the right, and the lowest strand is oriented to the
left, and tuck the middle strand under the others.

Reading down, we have on the left
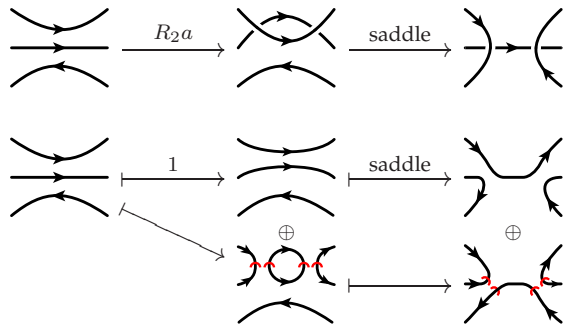
\begin{align*}
\newcommand{\pd}[1]{\mathfig{0.1}{movie_moves/MM15/#1}}
\xymatrix@R-7.5mm@C+5mm{
    \pd{MM15f1} \ar[r]^{R_2a} & \pd{MM15f2a} \ar[r]^{\text{saddle}} & \pd{MM15f3} \\
    \\
     \pd{MM15f1} \ar@{|->}[r]^{1} \ar@{|->}[dr] & %
     \pd{MM15f2a_trace1} \ar@{|->}[r]^{\text{saddle}} \ar@{}[d]|{\directSum} & %
     \pd{MM15f3_trace1} \ar@{}[d]|{\directSum} \\ %
     & \pd{MM15f2a_trace2} \ar@{|->}[r] & \pd{MM15f3_trace2} %
}
\end{align*}
and on the right
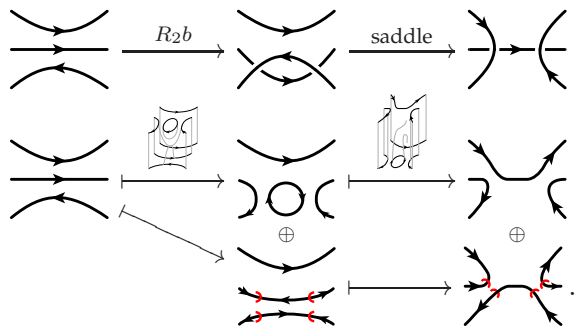
\begin{align*}
\newcommand{\pd}[1]{\mathfig{0.1}{movie_moves/MM15/#1}}
\newcommand{\pa}[1]{\mathfig{0.05}{movie_moves/MM15/#1}}
\xymatrix@R-7.5mm@C+5mm{
    \pd{MM15f1} \ar[r]^{R_2b} & \pd{MM15f2b} \ar[r]^{\text{saddle}} & \pd{MM15f3} \\
     & & \\
     \pd{MM15f1} \ar@{|->}[r]^{\pa{MM15_tuck_map}} \ar@{|->}[dr] & %
     \pd{MM15f2b_trace2} \ar@{|->}[r]^{\pa{elephant}} \ar@{}[d]|{\directSum} & %
     \pd{MM15f3_trace1} \ar@{}[d]|{\directSum} \\ %
     & \pd{MM15f2b_trace1} \ar@{|->}[r] & \pd{MM15f3_trace2}. %
}
\end{align*}

We've left some maps we don't need to know about unlabeled.

Looking only at the component of the maps going to
$\mathfig{0.06}{movie_moves/MM15/MM15f3_trace1}$, we see each side
of the movie move agrees on the nose; both maps are a saddle
involving the lower two strands.

Reading up, we have on the left
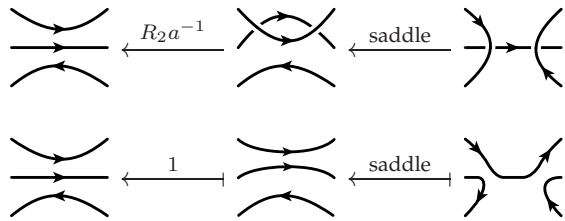
\begin{align*}
\newcommand{\pd}[1]{\mathfig{0.1}{movie_moves/MM15/#1}}
\xymatrix@R-7.5mm@C+5mm{
    \pd{MM15f1}& \pd{MM15f2a} \ar[l]_{{R_2a}^{-1}} & \pd{MM15f3} \ar[l]_{\text{saddle}} \\
    \\
     \pd{MM15f1} & %
     \pd{MM15f2a_trace1} \ar@{|->}[l]_{1} & %
     \pd{MM15f3_trace1} \ar@{|->}[l]_{\text{saddle}} \\ %
}
\end{align*}
and on the right
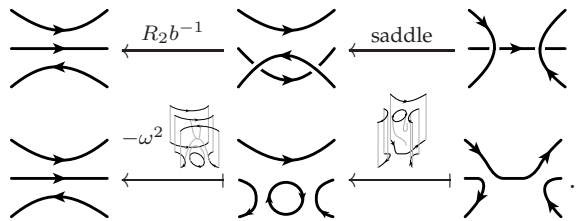
\begin{align*}
\newcommand{\pd}[1]{\mathfig{0.1}{movie_moves/MM15/#1}}
\newcommand{\pa}[1]{\mathfig{0.05}{movie_moves/MM15/#1}}
\xymatrix@R-7.5mm@C+5mm{
    \pd{MM15f1}& \pd{MM15f2b} \ar[l]_{{R_2b}^{-1}} & \pd{MM15f3} \ar[l]_{\text{saddle}} \\
    \\
     \pd{MM15f1} & %
     \pd{MM15f2b_trace2} \ar@{|->}[l]_{-\omega^{2} \pa{MM15_untuck_map}} & %
     \pd{MM15f3_trace1} \ar@{|->}[l]_{\pa{elephant_reverse}}. \\ %
}
\end{align*}
We see that the two movies differ by a sign of $-\omega^2$.

The other variations turn out exactly the same way. Changing the
orientations of the highest and lowest strand has no effect; we
simply interchange $R2a$ and $R2b$ maps throughout. Switching the height
ordering interchanges $R2al$ with $R2ar$, and $R2b+$ with $R2b-$, with no
net effect.

This concludes the proofs of Theorems
\ref{thm:isofunctoriality} and \ref{thm:cobfunctoriality}.

\section{Odds and ends}
\label{sec:oddsandends}
\subsection{Recovering Jacobsson's signs}
Summarizing the results of the above calculations at $\omega = 1$
(i.e. in the original unoriented theory), in Figure
\ref{fig:sign_table}, we see that in most cases we agree with the
signs Jacobsson observed \cite{MR2113903}. There are exceptions,
however (shown highlighted in the tables).

In particular, MM6 (Jacobsson's number 15) does not appear to
exhibit a sign problem in the unoriented theory, and the two mirror
images of MM12 (Jacobsson's number 12) both exhibit a sign problem,
one forwards in time, one backwards. These disagreements coincide
with calculations performed by the first author using Lee's
\cite{MR2173845} variant of Khovanov homology.

\begin{figure}[ht!]
\newcommand{\highlightentry}[1]{\colorbox{yellow}{$#1$}}
\begin{center}
\begin{tabular}{|l|c|c|} \hline
MM & J$\sharp$ & $\pm$ \\
\hline
\textbf{6}           & 15 & \highlightentry{+} \\
7           & 13 & - \\
7 (mirror)  & 13 & + \\
8           & 6  & - \\
8 (mirror)  & 6  & + \\
9           & 14 & + \\
9 (mirror)  & 14 & - \\
10          & 7  & + \\
\hline
\end{tabular}\qquad\qquad
\begin{tabular}{|l|c|c|c|}
\hline
MM & J$\sharp$ & $\downarrow$ & $\uparrow$ \\
\hline
11          & 9  & + & + \\
12          & 11 & - & + \\
\textbf{12 (mirror)} & 11 & + & \highlightentry{-} \\
13          & 12 & - & + \\
13 (mirror) & 12 & + & - \\
14          & 8  & + & - \\
15          & 10 & + & - \\
\hline
\end{tabular}
\end{center}
\caption{The signs observed in the unoriented theory.}
\label{fig:sign_table}
\end{figure}


\subsection{Relationship with the unoriented invariant}
In this section we'll prove that for knots and links (that is,
ignoring tangles and cobordisms), the disoriented and unoriented
invariants are equivalent. We'll write $\KhD{L}$ and $\KhU{L}$ for
the disoriented and unoriented invariants of $L$ respectively.

\begin{thm}
\label{thm:old-theory}%
There's an faithful functor $\Alt: \UnAb_0 \tensor \hZ[\omega] \Into
\DisAb_0$ (here the subscript $0$ denotes the part of the canopolis
with no boundary points), which `alternately orients' each
unoriented diagram. This induces another functor $\Alt:
\KomUnAb_0 \Into \KomDisAb_0$ such that
$$\Alt\left(\KhU{L}\right) \Iso \KhD{L},$$
although this isomorphism isn't canonical.
\end{thm}

\begin{proof}
We've already seen the forgetful map $\DisAb \To \UnAb$, setting
$\omega=1$ and forgetting orientation data. It's relatively easy to
see that this guarantees that we can reconstruct the unoriented
invariant from the disoriented one (for tangles too!). To see the
two invariants are actually equivalent, we'll introduce a new
canopolis of `alternately oriented cobordisms', $\AltAb$, a
subcanopolis of $\DisAb$. We'll construct an isomorphism $\UnAb_0
\tensor \hZ[\omega] \Iso \AltAb_0$, and additionally show that the
invariant of a knot or link (but not a tangle!), which is an
up-to-homotopy complex in $\DisAb_0$, always has a representative in
the subcategory $\AltAb_0$, which coincides with the image of the
unoriented invariant in $\AltAb_0$.

Thus the category $\AltAb_0$ consists of diagrams in the disc
comprised of oriented loops, such that all `outermost' loops are
oriented counterclockwise, and at each successive depth of nesting,
the orientations reverse. This is a subset of the objects of
$\DisAb_0$. The morphisms of $\AltAb_0$ are simply all the morphisms
of $\DisAb_0$ between these objects. In fact, $\AltAb_0$ is the
`boundaryless' part of a full canopolis $\AltAb$ defined in much the
same way.

The isomorphism $\UnAb_0 \tensor \hZ[\omega] \Iso \AltAb_0$ is easy;
simply orient the circles in an object of $\UnAb_0$ in the
prescribed manner, and note that for any cobordism, these
orientations always extend to an honest orientation of the
cobordism. It's an isomorphism because every cobordism in
$\AltAb_0$, which \emph{a priori} might have disorientation seams,
is actually a $\hZ[\omega]$ multiple of a properly oriented
cobordism, by the following Lemma and Corollary.

\begin{lem}
Reversing the fringe of a closed disorientation seam gives a sign of
$-1$.
\end{lem}
\begin{proof}
Use the neck cutting relation parallel to the seam.
\end{proof}

\begin{cor}
\label{cor:alt_dis} If $Y$ is a disoriented surface with all
disorientation seams closed, and with alternately oriented boundary
components, then $Y$ is equal to a multiple of the homeomorphic
oriented surface.
\end{cor}
\begin{proof}
By applying fringe moves, we can assume that the disorientation seam
is connected on each connected component of $Y$. (If necessary,
reverse fringe directions using the previous lemma.) The assumption
about boundary orientations now implies that the seam is
null-homologous, and so can be removed via further fringe moves.
\end{proof}

We next discover how to push a link complex $\Kh{L}$ in $\DisAb_0$
down into the subcategory $\AltAb_0$.

We begin with a quick statement about the disorientations that can
appear on a circle.

\begin{lem}
\label{lem:disorientation_number}%
Define the `disorientation number' of a `disoriented circle' to be
the number of counterclockwise facing disorientation marks minus the
number of clockwise facing disorientation marks. (See Figure
\ref{fig:disoriented_circle}.) Then two disoriented circles $C_1$
and $C_2$ are isomorphic in $\DisAb$ exactly if their disorientation
numbers agree.
\end{lem}

\begin{figure}[!ht]
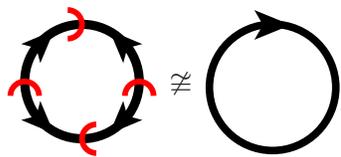

$$\mathfig{0.15}{knot_pieces/disoriented_circle} \ncong
\mathfig{0.135}{knot_pieces/clockwise_circle}$$
\caption{A disoriented circle with disorientation number $+2$ is not isomorphic in $\DisAb$ to an oriented circle.}%
\label{fig:disoriented_circle}
\end{figure}

We want to show that every circle appearing in an object of $\Kh{L}$
has disorientation number $0$. This is a conservation argument; near
each disorientation mark on the circle, there used to be a crossing
in $L$, either just inside or outside the circle. Whether the
disorientation mark faces counterclockwise or clockwise records
whether the two strands in the crossing were oriented `inwards' or
`outwards' across the circle. Since the original link must cross any
given circle a total of $0$ times, the signed count of
disorientation marks is $0$ as well. This shows that every object
appearing $\Kh{L}$ is actually isomorphic to the corresponding
object appearing in $\KhU{L}$ (this is Lemma
\ref{lem:disorientation_number}). In fact, the only choice in this
isomorphism is a multiple of $\pm1$ or $\pm\omega$. 

Thus we take the link complex $\Kh{L}$, and replace every disoriented
circle with the appropriately oriented circle. The complex now lies
entirely within the subcategory $\AltAb_0$. This complex agrees with the
unoriented link complex, thought of as living in $\AltAb_0$, except for
the fact that each morphism may be off by a unit, simply because the
underlying surfaces for each morphism are the same, and by Corollary
\ref{cor:alt_dis} above the morphisms are unit multiples of each other.

A little combinatorial lemma about sprinkling units in a complex
gets us to the desired result.

\begin{lem}[Sprinkling units]
Suppose we have two anticommutative cubes, with identical objects,
such that corresponding morphisms only ever differ by a unit.
Further suppose that the composition of any two `edges' of the cube
is nonzero. Then the two cubes are isomorphic, via a map which just
multiplies each object in the cube by some unit.
\end{lem}

\begin{rem}
The hypothesis that the composition of any two composable maps in
the cube is nonzero certainly holds in the case we're interested in.
The complex associated to a knot or link has as morphisms pairs of
pants and cylinders, and it's easy to see that any composition of a
pair of pairs of pants is nonzero.
\end{rem}

\begin{rem}
Something like this lemma is used in \cite{MR2124557} in describing
a categorification of the colored Jones polynomial, without the need
for functoriality. Note also that our construction of a properly
functorial version of Khovanov homology should make a more direct
construction of a categorification of the colored Jones polynomial
possible, and allow the possibility of this categorification itself
being functorial. See p. 20 of \cite{MR2124557}.
\end{rem}

\begin{proof}
An easy induction on the dimension of the cube. For one dimensional
cubes, the result is trivial. For any cube, by induction we can
choose an isomorphism $\phi_t$ between the top layers of the cubes,
and another $\phi_b$ between the bottom layers of the cubes. Now we
need to tweak the top layer isomorphism, so together the
isomorphisms give an isomorphism on the entire cube. Consider the
`highest' vertical differential $d_v$, between the initial objects
in the top and bottom layers, and define a unit $\epsilon$ by $d_v
\phi_t = \epsilon \phi_b d_v$. Now replace the isomorphism $\phi_t$
with $\epsilon \phi_t$. We now just need to check that our
isomorphism $\phi$ commutes with every vertical differential. Thus
consider a square of differentials in one cube,
$$\xymatrix{
    \bullet \ar[d]^{d^1_l} \ar[r]^{d^1_t} & \bullet \ar[d]^{d^1_r} \\
    \bullet                \ar[r]^{d^1_b} & \bullet
}$$ with $d^1_t$ a differential in the top layer, and $d^1_b$ a
differential in the bottom layer. There's a corresponding square of
differentials in the other cube, with differentials $d^2_t, d^2_b,
d^2_r$ and  $d^2_l$. By our construction $\phi d^1_t = d^2_t \phi$,
and $\phi d^1_b = d^2_b \phi$, and we'll assume further $\phi d^1_l
= d^2_l \phi$ (we're going to apply this piece of the argument to
every such square, starting with $d_l = d_v$, the `highest' vertical
differential described above). Now we know $\phi d^1_r = \zeta d^2_r
\phi$ for some unit $\zeta$; we just need to show $\zeta = 1$. We
then deduce the following equations
\begin{align*}
\phi d^1_r d^1_t & = \zeta d^2_r \phi d^1_t \\
                 & = \zeta d^2_r d^2_t \phi \\
\phi d^1_b d^1_l & = \zeta d^2_b d^2_l \phi \\
                 & = \zeta \phi d^1_b d^1_l
\end{align*}
and, making use of the hypothesis that the composition $d^1_b d^1_l$
is nonzero,  conclude that $\zeta$ is indeed $1$.
\end{proof}

That concludes the proof of Theorem \ref{thm:old-theory}.
\end{proof}

\subsection{Sliding a handle past a crossing}%
\label{ssec:handle-sliding}
In this section we give an example calculation in the new setup,
illustrating an interesting difference with the unoriented construction of
Khovanov homology.

Consider the following two cobordisms from
$\mathfig{0.06}{handle_crossing_lemma/pos_crossing}$ to itself:
the first is the identity except for a handle attached to the
over-sheet to the right of the crossing, and the second is the same except that
the handle is attached to the over-sheet to the left of the crossing. We'll
denote these schematically by
$F=\mathfig{0.06}{handle_crossing_lemma/morphism1}$ and
$G=\mathfig{0.06}{handle_crossing_lemma/morphism2}$.

\begin{prop}
\label{prop:handle_crossing}%
$F \htpy G$ are homotopic as maps in $\KomDisAb$, whereas in $\KomUnAb$,
we have $F \htpy -G$ instead.
\end{prop}

\begin{proof}
Note that these cobordisms are clearly isotopic, and so the
functoriality result above gives us an automatic proof. (Exercise: figure
out the sequence of movie moves relating them!) However, we will
construct an explicit homotopy: the arrow marked $h$ in Figure \ref{fig:handle
crossing diagram}.

\begin{figure}[ht]
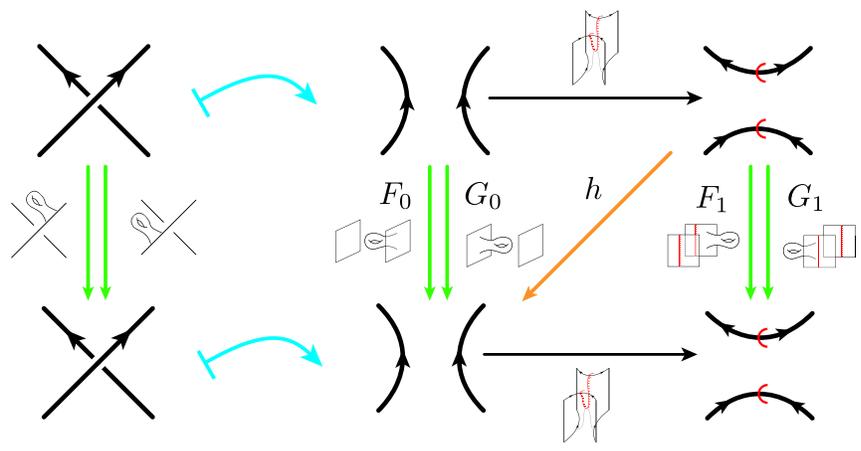

$$\mathfig{0.85}{handle_crossing_lemma/chain_map_diagram}$$
\caption{Chain maps and homotopy for the Proposition \ref{prop:handle_crossing}.} \label{fig:handle
crossing diagram}
\end{figure}

We need a homotopy $h$ such that $hd+dh =
F_{i}-G_{i}$, and propose $$h=2
\omega^{-1}\mathfig{0.1}{handle_crossing_lemma/saddle_seam_down}.$$

At height zero we then have
\begin{align*}
hd + dh & = 2 \omega^{-1} \mathfig{0.18}{handle_crossing_lemma/handle_crossing_3} \displaybreak[2] \\
    & = \omega^{-1} \left( \mathfig{0.18}{handle_crossing_lemma/handle_crossing_4} + \mathfig{0.18}{handle_crossing_lemma/handle_crossing_5} \right) \displaybreak[1] \\
    & = \omega^{-1} \left( \omega \mathfig{0.18}{cobordisms/two_sheets_handle_right} + \omega^{-1} \mathfig{0.18}{cobordisms/two_sheets_handle_left} \right) \displaybreak[2] \\
    & = \mathfig{0.18}{cobordisms/two_sheets_handle_right} + \omega^{-2} \mathfig{0.18}{cobordisms/two_sheets_handle_left},
\end{align*}

which is just $F_{0} - G_{0}$ when we set $\omega^{2}=-1$.
There's a similar computation at height 1:

\begin{align*}
hd + dh & = 2 \omega^{-1}\mathfig{0.18}{handle_crossing_lemma/handle_crossing_6} \displaybreak[2] \\
    & =  2 \mathfig{0.18}{handle_crossing_lemma/handle_crossing_7} \displaybreak[1] \\
    & =  \mathfig{0.18}{handle_crossing_lemma/handle_crossing_8} + \mathfig{0.18}{handle_crossing_lemma/handle_crossing_9} \displaybreak[2] \\
    & =  \mathfig{0.18}{handle_crossing_lemma/handle_crossing_8} + \omega^{-2} \mathfig{0.18}{handle_crossing_lemma/handle_crossing_10}.
\end{align*}

Again, setting $\omega^{2}=-1$ makes the last line $F_{1} -
G_{1}$, which gives the result.
\end{proof}

\subsection{Confusions}
\label{ssec:confusions}%
In this final section, we'll describe a
defect in the discussion so far, and say a little about a proposal
to fix it.

The construction we've proposed so far is a functor from the
category of oriented tangles, $\OrTang$, into the category of
complexes of disoriented flat tangles $\KomDisAb$. In particular, it
only gives maps for oriented cobordisms between oriented links. This
isn't really ideal; the old unoriented theory gave maps for
nonorientable cobordisms. For example, while a
M\"{o}bius band with positive $\frac{3}{2}$ twists
provides a generator of the Khovanov invariant of
the trefoil in the old theory, our construction doesn't know what to
do with nonorientable surfaces.

Thus we'd like to extend the theory to a functor from $\DisTang$,
the category of disoriented tangles. On the level of objects, this
is no problem; simply map disorientations to disorientations.
Unfortunately, there is now an additional Reidemeister move, namely
`sliding a disorientation through a crossing', which we'll name a
`vertigo', for which we need to provide an isomorphism between the
corresponding complexes. Further, we'd need to check additional
movie moves, relating this new Reidemeister move to the original
three.\footnote{There's actually a big incentive for this extension;
it turns out that all the different oriented versions of the usual 15
movie moves become equivalent modulo these extra movie moves
involving disorientations. This was actually our original motivation
for introducing confusions.}

\begin{figure}
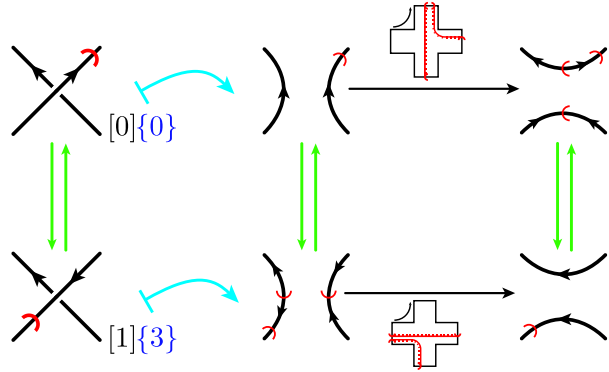

$$\mathfig{0.6}{vertigos/v_oxl_no_chain_maps}$$
\caption{A hypothetical isomorphism of complexes implementing a
particular case of the `vertigo' move.}%
\label{fig:vertigo_no_chain_maps}%
\end{figure}

However, it's easy to see that it just isn't possible to produce a
homotopy equivalence between the corresponding complexes in
$\KomDisAb$. To begin, such a homotopy equivalence would have to be
an isomorphism; using Lemma \ref{lem:loopless}, we see no homotopies
are possible in the complex for a single crossing, regardless of any
additional disorientations. Such an isomorphism would presumably be
of the form in Figure \ref{fig:vertigo_no_chain_maps}.

In particular, we'd need an isomorphism in $\DisAb$ reversing the
direction of a disorientation mark on a disoriented strand.

Such an isomorphism, which we'll dub a `confusion', would
necessarily be a troublesome thing; if the confusion were simply to
be some `local' structure on a surface, which I'll draw here as a
box labeled by $c$ (or a box labeled by $c^{-1}$ for its inverse),
we could perform the
calculation%
\begin{align}
\label{eq:confusion_contradiction}%
\omega
 = \mathfig{0.135}{disorientation_relations/inward_circle}
 = \mathfig{0.15}{confusions/contradiction2}
 = \mathfig{0.15}{confusions/contradiction1}
 = \mathfig{0.15}{disorientation_relations/outward_circle}
 = \omega^{-1}
\end{align}%
producing a contradiction with the requirement that $\omega^2 = -1$.

The way out of this seems to be to make the confusion a spinorial
object, so an extra sign gets introduced as we drag the confusion
around the circle, in the third equality in Equation
\ref{eq:confusion_contradiction}.

At this point it seems appropriate to apologise for having talked about a
particular diagrammatic model for such `spinorial confusions' at various
conferences, but to be omitting the details in this paper. We still
intend to write these details down!

We'll briefly list the improvements to the theory we anticipate
being able to make, after the introduction of confusions.
\begin{itemize}
\item Connecting the category $\DisAb$; in particular,
all disoriented circles would be isomorphic.
\item Extending the invariant to disoriented tangles, and
disoriented cobordisms between them.
\item Using the categorified Kauffman trick, to more easily describe
the Reidemeister 3 chain map.
\item After checking additional movie moves involving vertigos,
being able to reduce the computations required in \S
\ref{ssec:movie_moves}, by taking advantage of the fact that all
oriented versions of each oriented movie move become equivalent
module disoriented movie moves.
\end{itemize}

\appendix
\section{Boring technical details}
Throughout this appendix, we'll at times just write a `bullet', $\bullet$, for a matrix entry which we don't need to care about.

\subsection{Gaussian elimination}
\begin{lem}[Gaussian elimination for complexes]%
\label{lem:gaussian}%
Consider the complex
\begin{equation}
\label{eq:complex}
 \xymatrix@C+55pt@R+20pt{
    A                     \ar[r]^{\psmallmatrix{\bullet \\ \alpha}}              &
    \directSumStack{B}{C} \ar[r]^{\psmallmatrix{\varphi & \lambda \\ \mu & \nu}} &
    \directSumStack{D}{E} \ar[r]^{\psmallmatrix{\bullet & \epsilon}}             &
    F
 }
\end{equation}
in any additive category, where $\varphi: B \IsoTo D$ is an
isomorphism, and all other morphisms are arbitrary (subject to
$d^2=0$, of course). Then there is a homotopy equivalence with a
much simpler complex, `stripping off' $\varphi$.

\begin{equation*}
 \xymatrix@C+55pt@R+20pt{
    A                     \ar[r]^{\psmallmatrix{\bullet \\ \alpha}}              \ar@{<->}[d]^{\psmallmatrix{1}} &
    \directSumStack{B}{C} \ar[r]^{\psmallmatrix{\varphi & \lambda \\ \mu & \nu}} \ar@<-0.5ex>[d]_{\psmallmatrix{0 & 1}}  &
    \directSumStack{D}{E} \ar[r]^{\psmallmatrix{\bullet & \epsilon}}             \ar@<-0.5ex>[d]_{\psmallmatrix{-\mu \varphi^{-1} & 1}} &
    F                                                                            \ar@{<->}[d]^{\psmallmatrix{1}} \\
    A \ar[r]^{\psmallmatrix{\alpha}}                          &
    C \ar[r]^{\psmallmatrix{\nu - \mu \varphi^{-1} \lambda}} \ar@<-0.5ex>[u]_{\psmallmatrix{-\varphi^{-1} \lambda \\ 1}} &
    E \ar[r]^{\psmallmatrix{\epsilon}}                       \ar@<-0.5ex>[u]_{\psmallmatrix{0 \\ 1}} &
    F
 }
\end{equation*}
\end{lem}
\begin{rem}
Gaussian elimination is a strong deformation retract. In fact, it
preserves the simple homotopy type of the complex.
\end{rem}
\begin{proof}
This is simply Lemma 4.2 in \cite{math.GT/0606318} (see also Figure 2 there),
this time explicitly keeping track of the chain maps.
\end{proof}

We'll also state here the result of applying Gaussian elimination
twice, on two adjacent (but non-composable) isomorphisms. Having
these chain homotopy equivalences handy will tidy up the
calculations for the Reidemeister 2 and 3 chain maps.

\begin{lem}[Double Gaussian elimination]%
\label{lem:double_gaussian}%
When $\psi$ and $\varphi$ are isomorphisms, there's a homotopy
equivalence of complexes:
\begin{equation*}
 \xymatrix@C+45pt@R+20pt{
    A                                 \ar[r]^{\psmallmatrix{\bullet \\ \alpha}}                                    \ar@{<->}[d]^{\psmallmatrix{1}} &
    \directSumStack{B}{C}             \ar[r]^{\psmallmatrix{\psi & \beta \\ \bullet & \bullet \\ \gamma & \delta}} \ar@<-0.5ex>[d]_{\psmallmatrix{0 & 1}}  &
    \directSumStackThree{D_1}{D_2}{E} \ar[r]^{\psmallmatrix{\bullet & \varphi & \lambda \\ \bullet & \mu & \nu}}   \ar@<-0.5ex>[d]_{\psmallmatrix{-\gamma \psi^{-1} & 0 & 1}} &
    \directSumStack{F}{G}             \ar[r]^{\psmallmatrix{\bullet & \eta}}                                       \ar@<-0.5ex>[d]_<(0.2){\psmallmatrix{-\mu \varphi^{-1} & 1}} &
    H                                                                                                              \ar@{<->}[d]^{\psmallmatrix{1}} \\
    A \ar[r]^{\psmallmatrix{\alpha}}                          &
    C \ar[r]^{\psmallmatrix{\delta - \gamma \psi^{-1} \beta}} \ar@<-0.5ex>[u]_<(0.25){\psmallmatrix{-\psi^{-1} \beta \\ 1}} &
    E \ar[r]^{\psmallmatrix{\nu - \mu \varphi^{-1} \lambda}}  \ar@<-0.5ex>[u]_<(0.3){\psmallmatrix{0 \\ - \varphi^{-1} \lambda \\ 1}} &
    G \ar[r]^{\psmallmatrix{\eta}}                            \ar@<-0.5ex>[u]_<(0.35){\psmallmatrix{0 \\ 1}} &
    H
 }
\end{equation*}
\end{lem}
\begin{proof}
Apply Lemma \ref{lem:gaussian}, killing off the isomorphism $\psi$.
Notice that the isomorphism $\varphi$ survives unchanged in the
resulting complex, and apply the lemma again.
\end{proof}
\begin{rem}
Convince yourself that it doesn't matter in which order we cancel
the isomorphisms!
\end{rem}

\subsection{Calculations of Reidemeister chain maps}
\label{ssec:reidemeister_appendix}%
We can now go through the constructions of the Reidemeister chain
maps.

\begin{lem}
The chain maps displayed in Figures \ref{fig:R1a_maps} and
\ref{fig:R1b_maps} are homotopy equivalences.
\end{lem}
\begin{proof}
We'll just do the R1a move; the R1b is much the same.

The complex associated to
$\mathfig{0.06}{knot_pieces/positive_curl_right}$ is
$$
\xymatrix{
 q \mathfig{0.06}{knot_pieces/strand_with_circle} \ar[r]^d & q^2 \mathfig{0.06}{knot_pieces/disoriented_strand_bending_right}
}
$$
with $d$ simply the disoriented saddle. Delooping at homological
height $1$, and cancelling the disorientations at height $2$, using
the isomorphisms
\begin{align*}
 \zeta_1 & = \begin{pmatrix} \frac{1}{2}\mathfig{0.075}{cobordisms/curtain_and_handle_annihilation} \\ \mathfig{0.065}{cobordisms/curtain_and_disc_annihilation} \end{pmatrix} &
 \zeta_2 & = \omega^{-1} \mathfig{0.05}{cobordisms/fringe_semicircle_down} \\
\intertext{with inverses}
 \zeta_1^{-1} & = \begin{pmatrix} \mathfig{0.065}{cobordisms/curtain_and_disc_creation} & \frac{1}{2}\mathfig{0.075}{cobordisms/curtain_and_handle_creation} \end{pmatrix} &
 \zeta_2^{-1} & = \mathfig{0.045}{cobordisms/fringe_semicircle}, \\
\end{align*}
we obtain the complex
$$
\xymatrix@C+35mm{
 \directSumStack{q^2 \mathfig{0.06}{knot_pieces/strand_bending_right}}{\phantom{q^2} \mathfig{0.06}{knot_pieces/strand_bending_right}}
    \ar[r]^{\begin{pmatrix} \varphi = \Id & \lambda = \frac{\omega^{-2}}{2} \mathfig{0.05}{cobordisms/sheet_with_handle} \end{pmatrix}} &
 q^2 \mathfig{0.06}{knot_pieces/strand_bending_right}
}.
$$
The differential here is the composition $\zeta_2 d \zeta_1^{-1}$.
Stripping off the isomorphism $\varphi$, according to Lemma
\ref{lem:gaussian}, we see that the complex is homotopy equivalent
to the desired complex: a single strand. The `simplifying' homotopy
equivalence is
\begin{align*}
 s_1 & = \begin{pmatrix}0 & \Id\end{pmatrix} \compose \zeta_1 = \mathfig{0.065}{cobordisms/curtain_and_disc_annihilation} & s_2 & = 0 \\
\intertext{with inverse}%
 s_1^{-1} & = \zeta_1^{-1} \compose \begin{pmatrix}-\varphi^{-1} \lambda \\ \Id \end{pmatrix} = \frac{1}{2}\left(\mathfig{0.065}{cobordisms/curtain_and_handle_creation} - \omega^{-2} \mathfig{0.065}{cobordisms/curtain_handle_and_disc_creation}\right) & s_2^{-1} & = 0
\end{align*}
as claimed.
\end{proof}

\begin{lem}
The chain maps displayed in Figures \ref{fig:R2a_maps} and
\ref{fig:R2b_maps} are homotopy equivalences.
\end{lem}
\begin{proof}
We'll deal with the R2a move first.

The complex associated to $\mathfig{0.1}{knot_pieces/R2al}$ is
$$
\xymatrix{
 q^{-1} \mathfig{0.1}{reidemeister_maps/R2a_complex/R2a_cx1} \ar[r]^{d_{-1}} &
 \directSumStack{\mathfig{0.1}{reidemeister_maps/R2a_complex/R2a_cx2a}}{\mathfig{0.1}{reidemeister_maps/R2a_complex/R2a_cx2b}} \ar[r]^{d_0} &
 q \mathfig{0.1}{reidemeister_maps/R2a_complex/R2a_cx3}
}
$$
with differentials
\begin{align*}
 d_{-1} & =
 \begin{pmatrix}
    \mathfig{0.05}{reidemeister_maps/R2a_complex/R2a_cx1_da}
    \\ \mathfig{0.05}{reidemeister_maps/R2a_complex/R2a_cx1_db}
 \end{pmatrix} &
 d_0    & =
 \begin{pmatrix}
    \mathfig{0.05}{reidemeister_maps/R2a_complex/R2a_cx2a_d}
    & - \mathfig{0.05}{reidemeister_maps/R2a_complex/R2a_cx2b_d}
 \end{pmatrix}
\end{align*}

Applying the delooping isomorphism $\psmallmatrix{\frac{1}{2\omega}
\mathfig{0.05}{reidemeister_maps/R2a_complex/R2a_cx2b_delooping_handle}
\\ \mathfig{0.05}{reidemeister_maps/R2a_complex/R2a_cx2b_delooping}}$ (which has inverse
$\psmallmatrix{
\mathfig{0.05}{reidemeister_maps/R2a_complex/R2a_cx2b_looping} &
\frac{1}{2\omega}\mathfig{0.05}{reidemeister_maps/R2a_complex/R2a_cx2b_looping_handle}}$)
to the direct summand with a loop, we obtain the complex
$$
\xymatrix{
 q^{-1} \mathfig{0.1}{reidemeister_maps/R2a_complex/R2a_cx1} \ar[r]^{d_{-1}} &
 \directSumStackThree
    {       \mathfig{0.1}{reidemeister_maps/R2a_complex/R2a_cx2a}}
    {q      \mathfig{0.1}{reidemeister_maps/R2a_complex/R2a_cx2b_delooped}}
    {q^{-1} \mathfig{0.1}{reidemeister_maps/R2a_complex/R2a_cx2b_delooped}}
   \ar[r]^{d_0} &
 q \mathfig{0.1}{reidemeister_maps/R2a_complex/R2a_cx3}
}
$$
where
\begin{align*}
 d_{-1} & =
 \begin{pmatrix}\gamma = \mathfig{0.05}{reidemeister_maps/R2a_complex/R2a_cx1_da}
 \\ \bullet
 \\ \psi = \mathfig{0.05}{reidemeister_maps/R2a_complex/R2a_cx1}\end{pmatrix} &
 d_0    & =
 \begin{pmatrix}\lambda = \mathfig{0.05}{reidemeister_maps/R2a_complex/R2a_cx2a_d}
 & \varphi = - \mathfig{0.05}{reidemeister_maps/R2a_complex/R2a_cx2b_delooped}
 & \bullet\end{pmatrix}.
\end{align*}

Here we've named the entries of the differentials in the manner
indicated in Lemma \ref{lem:double_gaussian}. Applying that lemma
gives us chain equivalences with the desired one object complex. The
chain equivalences we're after are compositions of the chain
equivalences from Lemma \ref{lem:double_gaussian} with the delooping
isomorphism or its inverse.

Thus the R2a `untuck' chain map is
\begin{equation*}
 \begin{pmatrix} 1 & 0 & - \gamma\psi^{-1} \end{pmatrix} \compose
 \begin{pmatrix}
    1 & 0 \\
    0 & \frac{1}{2\omega} \mathfig{0.05}{reidemeister_maps/R2a_complex/R2a_cx2b_delooping_handle} \\
    0 & \mathfig{0.05}{reidemeister_maps/R2a_complex/R2a_cx2b_delooping}
 \end{pmatrix}
 =
 \begin{pmatrix}
    1 & - \mathfig{0.05}{reidemeister_maps/R2a_complex/R2a_cx1_da} \compose \mathfig{0.05}{reidemeister_maps/R2a_complex/R2a_cx2b_delooping}
 \end{pmatrix}
\end{equation*}
as claimed, and the `tuck' map is
\begin{equation*}
 \begin{pmatrix}
    1 & 0 & 0 \\
    0 & \mathfig{0.05}{reidemeister_maps/R2a_complex/R2a_cx2b_looping} & \frac{1}{2\omega} \mathfig{0.05}{reidemeister_maps/R2a_complex/R2a_cx2b_looping_handle}
 \end{pmatrix}  \compose
 \begin{pmatrix} 1 \\ - \varphi^{-1}\lambda \\ 0\end{pmatrix}
 =
 \begin{pmatrix}
    1 \\ \mathfig{0.05}{reidemeister_maps/R2a_complex/R2a_cx2b_looping} \compose \mathfig{0.05}{reidemeister_maps/R2a_complex/R2a_cx2a_d}
 \end{pmatrix}
\end{equation*}

Now the R2b move, in much the same way. The complex associated to
$\mathfig{0.1}{knot_pieces/R2b-}$ is
$$
\xymatrix{
 q^{-1} \mathfig{0.1}{reidemeister_maps/R2b_complex/R2b_cx1} \ar[r]^{d_{-1}} &
 \directSumStack{\mathfig{0.1}{reidemeister_maps/R2b_complex/R2b_cx2a}}{\mathfig{0.1}{reidemeister_maps/R2b_complex/R2b_cx2b}} \ar[r]^{d_0} &
 q \mathfig{0.1}{reidemeister_maps/R2b_complex/R2b_cx3}
}
$$
with differentials
\begin{align*}
 d_{-1} & =
 \begin{pmatrix}
    \mathfig{0.05}{reidemeister_maps/R2b_complex/R2b_cx1_da}
    \\ \mathfig{0.05}{reidemeister_maps/R2b_complex/R2b_cx1_db}
 \end{pmatrix} &
 d_0    & =
 \begin{pmatrix}
    - \mathfig{0.05}{reidemeister_maps/R2b_complex/R2b_cx2a_d}
    & \mathfig{0.05}{reidemeister_maps/R2b_complex/R2b_cx2b_d}
 \end{pmatrix}
\end{align*}
This time instead of just delooping, we'll also cancel the obvious
pairs of disorientation marks. The isomorphism we'll use is
\begin{align*}
\zeta^\bullet & =
\begin{pmatrix}
 \omega^{-1}
 \mathfig{0.025}{reidemeister_maps/R2b_complex/curtain_fringe_down}
 \mathfig{0.025}{reidemeister_maps/R2b_complex/curtain}
\end{pmatrix},
\begin{pmatrix}
 \omega^{-2} \mathfig{0.05}{reidemeister_maps/R2b_untuck_map_sheets} & 0 \\
 0 & \frac{1}{2}\mathfig{0.05}{reidemeister_maps/R2b_complex/R2b_cx2b_delooping_handle} \\
 0 & \mathfig{0.05}{reidemeister_maps/R2b_complex/R2b_cx2b_delooping} \\
\end{pmatrix},
\begin{pmatrix}
 \omega^{-1}
 \mathfig{0.025}{reidemeister_maps/R2b_complex/curtain}
 \mathfig{0.025}{reidemeister_maps/R2b_complex/curtain_fringe_down}
\end{pmatrix}, \\
\intertext{with inverses}%
(\zeta^{-1})^\bullet & =
\begin{pmatrix}
 \mathfig{0.025}{reidemeister_maps/R2b_complex/curtain_fringe_up}
 \mathfig{0.025}{reidemeister_maps/R2b_complex/curtain}
\end{pmatrix},
\begin{pmatrix}
 \mathfig{0.05}{reidemeister_maps/R2b_tuck_map_sheets} & 0 & 0 \\
 0 & \mathfig{0.05}{reidemeister_maps/R2b_complex/R2b_cx2b_looping} & \frac{1}{2}\mathfig{0.05}{reidemeister_maps/R2b_complex/R2b_cx2b_looping_handle}
\end{pmatrix},
\begin{pmatrix}
 \mathfig{0.025}{reidemeister_maps/R2b_complex/curtain}
 \mathfig{0.025}{reidemeister_maps/R2b_complex/curtain_fringe_up}
\end{pmatrix}.
\end{align*}

We obtain the complex
$$
\xymatrix{
 q^{-1} \mathfig{0.1}{reidemeister_maps/R2b_complex/R2b_cx1_oriented} \ar[r]^{d_{-1}} &
 \directSumStackThree
    {       \mathfig{0.1}{reidemeister_maps/R2b_complex/R2b_cx_simple}}
    {q      \mathfig{0.1}{reidemeister_maps/R2b_complex/R2b_cx2b_delooped}}
    {q^{-1} \mathfig{0.1}{reidemeister_maps/R2b_complex/R2b_cx2b_delooped}}
   \ar[r]^{d_0} &
 q \mathfig{0.1}{reidemeister_maps/R2b_complex/R2b_cx3_oriented}
}
$$
where
\begin{align*}
 d_{-1} & =
 \begin{pmatrix}\gamma = \omega^{-1} \mathfig{0.05}{reidemeister_maps/R2b_complex/R2b_cx1_oriented_d}
 \\ \bullet
 \\ \psi = \omega \mathfig{0.05}{reidemeister_maps/R2b_complex/R2b_cx1_oriented}\end{pmatrix} &
 d_0    & =
 \begin{pmatrix}\lambda = - \mathfig{0.05}{reidemeister_maps/R2b_complex/R2b_cx2a_d_oriented}
 & \varphi = \mathfig{0.05}{reidemeister_maps/R2b_complex/R2b_cx2b_delooped}
 & \bullet\end{pmatrix}.
\end{align*}

Thus the R2b `untuck' chain map is
\begin{equation*}
 \begin{pmatrix} 1 & 0 & - \gamma\psi^{-1} \end{pmatrix} \compose
 \begin{pmatrix}
    \omega^{-2} \mathfig{0.05}{reidemeister_maps/R2b_untuck_map_sheets} & 0 \\
    0 & \frac{1}{2} \mathfig{0.05}{reidemeister_maps/R2b_complex/R2b_cx2b_delooping_handle} \\
    0 & \mathfig{0.05}{reidemeister_maps/R2b_complex/R2b_cx2b_delooping}
 \end{pmatrix}
 =
 \begin{pmatrix}
    \omega^2 \mathfig{0.05}{reidemeister_maps/R2b_untuck_map_sheets} & - \omega^2 \mathfig{0.05}{reidemeister_maps/R2b_complex/R2b_cx1_oriented_d} \compose \mathfig{0.05}{reidemeister_maps/R2b_complex/R2b_cx2b_delooping}
 \end{pmatrix}
\end{equation*}
as claimed, and the `tuck' map is
\begin{equation*}
 \begin{pmatrix}
    \mathfig{0.05}{reidemeister_maps/R2b_tuck_map_sheets} & 0 & 0 \\
    0 & \mathfig{0.05}{reidemeister_maps/R2b_complex/R2b_cx2b_looping} & \frac{1}{2} \mathfig{0.05}{reidemeister_maps/R2b_complex/R2b_cx2b_looping_handle}
 \end{pmatrix}  \compose
 \begin{pmatrix} 1 \\ - \varphi^{-1}\lambda \\ 0\end{pmatrix}
 =
 \begin{pmatrix}
    \mathfig{0.05}{reidemeister_maps/R2b_tuck_map_sheets} \\ \mathfig{0.05}{reidemeister_maps/R2b_complex/R2b_cx2b_looping} \compose \mathfig{0.05}{reidemeister_maps/R2b_complex/R2b_cx2a_d_oriented}
 \end{pmatrix}
\end{equation*}
\end{proof}

\begin{proof}[Proof of Proposition \ref{prop:R3}]
Finally, we'll construct
explicit chain maps for the third Reidemeister move.

The complex associated to the left side is
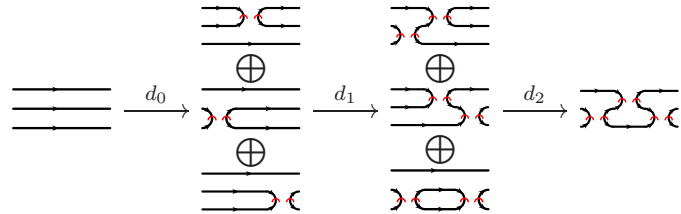
\newcommand{\Ral}[1]{\mathfig{0.1}{reidemeister_maps/R3_complex/R3al_#1}}
\newcommand{\Rar}[1]{\mathfig{0.1}{reidemeister_maps/R3_complex/R3ar_#1}}
\begin{equation*}
 \xymatrix{
    \Ral{000} \ar[r]^{d_0} &
    \directSumStackThree{\Ral{100}}{\Ral{010}}{\Ral{001}} \ar[r]^{d_1} &
    \directSumStackThree{\Ral{110}}{\Ral{101}}{\Ral{011}} \ar[r]^{d_2} &
    \Ral{111}
 }
\end{equation*}
with differentials
\begin{align*}
 d_0 & = \begin{pmatrix}s_1 \\ s_2 \\ s_3 \end{pmatrix} \\
 d_1 & = \begin{pmatrix}s_2 & -s_1 & 0 \\ s_3 & 0 & -s_1 \\ 0 & s_3 & -s_2\end{pmatrix} \\
 d_2 & = \begin{pmatrix}s_3 & -s_2 & s_1\end{pmatrix}.
\end{align*}

We now need to simplify the complex; first delooping the last object
at height two, and cancelling pairs of disorientations at height
three using the isomorphisms
\begin{align*}
 \zeta_{l2}      & = \begin{pmatrix} \Id & 0 & 0 \\ 0 & \Id & 0 \\ 0 & 0 & \frac{\omega^{-1}}{2} \mathfig{0.1}{reidemeister_maps/R3_complex/delooping_handle} \\ 0 & 0 & \mathfig{0.1}{reidemeister_maps/R3_complex/delooping} \end{pmatrix} &
 \zeta_{l3}      & = \mathfig{0.1}{cobordisms/fringe_double_semicircle_down} \displaybreak[1] \\
 \zeta_{l2}^{-1} & = \begin{pmatrix} \Id & 0 & 0 & 0\\ 0 & \Id & 0 & 0 \\ 0 & 0 & \mathfig{0.1}{reidemeister_maps/R3_complex/looping} & \frac{\omega^{-1}}{2} \mathfig{0.1}{reidemeister_maps/R3_complex/looping_handle} \end{pmatrix} &
 \zeta_{l3}^{-1} & = \omega^2 \mathfig{0.1}{cobordisms/fringe_double_semicircle}
\end{align*}

We obtain the complex
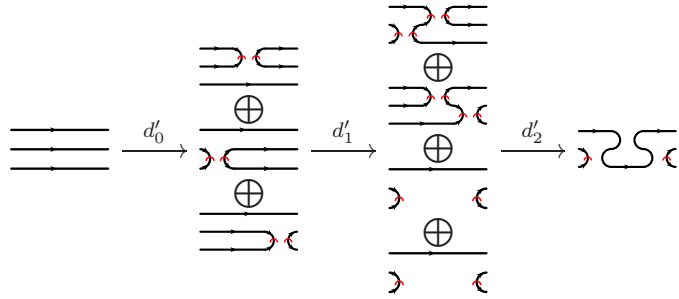
\newcommand{\directSumStackFour}[4]{{\begin{matrix}#1 \\ \DirectSum \\#2 \\ \DirectSum \\#3 \\ \DirectSum \\#4\end{matrix}}}
\begin{equation*}
 \xymatrix{
    \Ral{000} \ar[r]^{d_0'} &
    \directSumStackThree{\Ral{100}}{\Ral{010}}{\Ral{001}} \ar[r]^{d_1'} &
    \directSumStackFour{\Ral{110}}{\Ral{101}}{\Ral{011_delooped}}{\Ral{011_delooped}} \ar[r]^{d_2'} &
    \Ral{111_oriented}
 }
\end{equation*}
with differentials
\begin{equation*}
 d_0' = d_0 = \begin{pmatrix}s_1 \\ s_2 \\ s_3 \end{pmatrix}
\end{equation*}
\begin{align*}
 \begin{split}
  d_1' = \zeta_{l2} d_1 = \begin{pmatrix} \Id & 0 & 0 \\ 0 & \Id & 0 \\ 0 & 0 & \frac{\omega^{-1}}{2} \mathfig{0.1}{reidemeister_maps/R3_complex/delooping_handle} \\ 0 & 0 & \mathfig{0.1}{reidemeister_maps/R3_complex/delooping} \end{pmatrix} \begin{pmatrix}s_2 & -s_1 & 0 \\ s_3 & 0 & -s_1 \\ 0 & s_3 & -s_2\end{pmatrix} = \\ = \begin{pmatrix} \delta = \begin{pmatrix} s_2 & -s_1 \\ s_3 & 0 \end{pmatrix} & \gamma = \begin{pmatrix} 0 \\ -s_1\end{pmatrix} \\ \bullet & \bullet \\ \beta = \begin{pmatrix} 0 & \Id \end{pmatrix} & \psi = -\Id \end{pmatrix}\\
 \end{split} \displaybreak[1] \\
 \begin{split}
  d_2' = \zeta_{l3} d_2 \zeta_{l2}^{-1} = \mathfig{0.1}{cobordisms/fringe_double_semicircle_down} \begin{pmatrix}s_3 & -s_2 & s_1\end{pmatrix} \begin{pmatrix} \Id & 0 & 0 & 0\\ 0 & \Id & 0 & 0 \\ 0 & 0 & \mathfig{0.1}{reidemeister_maps/R3_complex/looping} & \frac{\omega^{-1}}{2} \mathfig{0.1}{reidemeister_maps/R3_complex/looping_handle} \end{pmatrix} = \\ = \begin{pmatrix} \lambda = \begin{pmatrix} \mathfig{0.05}{cobordisms/fringe_double_semicircle_down} s_3 & - \mathfig{0.05}{cobordisms/fringe_double_semicircle_down} s_2 \end{pmatrix} & \varphi = \omega^2 \Id & \bullet \end{pmatrix}.
 \end{split}
\end{align*}

Applying the double Gaussian elimination lemma, we reach the
homotopy equivalent complex
\begin{align}
\label{eq:R3_simplified_complex}%
\xymatrix{
   \Ral{000} \ar[r]^{d_0''} &
    \directSumStack{\Ral{100}}{\Rar{010}} \ar[r]^{d_1''} &
    \directSumStack{\Ral{110}}{\Ral{101}}
}
\end{align}
where
\begin{align}
\label{eq:R3_simplified_complex_differentials}%
d_0'' & = d_0' = \begin{pmatrix}\mathfig{0.1}{reidemeister_maps/R3_complex/simplified_d_0} \\ \\ \mathfig{0.1}{reidemeister_maps/R3_complex/simplified_d_1} \end{pmatrix} \\
d_1'' & = \delta - \gamma \psi^{-1} \beta = \begin{pmatrix}
\mathfig{0.1}{reidemeister_maps/R3_complex/simplified_d_00} &
-\mathfig{0.1}{reidemeister_maps/R3_complex/simplified_d_01} \\ & \\
\mathfig{0.1}{reidemeister_maps/R3_complex/simplified_d_10} &
-\mathfig{0.1}{reidemeister_maps/R3_complex/simplified_d_11}
\end{pmatrix}.
\end{align}
via the simplifying (and unsimplifying) maps
\begin{align*}
 s_{l0} & = \Id & s_{l0}^{-1} & = \Id \\
 s_{l1}   & = \begin{pmatrix} 1 & 0 & 0 \\ 0 & 1 & 0 \end{pmatrix} &
 s_{l1}^{-1} & = \begin{pmatrix} \begin{pmatrix} 1 & 0 \\ 0 & 1 \end{pmatrix} \\ -\psi^{-1} \beta \end{pmatrix} = \begin{pmatrix} 1 & 0 \\ 0 & 1 \\ 0 & 1 \end{pmatrix} \displaybreak[1] \\
 s_{l2}      & = \begin{pmatrix} \begin{pmatrix} 1 & 0 \\ 0 & 1 \end{pmatrix} & \begin{pmatrix} 0 \\ 0 \end{pmatrix} & -\gamma \psi^{-1} \end{pmatrix} \zeta_{l2} &
 s_{l2}^{-1} & = \zeta_{l2}^{-1} \begin{pmatrix} \begin{pmatrix} 1 & 0 \\ 0 & 1 \end{pmatrix} \\ -\varphi^{-1} \lambda \\ \begin{pmatrix} 0 & 0 \end{pmatrix} \end{pmatrix} \displaybreak[1] \\
             & = \begin{pmatrix} 1 & 0 & 0 \\ 0 & 1 & - c_1 \end{pmatrix}
             & & = \begin{pmatrix} 1 & 0  \\ 0 & 1 \\ -\omega^{-1} c_2 & \omega^{-1} c_3 \end{pmatrix} \\
 s_{l3} & = 0 & s_{l3}^{-1} & = 0.
\end{align*}
Here $c_1$ is the cobordism from $\Ral{011}$ to $\Ral{101}$ with three
components, a disc, a curtain, and a saddle, $c_2$ the similar cobordism from
$\Ral{110}$ to $\Ral{011}$ and $c_3$ is the similar cobordism from
$\Ral{101}$ to $\Ral{011}$ (the adjoint of $c_1$).

That's half the work! Now we need to do the same for the right side
of the third Reidemeister move, then compose a `simplifying map'
with an `unsimplifying map'.

Briefly, we calculate that the complex for the right side is
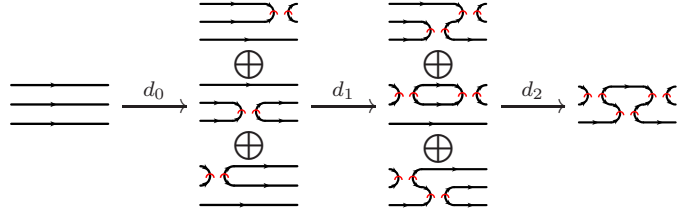
\begin{equation*}
 \xymatrix{
    \Rar{000} \ar[r]^{d_0} &
    \directSumStackThree{\Rar{100}}{\Rar{010}}{\Rar{001}} \ar[r]^{d_1} &
    \directSumStackThree{\Rar{110}}{\Rar{101}}{\Rar{011}} \ar[r]^{d_2} &
    \Rar{111}
 }
\end{equation*}
with differentials
\begin{align*}
 d_0 & = \begin{pmatrix}s_1 \\ s_2 \\ s_3 \end{pmatrix} \\
 d_1 & = \begin{pmatrix}s_2 & -s_1 & 0 \\ s_3 & 0 & -s_1 \\ 0 & s_3 & -s_2\end{pmatrix} \\
 d_2 & = \begin{pmatrix}s_3 & -s_2 & s_1\end{pmatrix}.
\end{align*}
and, applying the simplification algorithm, that this is homotopy
equivalent to the same complex as we obtained simplifying the other
side of the Reidemeister move (shown in Equation
\ref{eq:R3_simplified_complex}), but, somewhat tediously, with
slightly different differentials
\begin{align*}
 d_0'' & = \begin{pmatrix}\mathfig{0.1}{reidemeister_maps/R3_complex/simplified_d_0} \\ \\ \mathfig{0.1}{reidemeister_maps/R3_complex/simplified_d_1} \end{pmatrix} \\
 d_1'' & = \begin{pmatrix}
  \mathfig{0.1}{reidemeister_maps/R3_complex/simplified_d_00} &
  -\mathfig{0.1}{reidemeister_maps/R3_complex/simplified_d_01} \\ & \\
  -\mathfig{0.1}{reidemeister_maps/R3_complex/simplified_d_10} &
  \mathfig{0.1}{reidemeister_maps/R3_complex/simplified_d_11}
 \end{pmatrix}.
\end{align*}
These complexes thus differ by
\begin{align*}
 \xi_0 & = 1 &
 \xi_1 & = \begin{pmatrix} 1 &  0 \\ 0 & 1 \end{pmatrix} &
 \xi_2 & = \begin{pmatrix} 1 &  0 \\ 0 & -1 \end{pmatrix}.
\end{align*}

The simplifying and unsimplifying maps are
\begin{align*}
s_{r0} & = \Id & s_{r0}^{-1} & = \Id \\
s_{r1}      & = \begin{pmatrix} 1 & 0 & 0 \\ 0 & 1 & 0 \end{pmatrix}
&
s_{r1}^{-1} & = \begin{pmatrix} 1 & 0 \\ 0 & 1 \\ 1 & 0 \end{pmatrix} \displaybreak[1] \\
s_{r2}      & = \begin{pmatrix} 1 & 0 & 0 \\ 0 & -c_4 & 1
\end{pmatrix} &
s_{r2}^{-1} & = \begin{pmatrix} 1 & 0  \\ \omega^2 c_5 & \omega^2 c_6 \\ 0 & 1 \end{pmatrix} \\
s_{r3} & = 0 & s_{r3}^{-1} & = 0.
\end{align*}
Here
\begin{align*}
c_4 & :\Rar{101} \To \Rar{011} \\
c_5 & :\Rar{110} \To \Rar{101} \\
\intertext{and}%
c_6 & :\Rar{011} \To \Rar{101}
\end{align*}
are the obvious variations on $c_1, c_2$ and $c_3$.

The interesting compositions, which provide us with the chain map
between the two sides of the Reidemeister move, are
\begin{align*}
s_{r0}^{-1} \compose s_{l0} & = \begin{pmatrix} 1 \end{pmatrix} \\
s_{r1}^{-1} \compose s_{l1} & = \begin{pmatrix} 1 & 0 \\ 0 & 1 \\ 1 & 0 \end{pmatrix} \compose \begin{pmatrix} 1 & 0 & 0 \\ 0 & 1 & 0 \end{pmatrix} \\
                            & = \begin{pmatrix} 1 & 0 & 0 \\ 0 & 1 & 0 \\ 1 & 0 & 0 \end{pmatrix} \displaybreak[1] \\
s_{r2}^{-1} \compose \xi_2 \compose s_{l2} & = \begin{pmatrix} 1 & 0  \\ \omega^2 c_5 & \omega^2 c_6 \\ 0 & 1 \end{pmatrix} \compose \begin{pmatrix} 1 &  0 \\ 0 & -1 \end{pmatrix} \compose \begin{pmatrix} 1 & 0 & 0 \\ 0 & 1 & - c_1 \end{pmatrix} \\
                                           & = \begin{pmatrix} 1 & 0 & 0 \\ \omega^2 c_5 & -\omega^2 c_6 & \omega^2 c_6 c_1 \\ 0 & -1 & c_1 \end{pmatrix}\\
\end{align*}

The cobordism $c_6 c_1$ is the same `monkey saddle' appearing in
\cite{MR2174270}.

The maps described in Proposition \ref{prop:R3} describing the R3 chain map are simply a rearrangement of those presented here via matrices.
\end{proof}

\subsection{Proofs of the R3 variations lemmas}
\label{ssec:R3_variations}%
We now turn to the proofs of Lemmas \ref{lem:first-R3-general}, \ref{lem:second-R3-general} and \ref{lem:third-R3-general}.
As explained previously, in \S \ref{sec:R3-maps}, our strategy is to use the fact the Lemmas \ref{lem:first-R3}, \ref{lem:second-R3} and \ref{lem:third-R3}
are exactly the special case that the Reidemeister 3 move is $R3_{hml}$, and then to show that if two Reidemeister 3 variations are adjacent in the
cube of variations shown in Figure \ref{fig:cube-of-R3s}, and the spanning tree of definitions includes the connecting edge, then if the Lemmas hold for one
variation, they must hold for the other. However, this approach immediately requires two cases, depending on whether we are looking at
one of the four 'vertical' edges of the cube of R3 variations, or one of the eight 'horizontal edges'.

The case that the connecting edge is 'vertical', the formula defining one R3 variation in terms of the other (look back at Equation
\eqref{eq:edge-example}, for example, relating $R3_\circlearrowleft$ and $R3_{hlm}$) involves conjugation by an R2 move in the direction opposite the highest crossing.
On the other hand, when we look at the eight 'horizontal' edges, the R2 move conjugation takes place opposite either the middle or lowest crossing.
Because all of our lemmas are written describing the R3 moves in terms of resolutions of the highest crossing, it's unsurprising we need to treat these
cases separately.

It turns out that in order to prove Lemma \ref{lem:second-R3-general} for a given vertex $\star$, connected in the spanning tree to a vertex $\star'$,
we'll have to know slightly more about the $R3$ map for $\star'$ than is explicit in the Lemmas. This extra information follows from the Lemmas however, and so we'll
state it in the Corollary below. Once we have established the Lemmas for the vertex $\star'$ (starting at $\star'=hml$), we also know the Corollary
for $\star'$, and can use it in proving the Lemmas for $\star$.

\newcommand{\psuedoid}{\text{`$\Id$'}}

\begin{cor}[Corollary of Lemmas \ref{lem:first-R3-general} and \ref{lem:second-R3-general}]
\label{cor:R3-corollary}%
In the $\mathcal{P}$ layers of the cube of resolutions of $R3_\star$
there is exactly one resolution which appears for both the initial and final tangles of the R3 move. By grading considerations, the component of the
R3 map between these resolutions is some multiple of a map whose underlying unoriented surface is the identity.
There is always a unique configuration of seams on this surface without loops, and we will write $\psuedoid$ for such a `disorientation cylinder with minimal seams'.
Write $p_\star$ for the coefficient of this disoriented surface. Writing $\kappa_\star$ for the coefficient appearing
in the lowest homological height of the $\OO$ map, and $\lambda_\star$ for the coefficient in the highest height (so by Lemma \ref{lem:second-R3-general},
$\lambda_\star = \kappa_\star$ if $\star = hml, lhm, mhl$ or $lmh$, and $\lambda_\star = -\omega^2 \kappa_\star$ otherwise), we have

\begin{align*}
\frac{p_\star}{\kappa_\star} & =
\begin{cases}
-1          & \text{if $\star = hml, lmh, mlh$ or $\clockwise$} \\
\omega^2    & \text{if $\star = lhm, mhl, hlm$ or $\counterclockwise$}
\end{cases}
\end{align*}
\end{cor}

\begin{proof}
This is actually quite involved! Along the way, we'll also need to understand one of the coefficients in the 'downhill' map. We'll introduce some further
notation for particular resolutions of the R3 tangle, as follows: a symbol $abc$, with each of $a$ and $b$ either $>$ or $<$, and $c$ either
$\mathcal{O}$ or $\mathcal{P}$, refers to the resolution in which the first crossing is either in the higher or lower homological height
resolution, depending on $a$, the second crossing is again either in the higher or lower homological height
resolution, depending on $b$, and the third crossing is either in the orthogonal or parallel resolution relative to the triangle formed by the R3 tangle, depending on $c$.
Remember that the convention for the ordering of crossings is unobvious; before an R3 move (when the triangle is on the left of the lowest strand), the
crossings are ordered as 'middle' then 'low' then 'high', while after the R3 move the crossings are ordered as 'low' then 'middle' then 'high'.

It's easy to verify that in the $\OP$ cases, the ${><}\mathcal{O}$ resolution of the initial tangle is the same, ignoring orientation data, as the ${<<}\mathcal{P}$
resolution of the final tangle. Since these resolutions are in the same $q$-grading, the only maps between them are disoriented cylinders.
Taking into account orientation data, we claim that there is a unique allowed configuration of disorientation seams
on the cylinder with the appropriate boundary. We then define $q_\star$ to be the coefficient appearing on this map in the $R3_\star$ map.
Similarly, in the $\PO$ cases, the ${>>}\mathcal{P}$ resolution of the initial tangle is the same up to orientation data as the ${><}\mathcal{O}$ resolution
of the final tangle, and we define $q_\star$ to be the coefficient appearing on coefficient of the component of the R3 map between these resolutions.

The resolutions described in the statement of the corollary are in this notation ${><}\mathcal{P}$ (in the initial tangle) and ${<>}\mathcal{P}$ (in the final tangle), in the $\OP$ cases,
and the reverse in $\PO$ cases.

We now determine $q_\star$, and then $p_\star$ in terms of $\kappa_\star$, by considering the following two pairs of commuting squares coming from chain map conditions.

\begin{equation} 
\label{eq:OP-chain-map-condition}%
\xymatrix{
    {<<}\mathcal{O} \ar[r]^{\kappa_\star \Id} \ar[d]_{(-1)^{\star = lhm \text{ or } mlh}s} & {<<}\mathcal{O} \ar[d]^{+s} \\
    {><}\mathcal{O} \ar[r]^{q_\star \psuedoid} \ar[d]_{+s}                                           & {<<}\mathcal{P} \ar[d]^{(-1)^{\star=hml \text{ or } lhm}s} \\
    {><}\mathcal{P} \ar[r]^{p_\star \psuedoid}                                                       & {<>}\mathcal{P}
}
\end{equation}

\begin{equation} 
\label{eq:PO-chain-map-condition}%
\xymatrix{
    {>>}\mathcal{O} \ar[r]^{\lambda_\star \Id}                                                   & {>>}\mathcal{O} \\
    {>>}\mathcal{P} \ar[r]^{q_\star \psuedoid} \ar[u]^{+s}                                           & {><}\mathcal{O} \ar[u]_{(-1)^{\star=hlm \text{ or } \clockwise}s} \\
    {<>}\mathcal{P} \ar[r]^{p_\star \psuedoid} \ar[u]^{(-1)^{\star = \clockwise \text{ or } lmh}s}   & {><}\mathcal{P} \ar[u]_{+s}
}
\end{equation}

The signs appearing on saddles in Equations \eqref{eq:OP-chain-map-condition} and \eqref{eq:PO-chain-map-condition} are calculated by the usual rule
for sprinkling signs in tensor products of complexes (see \S \ref{ssec:appendix:tensor}), and the convention for ordering crossings before and after
Reidemeister 3 moves.

These calculations tell us that
\begin{align*}
p_\star & =
\begin{cases}
-\kappa_\star         & \text{if $\star = hml, mlh, lmh$ or $\clockwise$} \\
\omega^2 \kappa_\star & \text{if $\star = hlm, lhm, mhl$ or $\counterclockwise$}
\end{cases} \\
        & =
\begin{cases}
-1        & \text{if $\star = hml$ or $lmh$} \\
1         & \text{if $\star = hlm$ or $mlh$} \\
\omega^2  & \text{if $\star = lhm$ or $mhl$} \\
-\omega^2 & \text{if $\star = \clockwise$ or $\counterclockwise$.}
\end{cases}
\end{align*}

Unfortunately there's a small subtlety in extracting the relation between $k_\star$ and $p_\star$ from Equations \eqref{eq:OP-chain-map-condition} and \eqref{eq:PO-chain-map-condition};
the horizontal arrows are not labelled by multiples of the identity map, but by multiples of the `minimal seam' identity map. Depending on
the configuration of these seams, it might not be the case that $s \psuedoid = \psuedoid s$, but that they differ by a power of $\omega$.
This requires a case by case analysis. Defining $\sigma_\star$ and $\tau_\star$ so that $s \psuedoid = \sigma_\star \psuedoid s$ in the upper squares of
Equations \eqref{eq:OP-chain-map-condition} and \eqref{eq:PO-chain-map-condition}, and $s \psuedoid = \tau_\star \psuedoid s$ in the lower squares, we find
that all $\sigma_\star$ and $\tau_\star$ are equal to $1$, except that
\begin{equation*}
\sigma_{lhm} = \sigma_{mhl} = \tau_\clockwise = \tau_\counterclockwise = \omega^2
\end{equation*}

The corollary now follows.
\end{proof}

Let's begin the 'vertical' edge case by introducing some notation
for particular subspaces of the complexes associated to the four tangles appearing in our formula for one $R3$ move in terms of another.
The symbols $\mathcal{O}|$ and $\mathcal{P}|$ will denote the spaces of the Khovanov complex of the initial tangle in which the highest crossing has been resolved in the orthogonal and parallel
manners respectively. The symbols $|\mathcal{O}$ and $|\mathcal{P}$ will denote the corresponding subspaces of the final tangle.
The symbols $a|bc$ and $ab|c$, where $a,b,c = \mathcal{O}$ or $\mathcal{P}$ will denote subspaces of the two intermediate tangles (the lower left and lower right tangles in Equation \eqref{eq:edge-example}, respectively), in which the three crossings
not involving the lowest strand (that, the original highest crossing, and the two new crossings introduced by the $R2$ move) have been resolved
either orthogonal or parallel to the lowest strand, according to the values of $a, b$ and $c$, with $a$ referring to the original highest crossing, $b$
referring to the new crossing closest to the original tangle, and $c$ to the new crossing furthest away. (The vertical bar $|$ is meant
to denote the lowest strand.)

\begin{figure}[!ht]
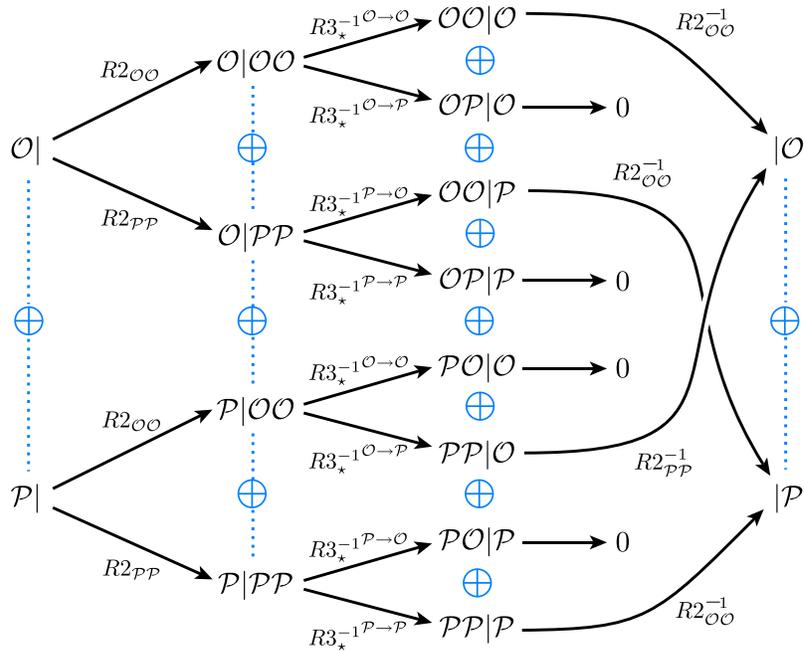

\begin{equation*}
\mathfig{0.8}{reidemeister_maps/R3_variations/lemma}
\end{equation*}
\caption{One R3 variation defined in terms of another via R2 moves.}
\label{fig:R3-variation-composition}
\end{figure}

We know, from Figures \ref{fig:R2a_maps} and \ref{fig:R2b_maps}, that the $R2$ maps only see those subspaces in which the two crossings involved in the $R2$ move have been resolved the same way, that is, with $b=c$.
Thus if $R3_{\star'}$ is defined in terms of $R3_\star$ as the composition of an $R2$ map, the map $R3_\star^{-1}$, and an inverse $R2$ map, as in the example in Equation \eqref{eq:edge-example},
it will have the form shown in Figure \ref{fig:R3-variation-composition}.

\begin{proof}[Proof of Lemma \ref{lem:first-R3-general}, vertical edge cases.]
If the statement of the Lemma holds for some move $R3_\star$, and we're defining
another $R3_{\star'}$ in terms of it via a vertical edge of the cube in Figure
\ref{fig:cube-of-R3s}, it must also hold for $R3_\star'$. (Recall that adjacent $R3$ moves in the cube have opposite arrangements of layers.)
This is the case simply because in Figure \ref{fig:R3-variation-composition} the map $R3_{\star'}^\PO$ factors through ${R3_{\star}^{-1}}^\OP$, and the map
$R3_{\star'}^\OP$ factors through ${R3_{\star}^{-1}}^\PO$.
\end{proof}

\begin{proof}[Proof of Lemma \ref{lem:second-R3-general}, vertical edge cases.]
Using Lemma \ref{lem:first-R3-general}, one sees that the component
$R3_\star^{\OO}$ is itself a chain map. This follows in the case that the layers are arranged as $\OP$ by writing
$d = d_\OO + d_\OP + d_\PP$; then the $\OO$ component of the equation $d R3_\star = R3_\star d$ simply says
$$d_\OO R3_\star^\OO = R3_\star^\OO d_\OO + R3_\star^\OP d_\PO$$
and since by Lemma \ref{lem:first-R3-general} $R3_\star^\PO = 0$,
$$d_\OO R3_\star^\OO = R3_\star^\OO d_\OO.$$
The other case, in which the layers are arranged as $\PO$, is the same.

We claim then that the only
chain maps from one orthogonal layer to another are multiples of the
identity in the $\star = hml, lhm, mhl$ or $lmh$ cases, and multiples of
the standard chain map described in the Lemma when $\star = hlm,
mlh, \circlearrowleft$ or $\circlearrowright$. The overall coefficient
is easily determined from Figure \ref{fig:R3-variation-composition}, and the formulas for the R2 maps in Figures \ref{fig:R2a_maps} and \ref{fig:R2b_maps}.
We obtain the result described in the Lemma, that for $\star = hml, lhm, mhl$ or $lmh$ the component in lowest homological height is actually the identity,
that for $\star = hlm$ and $\clockwise$ the coefficient of the component in the lowest homological height is $\omega^2$, and that for
$\star = mlh$ and $\counterclockwise$ that coefficient is $-1$.
\end{proof}

\begin{proof}[Proof of Lemma \ref{lem:third-R3-general}, vertical edge cases.]
Again looking at Figure \ref{fig:R3-variation-composition}, we see that $R3_{\star'}^\PP$ factors through ${R3_{\star}^{-1}}^\PP$. Thus
if the lemma holds for $R3_\star$ (which it does for $\star = hml$, by Lemma \ref{lem:third-R3}), it also holds for any adjacent $R3$ move which
we're defining in terms $R3_\star$. The second part of the Lemma, describing the normalisation, has already been proved as part of Corollary \ref{cor:R3-corollary}.
The final statement, about the other entries of the map in the middle homological height having disc components, follows immediately from grading considerations.
\end{proof}

We now deal with the cases involving a horizontal edge.

\begin{proof}[Proof of Lemma \ref{lem:first-R3-general}, horizontal edge cases.]
We'll introduce a new notational convention; when decorating a crossing with an $\mathcal{O}$ or a $\mathcal{P}$, to indicate a particular resolution,
we'll also draw a short squiggly line pointing towards the nearby strand with respect to which we mean `orthogonal' or `parallel'.
We can rotate this short squiggly line into a different region adjacent to the crossing, if at the same time we interchange the labels $\mathcal{O}$
and $\mathcal{P}$.

When we define an R3 map via a `horizontal' edge in the cube, in terms of some other map, say $R3_\star$, it has the form:
\begin{equation}
\label{eq:R3-horizontal-composition1}%
\mathfig{0.9}{reidemeister_maps/R3_variations/lemma_horizontal}.
\end{equation}
(Here the labelled crossings in the initial and final step are the highest crossings, as usual.)
Thus we see that the $\OP$ component factors through the original $\PO$ component of $R3_\star$, and similarly the $\PO$ component
factors through the $\OP$ component. Since the relative heights of the $\mathcal{O}$ and $\mathcal{P}$ maps are reversed in adjacent
R3 variations in the cube, this suffices to establish the lemma.
\end{proof}

\begin{proof}[Proof of Lemma \ref{lem:second-R3-general}, horizontal edge cases.]
This follows the argument above in the vertical edge case; Lemma \ref{lem:first-R3-general}, ensures that the
component $R3_\star^{\OO}$ is a chain map, and thus only multiples of the map described in this Lemma are possible. To check that the
multiple is the one described, we follow through the $\OO$ composition in Equation \eqref{eq:R3-horizontal-composition1} above.
Notice that this relies on Corollary \ref{cor:R3-corollary}, for the normalisation of the $\PP$ map appearing in Equation \eqref{eq:R3-horizontal-composition1}.
Further, in the cases where the R2 moves appearing are R2b moves, one must take into account a sign of homological origin, coming from reordering crossings.
\end{proof}

\begin{proof}[Proof of Lemma \ref{lem:third-R3-general}, horizontal edge cases.]
We now look in slightly more detail at Equation \eqref{eq:R3-horizontal-composition1}. The highest and lowest homological heights of the $\mathcal{P}$
layer consist of those resolutions in which the other two crossings (i.e., the middle and lowest crossing) have been resolved in opposite ways; one
as a $\mathcal{O}$, one as a $\mathcal{P}$. We look at one of the two cases, the other being essentially identical.
Making use of Lemma \ref{lem:second-R3-general} (in particular, that, ignoring all disorientation data and coefficients, the $\OO$ components
of all R3 variations are simply the identity), we see
\begin{equation}
\label{eq:R3-horizontal-composition3}%
\mathfig{0.95}{reidemeister_maps/R3_variations/lemma_horizontal2}.
\end{equation}
The second part of the Lemma, describing the normalisation, has already been proved as part of Corollary \ref{cor:R3-corollary}.
The final statement, about the other entries of the map in the middle homological height having disc components, follows immediately from grading considerations.
\end{proof}

\subsection{Planar algebras and canopolises}
\label{ssec:planar-algebras}%
A planar algebra is a gadget specifying how to combine objects in
planar ways. They were introduced in \cite{math.QA/9909027} to study
subfactors, and have since found more general use.

\newcommand{\pa}{\mathcal{P}}
In the simplest version, a planar algebra $\pa$ associates a vector space
$\pa_k$ to each natural number $k$ (thought of as a disc in the plane with
$k$ marked points on its boundary) and a linear map $\pa(T)  : \pa_{k_1}
\tensor \pa_{k_2} \tensor \cdots \tensor \pa_{k_r} \To \pa_{k_0}$ to each
planar tangle\footnote{Familiarly known as a 'spaghetti and meatballs'
diagram.} $T$, for example
$$\mathfig{0.2}{planar_algebras/unoriented_cubic_tangle_example},$$ with internal
discs with $k_1, k_2, \ldots, k_r$ marked points, and $k_0$ marked points
on the external disc. These maps (the 'planar operations') must satisfy
certain properties: ``radial'' tangles induce identity maps, and
composition of the maps $\pa(T)$ is compatible with the obvious
composition of planar diagrams by gluing one inside the other.

For the exact details, which are somewhat technical,
see \cite{math.QA/9909027}.

{
\newcommand{\labels}{\mathfrak{L}}
Planar algebras also come in more subtle flavors. Firstly, we can
introduce a label set $\labels$, and associate a vector space to
each disc with boundary points marked by this label set. (The
simplest version discussed above thus has a singleton label set, and
the discs are indexed by the number of boundary points.) The planar
tangles must now have arcs labeled using the label set, and the
rules for composition of diagrams require that labels match up.
Secondly, we needn't have vector spaces and linear maps between
them; a planar algebra can be defined over an arbitrary monoidal
category, associating objects to discs, and morphisms to planar
tangles. Thus we might say ``$\mathcal{P}$ is a planar algebra over
the category $\C$ with label set $\labels$.'' \footnote{A
``subfactor planar algebra'' is defined over $\operatorname{Vect}$,
and has a 2 element label set. We impose an additional condition
that only discs with an even number of boundary points and with
alternating labels have non-trivial vector spaces attached. There is
also a positivity
condition. See \cite[\S4]{math.GT/0606318}.}%
}

A ``canopolis'', introduced by Bar-Natan in
\cite{MR2174270}\footnote{He called it a `canopoly', instead, but
we're taking the liberty of fixing the name.}\footnote{See also \cite{math.GT/0610650} for a description of
Khovanov-Rozansky homology \cite{math.QA/0401268, math.QA/0505056} using canopolises.}, is simply a
planar algebra defined over some category of categories, with
monoidal structure given by cartesian product. Thus to each disc, we
associate some category of a specified type. A planar tangle then
induces a functor from the product of internal disc categories to
the outer disc category, thus taking a tuple of internal disc
objects to an external disc object, and a tuple of internal disc
morphisms to an external disc morphism. It is picturesque to think
of the objects living on discs, and the morphisms in `cans', whose
bottom and top surfaces correspond to the source and target objects.
Composition of morphisms is achieved by stacking cans vertically,
and the planar operations put cans side by side.

The functoriality of the planar algebra operations ensure that we
can build a `city of cans' (hence the name canopolis) any way we
like, obtaining the same result: either constructing several `towers
of cans' by composing morphisms, then combining them horizontally,
or constructing each layer by combining the levels of all the towers
using the planar operations, and then stacking the levels
vertically.

\subsection{Complexes in a canopolis form a planar algebra}
\label{ssec:appendix:derived_canopolis}%
{ 
\newsavebox{\qtbox}
\sbox{\qtbox}{\ensuremath{\mathfig{0.08}{complexes/quadratic_tangle}}}
\newcommand{\qt}{\usebox{\qtbox}}
\newcommand{\dt}{\mathfig{0.05}{complexes/disc_3}}
\newcommand{\df}{\mathfig{0.05}{complexes/disc_5}}

Given a quadratic tangle,
$\mathfig{0.16}{complexes/quadratic_tangle}$ and a pair of complexes
associated to the inner discs,
\begin{equation*}
\xymatrix@R-25pt{
 C_{\color{red}1}   = \Big( \dt \ar@[red][r]  & \dt \ar@[red][r]  & \dt \ar@[red][r]  & \dt \Big) \\
 C_{\color{blue}2}  = \Big( \df \ar@[blue][r] & \df \ar@[blue][r] & \df \ar@[blue][r] & \df \Big)\\
}
\end{equation*}
we need to define a new complex associated to the outer disc.

We'll imitate the usual construction for tensor product of
complexes, but use the quadratic tangle to combine objects and
morphisms. Form a double complex then collapse along the
anti-diagonal:
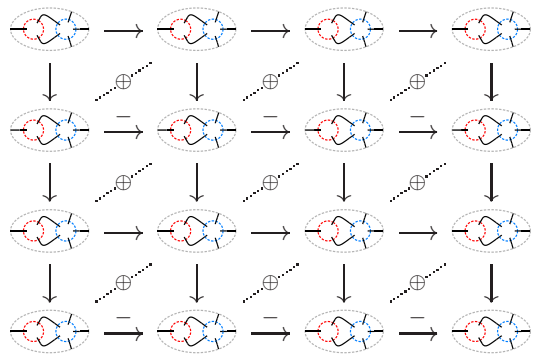
\begin{equation*}
\xymatrix@R-10pt@C-10pt{
 \qt \ar@[red][r]     \ar@[blue][d] &  \qt \ar@[red][r]     \ar@[blue][d] \ar@{.}[dl]|{\directSum} & \qt \ar@[red][r]     \ar@[blue][d] \ar@{.}[dl]|{\directSum} & \qt \ar@[blue][d] \ar@{.}[dl]|{\directSum} \\
 \qt \ar@[red][r]^{-} \ar@[blue][d] &  \qt \ar@[red][r]^{-} \ar@[blue][d] \ar@{.}[dl]|{\directSum} & \qt \ar@[red][r]^{-} \ar@[blue][d] \ar@{.}[dl]|{\directSum} & \qt \ar@[blue][d] \ar@{.}[dl]|{\directSum} \\
 \qt \ar@[red][r]     \ar@[blue][d] &  \qt \ar@[red][r]     \ar@[blue][d] \ar@{.}[dl]|{\directSum} & \qt \ar@[red][r]     \ar@[blue][d] \ar@{.}[dl]|{\directSum} & \qt \ar@[blue][d] \ar@{.}[dl]|{\directSum} \\
 \qt \ar@[red][r]^{-}               &  \qt \ar@[red][r]^{-}                                        & \qt \ar@[red][r]^{-}                                        & \qt               \\
}
\end{equation*}
Here each horizontal arrow is the planar composition of a morphism from
$C_1$, placed in the left disc, with the identity on the appropriate
object from $C_2$, in the right disc. Similarly, each vertical arrow
is the planar composition of a morphism from $C_2$ with an identity
morphism.

The extension to tangles with more than 2 internal discs is obvious.
Moreover, it's not hard to see that chain maps between complexes in
a canopolis also form a planar algebra, providing the morphism part
of `the canopolis of complexes and chain maps'.
} 

\section{Homological conventions}
\label{sec:appendix:homological_conventions}%
\subsection{Tensor product}
\label{ssec:appendix:tensor}%
In the next two sections we'll describe certain conventions to do
with tensoring complexes. (Please accept our apologies if they're
not what you're used to!) \cite{MR1438306} 

The tensor product of two complexes $(A^\bullet, d_A)$ and
$(B^\bullet, d_B)$ is defined to be
$$(A \tensor B)^\bullet = \DirectSum_{i+j=\bullet} A^i \tensor B^j,$$
and $$d_{(A \tensor B)^\bullet} = \sum_{i+j=\bullet} (-1)^j d_{A^i}
\tensor \Id_{B^j} + \Id_{A^i} \tensor d_{B^j}.$$

If you think of $A^\bullet$ as lying horizontally, and $B^\bullet$
as vertically, this rule says ``negate the differentials in every
odd row''. 

\subsection{Permuting tensor products}
\label{ssec:appendix:permuting}%
Unfortunately, while $A^\bullet \tensor B^\bullet \Iso B^\bullet
\tensor A^\bullet$ the isomorphism can't just be the identity.
Instead, we'll take it to be $A^i \tensor B^j \mapsto (-1)^{ij} B^j
\tensor A^i$; that is it negates anything in `doubly odd' degree.

The only complexes we ever take tensor products of are the complexes
associated to tangles. In the simplest case, where we are taking the
tensor product of two crossings, the `crossing reordering' map is
`negate doubly disoriented smoothings'. That is, objects in which
both crossings have been resolved in the disoriented direction get
negated when we change the ordering of the crossings.

\newcommand{\urlprefix}{}
\bibliographystyle{gtart}
\bibliography{bibliography/bibliography}

This paper is available online at arXiv:\arxiv{math.GT/0701339}, and
at \url{http://tqft.net/functoriality}.


\end{document}